\documentclass{article}

\usepackage{amsthm}
\usepackage{amsfonts}
\usepackage{amssymb}
\usepackage{verbatim}
\usepackage[all]{xy}
\usepackage{amsmath, amscd, amsthm}
\usepackage{longtable}
\usepackage{amscd}
\usepackage{mathrsfs}
\usepackage{graphicx}
\usepackage{latexsym}

\newcommand{\la}{\langle}
\newcommand{\ra}{\rangle}
\newtheorem{theorem}{Theorem}[section]
\newtheorem{lemma}[theorem]{Lemma}

\newtheorem{defn}[theorem]{Definition}
\newtheorem{rem}[theorem]{Remark}

\newtheorem{Condition}[theorem]{Condition}

\newtheorem{coro}[theorem]{Corollary}

\usepackage{latexsym}

\DeclareMathOperator{\Lab}{Lab}

\begin{document}
\title{Subnormal subgroups in free groups, their growth and cogrowth}
\author{A. Yu. Olshanskii \thanks{The
author was supported in part by the NSF grant DMS 1161294 and by the Russian Fund for Basic Research  grant 11-01-00945}}
\maketitle

\begin{abstract} In this paper, the author (1) compares subnormal closures of finite sets in a free group $F$; (2) obtains the limit for the series of subnormal closures of a single element in $F$;
(3) proves that the exponential growth rate (e.g.r.)  $\lim_{n\to \infty}\sqrt[n]{g_H(n)}$, where $g_H(n)$ is
the growth function of a subgroup $H$ with respect to a finite free basis of $F$, exists for
any subgroup $H$ of the free group $F$; (4) gives sharp estimates from below for the e.g.r. of
subnormal subgroups in free groups; and (5) finds the cogrowth of  the subnormal closures of free generators.
\end{abstract}

{\bf Key words:} free group, subnormal subgroup, normal closure, growth, random walk, van Kampen diagram,
wreath product.

\medskip

{\bf AMS Mathematical Subject Classification:}
20E05, 20D35, 20F69, 20P05, 20E22, 05C81, 05C38.

%\large

\section{Introduction}

Recall that a subgroup $H$ of a group $G$ is called $\ell$-{\it subnormal}
($\ell = 1, 2,\dots$) if there is a decreasing sequence
$G\; \triangleright H_1\; \triangleright \dots \triangleright H_{\ell}=H,$
where each term is a normal subgroup in the preceding one. For a subset $S\subset G$,
its $\ell$-subnormal closure $\la S\ra^G_{\ell}$ is defined by induction: $\la S\ra^G_1=\la S\ra^G $ is the normal closure
of $S$ in $G$, i.e., the smallest normal subgroup of $G$ containing the set $S,$ and $\la S\ra^G_{\ell}$
is the normal closure of $S$ in $\la S\ra^G_{\ell-1}.$ Clearly, the subgroup  $\la S\ra^G_{\ell}$ is
$\ell$-subnormal, and the obvious induction on $\ell$
shows that an $\ell$-subnormal subgroup $H$ of a group $G$ contains the $\ell$-subnormal closure in $G$
of arbitrary subset $S\subset H.$ It follows that $\la S\ra^G_{\ell}$ is the intersection of all ${\ell}$-subnormal
subgroups of $G$ containing $S$.

The concepts of subnormal series and subnormal subgroups are among the introductory ones in Group Theory.
The material devoted to subnormal subgroups  is exposed in the book \cite{SL} and in numerous papers; however,
our present topic, the subnormality in free groups, is not covered there.  The problems of asymptotic behavior raised in our paper are also new for this area. We start with the following
feature of normal and subnormal closures in free groups.

 \begin{theorem} \label{shm}   Let  $H$ be a normal subgroup in a free group $F$
with infinite factor group $F/H,$ and $S$ a finite subset of $H$.
Then the normal closure $N=\la S\ra ^H$ of $S$ in $H$ contains no nontrivial
normal subgroups of $F$.
\end{theorem}

\begin{coro}\label{t0601} Let $\ell\ge 1$, $F$ be a free group, $H$ a normal subgroup of infinite
index in $F$, and $H(\ell)=\la S \ra_{\ell}^H$  be the  $\ell$-subnormal closure of a finite subset $S$ in $H$. Then $H(\ell)$ does not contain nontrivial $\ell$-subnormal subgroups of $F$.
\end{coro}

\begin{coro}\label{g}

The $(\ell+1)$-subnormal closure of an element $g$ in a noncyclic free group $F$
contains no nontrivial $\ell$-subnormal in $F$ subgroups provided $\ell\ge 1$.
In particular, it contains no nontrivial normal subgroups of $F$.
\end{coro}

The next theorem exhibits one more characteristic  of
subnormal closures of a {\it single} element in
free groups. It is easy to see that the direct generalization of Theorem \ref{intersec} (1) to the normal closures of larger  sets fails.

\begin{theorem} \label{intersec}
 (1) If $g$ is an element of a free group $F$,
then $\cap_{\ell=1}^{\infty} \langle g\rangle_{\ell}^F=\langle g\rangle.$

(2) Moreover, let $D(g,\ell)$ denote the minimum of the lengths of the
elements in $\langle g\rangle^F_{\ell}\backslash \langle g \rangle $. Then
there is $c>0$ such that $D(g,\ell)\ge c 2^{\ell}$ for every $\ell\ge 1$, i.e.,
the normal closures $\langle g\rangle^F_{\ell} $ converge to the cyclic
subgroup $\langle g\rangle$ exponentially fast.
\end{theorem}

Then we focus on growth and cogrowth rates of subnormal subgroups in free groups.
Let $F_m=F(X)$ be a free group with a  free basis $X=\{x_1,\dots,x_m\}$. The (relative) growth
function $g_U(n)$ ($n=0,1,\dots$) of arbitrary subset $U\subset F_m$ with respect to $X$ is given
by the formula $g_U(n)=\#\{g\in U\mid |g|\le n\}$, where the sign  $\#$ is placed for the number of
elements and $|g|=|g|_X$ denote the length of an element with respect to the basis $X$.

Clearly, the relative growth function of any non-cyclic subgroup $H$ in  $F$ is at least exponential
(since $H$ is a free group itself) and at most exponential (since $H\le F$). Can one assert that
for any subgroup $H\le F_m$, its growth is exponential in the stronger sense that there is a limit
$a=\lim_{n\to\infty}\sqrt[n]{g_H(n)}$ (i.e., for any $\varepsilon>0,$ the function $g_H(n)$ is
between $(a-\varepsilon)^n$ and  $(a+\varepsilon)^n$ for every large enough $n$ )? Probably this
question has been open (instead the symbol `$\limsup$' is used in literature). The affirmative answer for normal subgroups was given by R.Grigorchuk \cite{Gr}, and for finitely generated subgroups, it can be derived from \cite{Gr}. Nevertheless the limit exists for arbitrary subgroup $H\le F_m$:

\begin{theorem} \label{rate} For any subgroup $H$ of a finitely generated free group $F(X)$, there exists
the exponential growth rate (or e.g. rate) $\alpha_H = \lim_{n\to\infty}\sqrt[n]{g_H(n)}$, where $g_H(n)$
is the growth function of $H$ with respect to the free basis $X$ of $F(X)$.
\end{theorem}

Note that this statement fails if one replaces the free group $F(X)$ by arbitrary finitely generated group; see Remark \ref{fa}.

Theorem \ref{rate} follows from a general assertion about the number of reduced paths of given length
in arbitrary graph. (A combinatorial  path $p=e_1\dots e_k$ of length $k$ is called {\it reduced} if no edge
$e_i$ is the inverse one for the preceding  edge $e_{i-1}$, $i=2,\dots,k$.)

\begin{lemma} \label{graph} Let $\Gamma$ be  a graph with bounded degrees of vertices, $o$ a vertex of $\Gamma$,
and $\nu(n)$ the number of reduced paths of length at most $n$ starting and terminating at $o$.
Then there is a finite limit $\lim_{n\to\infty}\sqrt[n]{\nu(n)}.$

For a connected $\Gamma$, this limit does not depend on the choice of $o$.
\end{lemma}

We prove Lemma \ref{graph} in this paper since usually graph theorists count all paths (reduced or not) between two vertices using the adjacency matrix. (For example, see \cite{GM}, where some results on finite graphs are
extended to infinite graphs with bounded vertex degrees.) But we need to know the numbers of {\it reduced} paths.

Note that, on the one hand, the e.g. rate $\alpha_H$ does not exceed $2m-1$ for any subgroup
since the number of reduced words of length $n\ge 1$ over the alphabet $X^{\pm 1}$ is equal to $2m(2m-1)^{n-1}$.
On the other hand, if $H$ is a nontrivial normal subgroup, then $\alpha_H\ge \sqrt{2m-1}$. This is easy to see just
counting the number of distinct conjugates $vgv^{-1}$ of a non-trivial $g\in H,$ where $v\in F_m$ and $|v|<(n-|g|)/2$.
In fact, the inequality for $\alpha_H$ is strict:

\begin{lemma} \label{gri} (R.I.Grigorchuk \cite{Gr}). If $m\ge 2$, then for any nontrivial normal subgroup $H\le F_m$, we have $\alpha_H > \sqrt{2m-1}$.
\end{lemma}

The sharp estimate of growth exponents for subnormal subgroups of $F_m$ is given
by the following theorem, where the second part is harder.

\begin{theorem}\label{base} (1) Let $m\ge 2$ and $N$ be a non-trivial $\ell$-subnormal subgroup of $F_m$ for some $\ell\ge 1.$
Then the e.g. rate $\alpha_N$ of $N$ with respect to the free basis of $F_m$ is greater than $(2m-1)^{2^{-\ell}}.$

(2) For any $\varepsilon>0$, there is a nontrivial $\ell$-subnormal subgroup $N$ of $F_m$ with $\alpha_N < (2m-1)^{2^{-\ell}}+\varepsilon.$
\end{theorem}

The {\it cogrowth function} $f_H(n)$ of a subgroup $H\le F_m$ with respect to the basis $X$ counts the number of distinct cosets  of $H$ in $F_m$ with length at most $n$, i.e., $f_H(n)=\#\{Hg\mid |g|\le n\}.$ So the cogrowth
function is equal to the growth function of the transitive action of $F_m$ with the stabilizer of a distinguished
point equal to $H$.
%There is no simple connection of this concept to growth,
In  contrast with the situation in algebras, where the dimension growth of a subalgebra plus its cogrowth is just equal to the dimension growth of the whole algebra, the rate of growth of $H\le F_m$ does not determine the rate of its
cogrowth, nor vice versa.

If $H$ is a normal subgroup of $F_m$, then the cogrowth function $f_H(n)$ is equal to the growth function
of the factor group $F_m/H$ with respect to (the canonical image of) $X$. Therefore it is exponentially
negligible in comparison with the growth of $F_m$ if $H\ne \{1\}$, that is, for every $n$,
 $f_H(n)\le C(2m-1-\varepsilon)^n,$ where $C$ and $\varepsilon$ are some positive constants. This
 follows from the fact that non-empty reduced words $w\in H$ cannot be subwords of the shortest
 coset representatives for $H$ (see \cite{GH} and also \cite{S} or \cite{BO}). Do nontrivial subnormal subgroups
 of $F_m$ have an exponentially negligible cogrowth too?
 It turns out that the answer given by the next theorem is negative (in contrast with the cogrowth of subideals
 in free Lie algebra \cite{BO1}), and the behavior
 of the cogrowth functions of $\ell$-subnormal subgroups is different for $\ell>1.$

  We shall use the symbol $\Theta$ from Computational Complexity Theory for the following equivalence
     of two real-valued functions $\phi(n)$ and $\psi(n)$, namely, $\phi(n)$ is $\Theta$-equivalent to $\psi(n)$ if
     both conditions $\phi(n) =O(\psi(n))$ and $\psi(n)=O(\phi(n))$ hold. In other words, $\phi(n)=\Theta(\psi(n))$.
   As in \cite{BO}, we say that the cogrowth of a subgroup $H\le F_m$ is {\it maximal} if the cogrowth function
   $f_H(n)$ is  $\Theta$-equivalent to the cogrowth function of the trivial subgroup, that is to
   the growth function of the group $F_m,$ i.e., $f_H(n)=\Theta((2m-1)^n)$. For example, every finitely
   generated subgroup of infinite index in $F_m$ has maximal cogrowth. (See \cite{BO} for this and other examples
   of maximal cogrowth in $F_m$.)

 \begin{theorem} \label{cogr} (1) The cogrowth function $f_H(n)$ of the $\ell$-subnormal closure $\la x \ra_{\ell}^{F_m},$
 where $x\in X,$

 (1) is maximal if  $\ell \ge 2$, and $m \ge 3$.

 (2) It is also maximal if $\ell\ge 3$ and $m\ge 2$.

(3) If $\ell=m=2$, then $f_H (n) =\Theta\left(\frac{3^n}{\sqrt n}\right)$.

 \end{theorem}

 The proofs of the five theorems formulated above are located in different
 sections since they are based on different ideas and can be read (almost) independently.

 In Section \ref{AL},
 we apply nilpotent wreath products and Shmelkin's embedding construction to prove Theorem \ref{shm}.
 The verbal wreath products had been invented in \cite{Sh} for the exploration of the products of group varieties.
 Hereby we show that they are also helpful aside from the theory of group varieties.

 Lemma \ref{graph} and Theorem \ref{rate} are proved in Section \ref{eg}, where the estimates are based
 on the well-known choice of the free basis in the fundamental group of a graph.

 To explore the subgroups of the form $\la R\ra^H$, where $H$ is a subgroup of $F$ and $R\subset H$, we introduce the concept of $H$-diagrams and apply it
 to the proof of Theorem \ref{intersec} in Section \ref{dia}.
 The proof of Theorem \ref{base}  is placed in Section \ref{subgrowth}. Here we also need $H$-diagrams. However
  the main issue, in comparison with the case $H=F_m$, is
 that the maximal prefix  of two reduced subwords from $H$ (with respect to the generators of $F_m$)
 does not belong to $H$, and this disturb the control of cancelations in the products of the $H$-conjugates of the elements from $R^{\pm 1}$. The  condition $P(d,\mu, \rho)$ is introduced
 to carry out the induction on the index of subnormality $\ell$.

 Sections \ref{cos} - \ref{part2} occupy the larger half of the paper. Here we prove the statements of Theorem \ref{cogr}  using more laborious combinatorial and probabilistic estimates. In particular, the proof of the strong $\Theta$-asymptotics in Theorem \ref{cogr} (3) is obtained due to the lucky opportunity of encoding coset representatives into 2-dimensional random trajectories
 lying above the $x$-axis. Other difficulties and tricks in these sections are caused by the presence of correlation in the 'reduced' random walk on an integer lattice.

\section{Comparison of subnormal closures in free groups}\label{AL}

{\bf Proof of Theorem \ref{shm}.} Assume by the contrary, that
the subgroup $N$ contains a nontrivial subgroup $M$ normal in $F$.

Since $H$ is a free group, we have $\cap_{i=1}^{\infty}\gamma_i(H)=\{1\},$
where $\gamma_i(K)$ denotes the $i$'th term of the lower central series of a group $K$
(see \cite{MKS}, Theorem 5.7 and Corollary 5.7). Hence we can find $c\ge 1$ such that $M\le \gamma_c(H)$ but $M$ is not
contained in $\gamma_{c+1}(H).$ Therefore the canonical image $\tilde M$ of $M$ in $\tilde F=F/\gamma_{c+1}(H)$
is nontrivial and  lies in the center of the image $\tilde H$ of $H$ in $\tilde F$.  Also
the image $\tilde N$ of $N$ is the normal closure  in $\tilde H$ of the image $\tilde S$ of $S$, and $\tilde M\le \tilde N$.

We now recall Shmelkin's embedding theorem \cite{Sh}. Shmelkin introduced verbal wreath products
and embeded groups of the form $F/V(H)$ into them. Here we adopt his construction
to the variety $\underline{\underline{V}}=\underline{\underline{N}}_c$ of all $c$-nilpotent groups.

Let $X=\{x_i\}_{i\in I}$ be a free basis of $F$, $H$ a normal subgroup in $F$, and $G=F/H$. We denote by $Y$ the alphabet $\{y_{i,g}\}_{i\in I, g\in G}.$
Let $E$ denote the free $c$-nilpotent group with basis $Y$, i.e., the factor group
of the (absolutely) free group $F(Y)$ over  $\gamma_{c+1}(F(Y))$. The group $E$ is free in the variety
$\underline{\underline{N}}_c$ of all nilpotent groups of class $\le c$, and we regard $Y$ as the free basis of $E$ too. We have the right action $\circ$ of $G$ on $Y$ given by the
rule  $y_{i,g}\circ g' = y_{i,gg'}$ for every $g'\in G$. It defines the semidirect product
$W= G\cdot E$ with the action of $G$ on $E$ by conjugation: $(g')^{-1} y_{i,g}g' = y_{i,gg'}.$
($W$ is called the $c$-nilpotent wreath product of $F/\gamma_{c+1}(F)$ and $G$.)

Shmelkin considers the
homomorphism $\mu: F\to W$ defined on  the free generators by $x_i\mapsto \bar x_i y_{i,1}$, where
$\bar x_i$ is the canonical image of $x_i$ in  $G$ and the second subscript in $y_{i,1}$ is the identity of the group $G$. His theorem asserts that
$\ker(\mu)=\gamma_{c+1}(H)$, and so we may regard $\mu$ as an embedding $\tilde F\hookrightarrow W.$

It follows from the above definition that (1) $\mu(\tilde H)$ is a subgroup of $E,$
(2) $\mu(\tilde F)E=W,$ and so the following diagram is commutative.
\medskip
$$
\begin{array}{ccccccccc}
1 &\longrightarrow & \tilde H &
\longrightarrow & \tilde F &\longrightarrow & G
& \longrightarrow 1\\
&&\downarrow &&\downarrow\lefteqn{\mu}&&\downarrow\lefteqn{id} &&\\
1 &\longrightarrow & E &
\longrightarrow & W &\longrightarrow & G
& \longrightarrow 1
\end{array}
$$
\medskip

  Property (1) implies that $\mu(\tilde M)$ belongs to the normal closure of the finite set $\mu(\tilde S)$ in $E.$ Property (2)  means that $\mu(\tilde M)$ is a normal
subgroup of $W$ because it is normalized by $\mu(\tilde F)$ and, being a part of the central subgroup $\mu(\gamma_c(\tilde H))\le\gamma_c(E)$ of $E$, it is also normalized by $E$.

Since $E$ is free nilpotent, for every subset $Z\subset Y$, we have the endomorphism
$\alpha_Z$ of $E$ killing all the generators from $Z$ and leaving fixed the generators
from $Y\backslash Z$. Since the set $\mu(\tilde S)$ is finite, we can fix a finite subset $Z$
so that  $\alpha_Z(\la\mu(\tilde S)\ra^E)=1$, and therefore $\alpha_Z(\mu(\tilde M))=1$.

Let $J\subset G$ denote  the finite set of second indices in the generators $y_{i,g}\in Z.$ Since the group
$G$ is infinite, one can find an element $g'\in G$ such that the set of the second indices $Jg'$ for
all $y_{i,gg'}\in Z'=g'^{-1}Zg'$ is disjoint with $J.$
However, being normal in $W,$
the subgroup $\mu(\tilde M)$ belongs to the normal closure in $E$ of both $\tilde S$ and $g'^{-1}\tilde S g'$,
and so both endomorphisms $\alpha_Z$ and $\alpha_{Z'}$ kill $\mu(\tilde M)$.

On the one hand, arguing in this way we can find infinitely many disjoint subsets $Z, Z', Z'',\dots$ of $G$
such that each of the endomorphism  $\alpha_Z, \alpha_{Z'}, \alpha_{Z'},\dots$ kills
the subgroup $\mu(\tilde M).$ On the other hand, every word $w$ in $Y$ involves no letters from one
of these sets, and so it is left fixed by one of the listed endomorphisms. Thus
$\mu(\tilde M)$ and $\tilde M$ are trivial groups, a contradiction. $\Box$

\begin{rem} \label{inf} Below we will use that for the subgroup $N$ defined in the formulation of
Theorem \ref{shm} and for any nontrivial normal subgroup $M$ of $F$, we have $[M:M\cap N]=\infty$.
Indeed, if this index were $n<\infty$, then the subgroup generated by all the $(n!)$-th powers of the elements
of $M$ would be a nontrivial normal in $F$ subgroup contained in $N,$ contrary to Theorem \ref{shm}.
\end{rem}

\begin{rem}  If $|F/H|<\infty,$ then for any $g\in H$ and a transversal $T$ of $H$ in $F$ (where $F$
is any group), we can take the finite set $S=\{t^{-1}gt\mid t\in T\}.$ Then the normal closure
of $S$ in $H$ is normal in $F$ as well.

If, in addition, $F$ is free, then the normal closure $L$ in $H$ of any single element $g\in H\backslash\{1\}$
contains a nontrivial normal in $F$ subgroup, namely $\cap_{t\in T}t^{-1}Lt$ since the intersection
of finitely many nontrivial normal subgroups of the free group $H$ is nontrivial.
\end{rem}

\medskip

{\bf Proof of Corollary \ref{t0601}.} We will use that the intersections $K\cap L$ of two  nontrivial subnormal subgroups of a free group $F$ is nontrivial. Indeed, let $K=K_s$ and $L=L_t$ be the members of two subnormal series $F =K_0\triangleright K_1\triangleright \dots $ and $F=L_0 \triangleright L_1\triangleright \dots .$ We will induct on
$s+t$ with obvious base if $\min(s,t)=0.$
Assume now that $s,t\ge 1$. Then the intersection $K_s\cap L_t$ is equal to the intersection $R\cap Q$ of $R= K_{s-1}\cap L_t$ and  $Q=K_s\cap L_{t-1}$. Both $R$ and $Q$ are non-trivial by the inductive hypothesis
and normal in the free group $K_{s-1}\cap L_{t-1}$. Since the intersection of nontrivial normal subgroups is nontrivial in a free group, the induction is completed.

 Now let $F \triangleright M_1\triangleright \dots \triangleright M_{\ell}$ be an arbitrary subnormal series in $F$ with nontrivial terms. To prove Corollary \ref{t0601}, we need to show that $M_{\ell}$ is not
contained in $H(\ell)$.

We first observe that setting $N_i = H(i-1) \cap M_i,$ ($i = 1, 2,\dots,$ and $H(0)=H$), the series
$F \triangleright N_1 \triangleright N_2\triangleright \dots$ is also a subnormal series. It was notices
above that the subgroups $N_i$-s are all non-trivial.
As a result, proving Corollary \ref{t0601}, we may assume that
$M_i\le H(i-1)$ for every $i.$ So we need to prove that $M_{\ell}$ is not contained in $H(\ell)$.
Moreover, we will show (to use induction on $\ell$) that $[M_{\ell}:M_{\ell}\cap H(\ell)]=\infty$.
This is true for $\ell=1$ by Remark \ref{inf}. Then we assume that $\ell\ge 2$ and the
inductive hypothesis is true for $\ell-1$. Proving by contradiction, assume that $[M_{\ell}:M_{\ell}\cap H(\ell)]<\infty$.

Let $R$ be the normal closure of $M_{\ell}$ in $H(\ell-1)$. Note that the normal closure of the intersection $M_{\ell}\cap H(\ell)$ is contained in the normal closure $H(\ell)$ of the finite set $S$ in $H(\ell-1)$.
Since the index $[M_{\ell}:M_{\ell}\cap H(\ell)]$ is finite, it follows that $R$ is also contained in the normal
closure of a finite set $T$ in $H(\ell-1).$

Taking into account the inclusion
$M_{\ell-1}\le H(\ell-2)$ we see that $M_{\ell-1}$ normalizes both $M_{\ell}$ and $H(\ell-1)$.
Hence $M_{\ell-1}$ is contained in the normalizer of $R$. Therefore $R$ is a normal subgroup in
the free subgroup  $E=M_{\ell-1}H(\ell-1).$ By Theorem \ref{shm} applied to $E$ and to its normal
subgroup $H(\ell-1)$, the factorgroup  $M_{\ell-1}H(\ell-1)/H(\ell-1)\cong M_{\ell-1}/(M_{\ell-1}\cap H(\ell-1))$ is finite, a contradiction with
the inductive hypothesis. The corollary is proved.
$\Box$
\medskip

{\bf Proof of Corollary \ref{g}.} The normal closure $H$ of $g$  has infinite index in $F$ (even modulo
the derived subgroup $F'$) since the free rank of $F$ is greater than $1$. So $\la g\ra_{\ell+1}^{F}= \la g\ra_{\ell}^{H}$ contains no nontrivial $\ell$-subnormal in $F$ subgroups by Corollary \ref{t0601}. $\Box$

\section{Existence of the e.g. rate.}\label{eg}

{\bf Proof of Lemma \ref{graph}.} Let $\Gamma$ be a connected graph, where for every edge $e$, there exists a unique
inverse edge $e^{-1}\ne e$, and $(e^{-1})^{-1}=e$. We will fix any vertex $o$ in $\Gamma$ and denote by
 $\varphi(n)=\varphi_o(n)$ the number of reduced closed paths $p=e_1\dots e_k$ of lengths $|p|=k\le n$ with $p_-=p_+=o$. Here and
 further in the paper, we denote by
 $p_-$ the original vertex $(e_1)_-$ of the edge $e_1$ and denote by $p_+$ the terminal vertex $(e_k)_+$ of $e_k$.

The well-known construction of the free basis in the fundamental group $\pi_1(\Gamma)$ is the following
(e.g., see \cite{LS}, III.2). At first one chooses a maximal subtree  $T$  in $\Gamma$; there are no reduced closed
paths of positive length in $T$. The tree $T$ contains all the vertices of $\Gamma$.
For every vertex $v$, we denote by $p(v)$ the unique reduced path in $T$ connecting the base
point $o$ and $v.$ Denote by $E$ the set of edges of $\Gamma$ which are not in the tree $T$. The free basis of the fundamental
group $\pi_1(\Gamma)$  is given by the paths $q_e = p(e_-) e p(e_+)^{-1}$, where $e\in E$.
Note that $q_{e^{-1}} = (q_e)^{-1},$ and we call a product $q_{e_1}\dots q_{e_s}$ {\it proper} if
$e_{i-1}\ne e_i^{-1}$ for $i=2,\dots,s$. Thus, every closed reduced path $q$ starting at $o$
is $1$-homotopic to a unique proper product $q_{e_1}\dots q_{e_s}$; in other words, $q$ results
from such a product after a number of cancelations of the edges from $T$ while the edges
$e_1,\dots, e_s$ remain untouched by the cancelations; in particular  $|q|\ge s.$

We may assume that the set $E$ has at least two edges $e$ and $e'$, where $e'\ne e^{\pm 1}$. Indeed, otherwise
the group $\pi_1(\Gamma)$ is cyclic, arbitrary proper product is a power of $q_e,$ the function
$\varphi(n)$ is bounded by a linear function, and so $\lim_{n\to\infty}\sqrt[n]{\varphi(n)}=0.$
Now we fix those two edges and choose an integer  $c \ge \max(|q_e|, |q_{e'}|).$

Let $s$ and $t$ be nonnegative integers and the proper products $q_{e_1}...q_{e_m}$ and $q_{f_1}...q_{f_n}$ represent two reduced loops $p$ and $q$ of
lengths $\le s$ and $\le t,$ respectively. One can choose an edge $f\in \{e^{\pm 1},e'^{\pm 1}\}$
so that $f\ne e_m^{-1}$ and $f\ne f_1^{-1}.$ Then the product $q_{e_1}...q_{e_m}q_f q_{f_1}...q_{f_n}$ is
proper and represents a loop of length $\le s+t+c.$ Thus, given a pair $(p,q)$ of reduced paths of
length $\le s$ and $\le t$ in $\Gamma$, we correspond a reduced form $z$ of $q_{e_1}...q_{e_m}q_f q_{f_1}...q_{f_n}$
 with length $\le s+t+c$. This mapping is not injective, since although the proper decomposition
of a proper product representing $z$ is unique, still there is a choice for
the middle factor $q_f$ in a product $q_{e_1}...q_{e_m}q_f q_{f_1}...q_{f_n}.$ But the
number of such options does not exceed $s+1$ since $m\le s.$ It follows that
\begin{equation}\label{sub}
\varphi(s)\varphi(t)\le (s+1) \varphi(s+t+c)
\end{equation}
for a constant $c$ and any  $s,t.$

On the one hand, the number of arbitrary paths of length at most $n$ starting at $o$ is bounded by
an exponential function of $n$ since the degrees of the vertices
in $\Gamma$ are bounded. Therefore the function $\varphi(n)$ is bounded from above by an exponential one.
On the other hand, the number of the proper products of the paths $q_{e^{\pm 1}}$  and $q_{(e')^{\pm 1}}$
with $n$ factors is bounded from below by $3^n$. Since different proper products provide us with
different reduced forms of length $\le cn$, we have the inequalities $\varphi(n)\ge 3^{\lfloor c^{-1}n\rfloor}$ for all $n\ge 0$.
It follows that there exists a finite upper limit $a=\limsup_{n\to \infty} \sqrt[n]{\varphi(n)}$ and $a>1$.
We set $a_n= \sqrt[n]{\varphi(n)}$, and therefore  $\limsup_{n\to \infty} a_n=a>1.$

Now, given any $\varepsilon\in (0,a-1)$, one can find and fix an integer $s\ge 1$ such that
\begin{equation}\label{as}
|(s+1)^{-1/s}a_s-a|<\varepsilon/3 \;\;\; and \;\;\; (a-2\varepsilon/3)^{\frac{s+c}{s}}< a-\varepsilon/3
\end{equation}
because $\lim _{s\to\infty}(s+1)^{-1/s}=1$ and $\lim _{s\to\infty}\frac{s+c}{s}=1$.

Arbitrary integer $n$  can be presented in the form $n=q(s+c)+r$ for some integers $q$ and $r$
with $0\le r<s+c$. Observe that for any large enough $n$, we have
\begin{equation}\label{n}
(a-2\varepsilon/3)^{n-s-c}>(a-\varepsilon)^n,
\end{equation}

Since $\phi(t)\ge \phi(t')$ for $t\ge t'$, we obtain from the inequality (\ref{sub}):
\begin{equation}\label{power}
\varphi(n)\ge \varphi(q(s+c))\ge (s+1)^{-1}\varphi(s)\varphi((q-1)(s+c))\ge\dots\ge (s+1)^{-q}\varphi(s)^q
\end{equation}

The right-hand side of (\ref{power}) is equal to $(a_s(s+1)^{-1/s})^{sq}$ by the definition of $a_s$.
So for every sufficiently large $n$, it follows from (\ref{power}), (\ref{as}), (\ref{n}), and from the inequality
$r<s+c$ that

$$\varphi(n)\ge (a_s(s+1)^{-1/s})^{sq}\ge (a-\varepsilon/3)^{sq}\ge (a-2\varepsilon/3)^{(s+c)q}$$
$$=(a-2\varepsilon/3)^{n-r} >(a-2\varepsilon/3)^{n-s-c}> (a-\varepsilon)^n$$

Thus $a_n> a-\varepsilon$ for every large enough $n.$ Since $\varepsilon$ can be arbitrary small,
we conclude that $a=\lim_{n\to \infty} a_n$. The first assertion of Lemma \ref{graph} is proved.

If there is a path $p$ connecting two vertices $o$ and $o'$ of $\Gamma$, then to any reduced
closed path $q'$ at $o'$ (i.e., starting and terminating at $o'$), one can correspond the reduced form $q$ of the path $pq'p^{-1}$ at $o$.
Since this mapping is injective and $|q|\le |q'|+2\ell$, where $\ell=|p|$, we have $\varphi_{o'}(n)\le \varphi_o(n+2\ell),$
whence $$\lim_{n\to\infty}\sqrt[n]{\varphi_{o'}(n)}\le \limsup_{n\to\infty}\sqrt[n]{\varphi_o(n+2\ell)}=
\limsup_{n\to\infty}\sqrt[n+2\ell]{\varphi_o(n+2\ell)} =$$ $$= \limsup_{n\to\infty}\sqrt[n]{\varphi_o(n)}=\lim_{n\to\infty}\sqrt[n]{\varphi_o(n)}$$
Similarly we have the opposite inequality, and the second claim of the lemma is proved too.
$\Box$
\medskip

{\bf Proof of Theorem \ref{rate}.} Recall that the cosets $Hg$ are the vertices of the coset graph $\Gamma$ of a
subgroup $H\le F(X)$, and for every $x\in X^{\pm 1}$, every  coset $Hg$ is connected by an edge $e$
with the coset $Hgx$. The edge $e$ is labeled by $x$, and the inverse edge is labeled by $x^{-1}$.
Hence a path $p=e_1\dots e_s$ is reduced in $\Gamma$ if and only if its label $\Lab(p)\equiv \Lab(e_1)\dots \Lab(e_s)$
is a reduced word. (We use the sign '$\equiv$' for letter-by-letter equality of words.) A path $p$ starting at the vertex $H$ is closed if and only if $H\Lab(p)=H,$ i.e.,
iff the label of $p$ represents an element from $H$. Therefore the number of elements in $H$ with
length at most $n$ is equal to the number of reduced closed paths in $\Gamma$ starting at $H$ and having
length at most $n$. Now the statement of Theorem \ref{rate} follows from Lemma \ref{graph} since the
degrees off all vertices in $\Gamma$ are equal to $2m$, where $m=\#X.$
$\Box$

\begin{rem}\label{fa}  %If $X$ is a finite set of generators in a non-free group $G$, then
It may happen in a non-free finitely generated group $G$ that
$\limsup_{n\to\infty}\sqrt[n]{g_H(n)} >1$ while $\liminf_{n\to\infty}\sqrt[n]{g_H(n)} =1$
for the growth function
$g_H(n)$ of a subgroup $H\le G$ with respect to arbitrary finite set of generators of $G$, and so the e.g. rate $\lim_{n\to\infty}\sqrt[n]{g_H(n)}$ does not exist.

This follows from the description of growth functions
given in Theorem 2.1 (3) of \cite{DO}; moreover the counter-example $H$ can be chosen as a cyclic subgroup of a solvable group $G$ with $\liminf_{n\to\infty}g_H(n)/n^{1+\varepsilon} =0$ for every $\varepsilon>0$ .
\end{rem}

\section{$H$-diagrams and subnormal closures of a single element.}\label{dia}

Recall that a van Kampen {\it diagram} $\Delta $ over an alphabet $X$
is a finite, oriented, connected and simply-connected, planar 2-complex endowed with a
labeling function $\Lab : E(\Delta )\to X^{\pm 1}$, where $E(\Delta
) $ denotes the set of oriented edges of $\Delta $, such that $\Lab
(e^{-1})\equiv \Lab (e)^{-1}$. Given a face (that is a 2-cell) $\Pi $ of $\Delta $,
we denote by $\partial \Pi$ the boundary of $\Pi $; similarly,
$\partial \Delta $ denotes the boundary of $\Delta $.
An additional requirement for a diagram over a presentation $G=\langle X\; | \; \mathcal R\rangle$ (or just over the group $G$ given by this presentation) is that the label of any
face $\Pi $ of $\Delta $ is equal to
a word $R^{\pm 1}$, where $R\in \mathcal R$.
This implies that there is a base point $o(\Pi)$ on the boundary $\partial\Pi$ of every
face $\Pi$, and one read the boundary label of $\Pi$ starting with the
vertex $o(\Pi)$.
Labels and lengths of
paths are defined as in Section \ref{eg}.

The van Kampen Lemma states that a word $w$ over the alphabet $X^{\pm 1}$
represents the identity in the group $G$ (i.e., $w$ belongs to the normal closure
of the set $\cal R$ in the free group $F(X)$) if and only
if there exists a diagram $\Delta
$ over $G$ such that
$\Lab (\partial \Delta )\equiv w.$ This implies that one reads the boundary
label of $\Delta$ starting with the base point $o(\Delta)$ which is a vertex
on the boundary of $\Delta$.
(See \cite{LS}, Ch. 5, Theorem 1.1 or \cite{O}, Section 11. The above definition is closer to \cite{O}
since the  edge labels are just letters, whereas such labels are words in \cite{LS}.)

If a diagram $\Delta$ has a pair of faces $\Pi_1$ and $\Pi_2$ whose boundaries share an
edge $e,$ and the labels of the boundaries $\partial\Pi_1$ and $\partial\Pi_2$ coincide
when they are read, respectively, clockwise and counter-clockwise starting with $e,$ then
the diagram $\Delta$ is non-reduced. A diagram is called {\it reduced} if it has no such pairs of faces.
One can replace the word ``diagram'' by ``reduced diagram'' in the formulation of the van Kampen Lemma
(\cite{LS}, Ch. 5 or \cite{O}, Theorem 11.1).

However we want to analyze the words in the alphabet $X^{\pm 1}$ belonging to the normal closure of a subset $\cal R$
in a {\it subgroup}  $H\le F(X)$, not in the whole $F(X)$.
With this purpose, we will modify the classical notion of diagram
(where $H=F(X)$) as follows.

\begin{defn}
Let $\Delta$ be a diagram over a presentation $G=\langle X\; | \; \mathcal R\rangle$ and let
${\cal R}\subset H$, where $H$ is a subgroup of the free group $F=F(X).$ Then we say $\Delta$ is
an $H$-{\it diagram} if for arbitrary its face $\Pi$, every path $p(\Pi)$ connecting
the base points $o(\Delta)$ and $o(\Pi)$ has the label $\Lab(p(\Pi))$ equal in $F$ to an element
of the subgroup $H$.
\end{defn}
It follows that every path of an $H$-diagram connecting two vertices
$o(\Pi)$ and $o(\Pi')$ has label in $H$ as well.

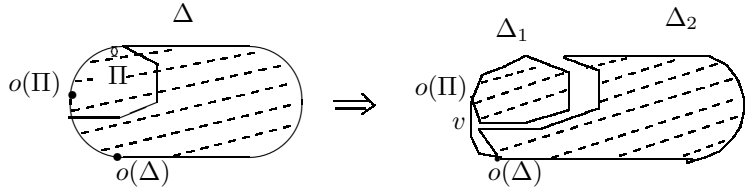
\begin{figure}[h!]
\begin{center}
%TeXCAD (http://texcad.sf.net/) Picture. File: [subsub8.pic]. Options on following lines.
%\grade{\on}
%\emlines{\off}
%\epic{\off}
%\beziermacro{\on}
%\reduce{\on}
%\snapping{\off}
%\pvinsert{% Your \input, \def, etc. here}
%\quality{8.000}
%\graddiff{0.005}
%\snapasp{1}
%\zoom{4.0000}
\unitlength 1mm % = 2.845pt
\linethickness{0.4pt}
\ifx\plotpoint\undefined\newsavebox{\plotpoint}\fi % GNUPLOT compatibility
\begin{picture}(103.25,40)(5,60)
\put(29.125,85.375){\oval(30.75,14.75)[]}
\put(19.625,92.625){\oval(.75,.25)[]}
%\emline(20.75,92.75)(25.5,90.25)
\multiput(20.75,92.75)(.06333333,-.03333333){75}{\line(1,0){.06333333}}
%\end
\put(25.25,90.25){\line(0,-1){5}}
\put(25.25,85.25){\line(-5,-2){5}}
\put(20.25,83.25){\line(-1,0){6.75}}
%\emline(71,78)(68,78.5)
\multiput(71,78)(-.2,.0333333){15}{\line(-1,0){.2}}
%\end
%\emline(68,78.5)(67,80.75)
\multiput(68,78.5)(-.0333333,.075){30}{\line(0,1){.075}}
%\end
\put(67,80.75){\line(0,1){4.5}}
%\emline(67,85.25)(68.5,89)
\multiput(67,85.25)(.03333333,.08333333){45}{\line(0,1){.08333333}}
%\end
%\emline(68.5,89)(71.25,90)
\multiput(68.5,89)(.0916667,.0333333){30}{\line(1,0){.0916667}}
%\end
%\emline(71.25,90)(74.25,91.5)
\multiput(71.25,90)(.06666667,.03333333){45}{\line(1,0){.06666667}}
%\end
%\emline(74.25,91.5)(80,89.25)
\multiput(74.25,91.5)(.0858209,-.03358209){67}{\line(1,0){.0858209}}
%\end
\put(80,89.25){\line(0,-1){5}}
%\emline(80,84.25)(74,82.5)
\multiput(80,84.25)(-.11538462,-.03365385){52}{\line(-1,0){.11538462}}
%\end
\put(74,82.5){\line(-1,0){5.75}}
%\emline(68.25,82.5)(67,86)
\multiput(68.25,82.5)(-.03289474,.09210526){38}{\line(0,1){.09210526}}
%\end
%\emline(79.25,91.25)(84,89.5)
\multiput(79.25,91.25)(.09134615,-.03365385){52}{\line(1,0){.09134615}}
%\end
\put(84,89.5){\line(0,-1){5}}
%\emline(84,84.5)(76.25,81.75)
\multiput(84,84.5)(-.0945122,-.03353659){82}{\line(-1,0){.0945122}}
%\end
\put(76.25,81.75){\line(-1,0){8.25}}
%\emline(68,81.75)(70.5,78.25)
\multiput(68,81.75)(.03333333,-.04666667){75}{\line(0,-1){.04666667}}
%\end
\put(79.5,91.5){\line(1,0){19}}
%\emline(98.5,91.5)(101,90.25)
\multiput(98.5,91.5)(.06578947,-.03289474){38}{\line(1,0){.06578947}}
%\end
%\emline(101,90.25)(102.5,88.25)
\multiput(101,90.25)(.03333333,-.04444444){45}{\line(0,-1){.04444444}}
%\end
%\emline(102.5,88.25)(103.25,86)
\multiput(102.5,88.25)(.0326087,-.0978261){23}{\line(0,-1){.0978261}}
%\end
\put(103.25,86){\line(0,-1){2}}
%\emline(103.25,84)(102.75,82)
\multiput(103.25,84)(-.0333333,-.1333333){15}{\line(0,-1){.1333333}}
%\end
%\emline(102.75,82)(101.75,80.5)
\multiput(102.75,82)(-.0333333,-.05){30}{\line(0,-1){.05}}
%\end
%\emline(101.75,80.5)(100,78.75)
\multiput(101.75,80.5)(-.03365385,-.03365385){52}{\line(0,-1){.03365385}}
%\end
%\emline(100,78.75)(98.25,78)
\multiput(100,78.75)(-.076087,-.0326087){23}{\line(-1,0){.076087}}
%\end
%\emline(98.25,78)(95.75,77.25)
\multiput(98.25,78)(-.1086957,-.0326087){23}{\line(-1,0){.1086957}}
%\end
\put(95.75,77.25){\line(0,1){.5}}
\put(96.5,77.75){\line(-1,0){25.5}}
\put(18.75,87.75){$\Pi$}
\put(20,78){\circle*{1.118}}
\put(70.5,77.75){\circle*{.707}}
%\dashline{1}(16.5,90.75)(28.25,92.5)
\multiput(16.43,90.68)(.195833,.029167){5}{\line(1,0){.195833}}
\multiput(18.388,90.971)(.195833,.029167){5}{\line(1,0){.195833}}
\multiput(20.346,91.263)(.195833,.029167){5}{\line(1,0){.195833}}
\multiput(22.305,91.555)(.195833,.029167){5}{\line(1,0){.195833}}
\multiput(24.263,91.846)(.195833,.029167){5}{\line(1,0){.195833}}
\multiput(26.221,92.138)(.195833,.029167){5}{\line(1,0){.195833}}
%\end
%\dashline{1}(14.5,88.25)(16.75,88.5)
\put(14.43,88.18){\line(1,0){.75}}
\put(15.93,88.346){\line(1,0){.75}}
%\end
%\dashline{1}(21.25,88.75)(38.25,92)
\multiput(21.18,88.68)(.149123,.028509){6}{\line(1,0){.149123}}
\multiput(22.969,89.022)(.149123,.028509){6}{\line(1,0){.149123}}
\multiput(24.759,89.364)(.149123,.028509){6}{\line(1,0){.149123}}
\multiput(26.548,89.706)(.149123,.028509){6}{\line(1,0){.149123}}
\multiput(28.338,90.048)(.149123,.028509){6}{\line(1,0){.149123}}
\multiput(30.127,90.39)(.149123,.028509){6}{\line(1,0){.149123}}
\multiput(31.917,90.732)(.149123,.028509){6}{\line(1,0){.149123}}
\multiput(33.706,91.074)(.149123,.028509){6}{\line(1,0){.149123}}
\multiput(35.495,91.417)(.149123,.028509){6}{\line(1,0){.149123}}
\multiput(37.285,91.759)(.149123,.028509){6}{\line(1,0){.149123}}
%\end
%\dashline{1}(14,84.25)(41.75,90.25)
\multiput(13.93,84.18)(.1367,.029557){7}{\line(1,0){.1367}}
\multiput(15.844,84.594)(.1367,.029557){7}{\line(1,0){.1367}}
\multiput(17.757,85.007)(.1367,.029557){7}{\line(1,0){.1367}}
\multiput(19.671,85.421)(.1367,.029557){7}{\line(1,0){.1367}}
\multiput(21.585,85.835)(.1367,.029557){7}{\line(1,0){.1367}}
\multiput(23.499,86.249)(.1367,.029557){7}{\line(1,0){.1367}}
\multiput(25.412,86.662)(.1367,.029557){7}{\line(1,0){.1367}}
\multiput(27.326,87.076)(.1367,.029557){7}{\line(1,0){.1367}}
\multiput(29.24,87.49)(.1367,.029557){7}{\line(1,0){.1367}}
\multiput(31.154,87.904)(.1367,.029557){7}{\line(1,0){.1367}}
\multiput(33.068,88.318)(.1367,.029557){7}{\line(1,0){.1367}}
\multiput(34.981,88.731)(.1367,.029557){7}{\line(1,0){.1367}}
\multiput(36.895,89.145)(.1367,.029557){7}{\line(1,0){.1367}}
\multiput(38.809,89.559)(.1367,.029557){7}{\line(1,0){.1367}}
\multiput(40.723,89.973)(.1367,.029557){7}{\line(1,0){.1367}}
%\end
%\dashline{1}(15.25,81.25)(43.5,88)
\multiput(15.18,81.18)(.134524,.032143){7}{\line(1,0){.134524}}
\multiput(17.063,81.63)(.134524,.032143){7}{\line(1,0){.134524}}
\multiput(18.946,82.08)(.134524,.032143){7}{\line(1,0){.134524}}
\multiput(20.83,82.53)(.134524,.032143){7}{\line(1,0){.134524}}
\multiput(22.713,82.98)(.134524,.032143){7}{\line(1,0){.134524}}
\multiput(24.596,83.43)(.134524,.032143){7}{\line(1,0){.134524}}
\multiput(26.48,83.88)(.134524,.032143){7}{\line(1,0){.134524}}
\multiput(28.363,84.33)(.134524,.032143){7}{\line(1,0){.134524}}
\multiput(30.246,84.78)(.134524,.032143){7}{\line(1,0){.134524}}
\multiput(32.13,85.23)(.134524,.032143){7}{\line(1,0){.134524}}
\multiput(34.013,85.68)(.134524,.032143){7}{\line(1,0){.134524}}
\multiput(35.896,86.13)(.134524,.032143){7}{\line(1,0){.134524}}
\multiput(37.78,86.58)(.134524,.032143){7}{\line(1,0){.134524}}
\multiput(39.663,87.03)(.134524,.032143){7}{\line(1,0){.134524}}
\multiput(41.546,87.48)(.134524,.032143){7}{\line(1,0){.134524}}
%\end
%\dashline{1}(18,79)(44,85.25)
\multiput(17.93,78.93)(.132653,.031888){7}{\line(1,0){.132653}}
\multiput(19.787,79.376)(.132653,.031888){7}{\line(1,0){.132653}}
\multiput(21.644,79.823)(.132653,.031888){7}{\line(1,0){.132653}}
\multiput(23.501,80.269)(.132653,.031888){7}{\line(1,0){.132653}}
\multiput(25.358,80.715)(.132653,.031888){7}{\line(1,0){.132653}}
\multiput(27.215,81.162)(.132653,.031888){7}{\line(1,0){.132653}}
\multiput(29.073,81.608)(.132653,.031888){7}{\line(1,0){.132653}}
\multiput(30.93,82.055)(.132653,.031888){7}{\line(1,0){.132653}}
\multiput(32.787,82.501)(.132653,.031888){7}{\line(1,0){.132653}}
\multiput(34.644,82.948)(.132653,.031888){7}{\line(1,0){.132653}}
\multiput(36.501,83.394)(.132653,.031888){7}{\line(1,0){.132653}}
\multiput(38.358,83.84)(.132653,.031888){7}{\line(1,0){.132653}}
\multiput(40.215,84.287)(.132653,.031888){7}{\line(1,0){.132653}}
\multiput(42.073,84.733)(.132653,.031888){7}{\line(1,0){.132653}}
%\end
%\dashline{1}(27.5,78.25)(43,82.25)
\multiput(27.43,78.18)(.130252,.033613){7}{\line(1,0){.130252}}
\multiput(29.253,78.65)(.130252,.033613){7}{\line(1,0){.130252}}
\multiput(31.077,79.121)(.130252,.033613){7}{\line(1,0){.130252}}
\multiput(32.9,79.591)(.130252,.033613){7}{\line(1,0){.130252}}
\multiput(34.724,80.062)(.130252,.033613){7}{\line(1,0){.130252}}
\multiput(36.547,80.533)(.130252,.033613){7}{\line(1,0){.130252}}
\multiput(38.371,81.003)(.130252,.033613){7}{\line(1,0){.130252}}
\multiput(40.194,81.474)(.130252,.033613){7}{\line(1,0){.130252}}
\multiput(42.018,81.944)(.130252,.033613){7}{\line(1,0){.130252}}
%\end
%\dashline{1}(68.5,88)(76.5,90.25)
\multiput(68.43,87.93)(.114286,.032143){7}{\line(1,0){.114286}}
\multiput(70.03,88.38)(.114286,.032143){7}{\line(1,0){.114286}}
\multiput(71.63,88.83)(.114286,.032143){7}{\line(1,0){.114286}}
\multiput(73.23,89.28)(.114286,.032143){7}{\line(1,0){.114286}}
\multiput(74.83,89.73)(.114286,.032143){7}{\line(1,0){.114286}}
%\end
%\dashline{1}(67.5,85.25)(79.25,89.25)
\multiput(67.43,85.18)(.090385,.030769){10}{\line(1,0){.090385}}
\multiput(69.237,85.795)(.090385,.030769){10}{\line(1,0){.090385}}
\multiput(71.045,86.41)(.090385,.030769){10}{\line(1,0){.090385}}
\multiput(72.853,87.026)(.090385,.030769){10}{\line(1,0){.090385}}
\multiput(74.66,87.641)(.090385,.030769){10}{\line(1,0){.090385}}
\multiput(76.468,88.257)(.090385,.030769){10}{\line(1,0){.090385}}
\multiput(78.276,88.872)(.090385,.030769){10}{\line(1,0){.090385}}
%\end
%\dashline{1}(68.5,83.25)(79.75,87)
\multiput(68.43,83.18)(.096154,.032051){9}{\line(1,0){.096154}}
\multiput(70.16,83.757)(.096154,.032051){9}{\line(1,0){.096154}}
\multiput(71.891,84.334)(.096154,.032051){9}{\line(1,0){.096154}}
\multiput(73.622,84.91)(.096154,.032051){9}{\line(1,0){.096154}}
\multiput(75.353,85.487)(.096154,.032051){9}{\line(1,0){.096154}}
\multiput(77.084,86.064)(.096154,.032051){9}{\line(1,0){.096154}}
\multiput(78.814,86.641)(.096154,.032051){9}{\line(1,0){.096154}}
%\end
%\dashline{1}(82.5,90)(87.75,91.25)
\multiput(82.43,89.93)(.125,.029762){6}{\line(1,0){.125}}
\multiput(83.93,90.287)(.125,.029762){6}{\line(1,0){.125}}
\multiput(85.43,90.644)(.125,.029762){6}{\line(1,0){.125}}
\multiput(86.93,91.001)(.125,.029762){6}{\line(1,0){.125}}
%\end
%\dashline{1}(84.5,88)(95,91.25)
\multiput(84.43,87.93)(.097222,.030093){9}{\line(1,0){.097222}}
\multiput(86.18,88.471)(.097222,.030093){9}{\line(1,0){.097222}}
\multiput(87.93,89.013)(.097222,.030093){9}{\line(1,0){.097222}}
\multiput(89.68,89.555)(.097222,.030093){9}{\line(1,0){.097222}}
\multiput(91.43,90.096)(.097222,.030093){9}{\line(1,0){.097222}}
\multiput(93.18,90.638)(.097222,.030093){9}{\line(1,0){.097222}}
%\end
%\dashline{1}(84.25,85.5)(99.75,90.5)
\multiput(84.18,85.43)(.101307,.03268){9}{\line(1,0){.101307}}
\multiput(86.003,86.018)(.101307,.03268){9}{\line(1,0){.101307}}
\multiput(87.827,86.606)(.101307,.03268){9}{\line(1,0){.101307}}
\multiput(89.65,87.194)(.101307,.03268){9}{\line(1,0){.101307}}
\multiput(91.474,87.783)(.101307,.03268){9}{\line(1,0){.101307}}
\multiput(93.297,88.371)(.101307,.03268){9}{\line(1,0){.101307}}
\multiput(95.121,88.959)(.101307,.03268){9}{\line(1,0){.101307}}
\multiput(96.944,89.547)(.101307,.03268){9}{\line(1,0){.101307}}
\multiput(98.768,90.136)(.101307,.03268){9}{\line(1,0){.101307}}
%\end
%\dashline{1}(71,79)(101.75,88.5)
\multiput(70.93,78.93)(.103535,.031987){9}{\line(1,0){.103535}}
\multiput(72.793,79.505)(.103535,.031987){9}{\line(1,0){.103535}}
\multiput(74.657,80.081)(.103535,.031987){9}{\line(1,0){.103535}}
\multiput(76.521,80.657)(.103535,.031987){9}{\line(1,0){.103535}}
\multiput(78.384,81.233)(.103535,.031987){9}{\line(1,0){.103535}}
\multiput(80.248,81.808)(.103535,.031987){9}{\line(1,0){.103535}}
\multiput(82.112,82.384)(.103535,.031987){9}{\line(1,0){.103535}}
\multiput(83.975,82.96)(.103535,.031987){9}{\line(1,0){.103535}}
\multiput(85.839,83.536)(.103535,.031987){9}{\line(1,0){.103535}}
\multiput(87.702,84.112)(.103535,.031987){9}{\line(1,0){.103535}}
\multiput(89.566,84.687)(.103535,.031987){9}{\line(1,0){.103535}}
\multiput(91.43,85.263)(.103535,.031987){9}{\line(1,0){.103535}}
\multiput(93.293,85.839)(.103535,.031987){9}{\line(1,0){.103535}}
\multiput(95.157,86.415)(.103535,.031987){9}{\line(1,0){.103535}}
\multiput(97.021,86.99)(.103535,.031987){9}{\line(1,0){.103535}}
\multiput(98.884,87.566)(.103535,.031987){9}{\line(1,0){.103535}}
\multiput(100.748,88.142)(.103535,.031987){9}{\line(1,0){.103535}}
%\end
%\dashline{1}(77.75,78)(102.75,86)
\multiput(77.68,77.93)(.099206,.031746){9}{\line(1,0){.099206}}
\multiput(79.465,78.501)(.099206,.031746){9}{\line(1,0){.099206}}
\multiput(81.251,79.073)(.099206,.031746){9}{\line(1,0){.099206}}
\multiput(83.037,79.644)(.099206,.031746){9}{\line(1,0){.099206}}
\multiput(84.823,80.215)(.099206,.031746){9}{\line(1,0){.099206}}
\multiput(86.608,80.787)(.099206,.031746){9}{\line(1,0){.099206}}
\multiput(88.394,81.358)(.099206,.031746){9}{\line(1,0){.099206}}
\multiput(90.18,81.93)(.099206,.031746){9}{\line(1,0){.099206}}
\multiput(91.965,82.501)(.099206,.031746){9}{\line(1,0){.099206}}
\multiput(93.751,83.073)(.099206,.031746){9}{\line(1,0){.099206}}
\multiput(95.537,83.644)(.099206,.031746){9}{\line(1,0){.099206}}
\multiput(97.323,84.215)(.099206,.031746){9}{\line(1,0){.099206}}
\multiput(99.108,84.787)(.099206,.031746){9}{\line(1,0){.099206}}
\multiput(100.894,85.358)(.099206,.031746){9}{\line(1,0){.099206}}
%\end
%\dashline{1}(86.75,78.25)(102.5,83.5)
\multiput(86.68,78.18)(.092647,.030882){10}{\line(1,0){.092647}}
\multiput(88.533,78.797)(.092647,.030882){10}{\line(1,0){.092647}}
\multiput(90.386,79.415)(.092647,.030882){10}{\line(1,0){.092647}}
\multiput(92.239,80.033)(.092647,.030882){10}{\line(1,0){.092647}}
\multiput(94.091,80.65)(.092647,.030882){10}{\line(1,0){.092647}}
\multiput(95.944,81.268)(.092647,.030882){10}{\line(1,0){.092647}}
\multiput(97.797,81.886)(.092647,.030882){10}{\line(1,0){.092647}}
\multiput(99.65,82.503)(.092647,.030882){10}{\line(1,0){.092647}}
\multiput(101.503,83.121)(.092647,.030882){10}{\line(1,0){.092647}}
%\end
%\dashline{1}(95.25,78.25)(100.75,80.25)
\multiput(95.18,78.18)(.087302,.031746){9}{\line(1,0){.087302}}
\multiput(96.751,78.751)(.087302,.031746){9}{\line(1,0){.087302}}
\multiput(98.323,79.323)(.087302,.031746){9}{\line(1,0){.087302}}
\multiput(99.894,79.894)(.087302,.031746){9}{\line(1,0){.087302}}
%\end
\put(14,86.5){\circle*{.5}}
\put(27.25,96){$\Delta$}
\put(70.25,94){$\Delta_1$}
\put(92.75,95.5){$\Delta_2$}
\put(19.75,75){$o(\Delta)$}
\put(5.5,86.75){$o(\Pi)$}
\put(64.5,81.75){$v$}
\put(14,86.25){\circle*{1.118}}
\put(48.75,85.5){\line(1,0){4.75}}
\put(49,84.5){\line(1,0){4}}
%\emline(52.5,86.5)(54,85)
\multiput(52.5,86.5)(.03333333,-.03333333){45}{\line(0,-1){.03333333}}
%\end
%\emline(52.5,83.75)(54.25,85)
\multiput(52.5,83.75)(.04605263,.03289474){38}{\line(1,0){.04605263}}
%\end
\put(69.25,75){$o(\Delta)$}
\put(59.5,86.25){$o(\Pi)$}
\end{picture}
\end{center}
\caption{Splitting of the diagram $\Delta$}\label{spl}
\end{figure}

\begin{lemma} \label{vK}(1) The boundary label of an $H$-diagram $\Delta$ represents an element of the
normal closure $N=\la{\cal R}\ra^H$ of the set $\cal R$ in the subgroup $H$.

(2) For every word $w$ representing an element of $N$, there is a reduced $H$-diagram over the presentation
$G=\langle X\; | \; \mathcal R\rangle$ whose
boundary label is letter-by-letter equal to $w$.
\end{lemma}
\proof (1) If $\Delta$ has no faces, then its boundary label is trivial in $F$, and it is nothing
to prove. Then we induct on the number of faces $f$ in $\Delta$ assuming that $f\ge 1$.
One can cut off one face from $\Delta$; in other words, there are two subdiagrams $\Delta_1$
and $\Delta_2$ with the base points $o(\Delta_1)=o(\Delta_2)=o(\Delta)$ such that $\Delta_1$
has one face, $\Delta_2$ has $f-1$ faces, and the boundary label $w$ of $\Delta$ is equal in $F$
to the product of the boundary labels $w_1$ and $w_2$ of $\Delta_1$ and $\Delta_2$, respectively
(see Fig.\ref{spl}).

Moreover, $w_1$ is the product $vR^{\pm 1}v^{-1}$, where $R\in \cal R$ and $v$ is
the label of a simple path connecting $o(\Delta)$ with the vertex $o(\Pi)$ of the face $\Pi$  of $\Delta_1$.
Hence $v\in H$ and $w_1\in N$. The word $w_2$ is in $N$ by the inductive hypothesis, and so
$w=w_1w_2\in N$ too.

(2) Since $w\in N$, the word $w$ must be equal in $H$ (and in $F$) to a product $v\equiv\prod_{i=1}^m v_iR_i^{\pm 1}v_i^{-1},$
where $R_i\in\cal R$ and $v_i\in H$. Therefore one can construct a diagram $\Gamma$ as a ``bouquet'' of
one-faced subdiagrams $\Gamma_i$ ($i=1,\dots,m$) with faces $\Pi_i$ having boundary labels $R_i^{\pm 1}$,
where $o(\Gamma$) is connected with $o(\Pi_i)$ by a simple path $p(o(\Pi_i))$ labeled by the word $v_i$.

Since arbitrary closed path of $\Gamma$ starting at $o(\Gamma)$ is homotopic in the $1$-skeleton of $\Gamma,$
to some product of the boundary labels (and inverses) of the subdiagrams $\Gamma_i$-s, its label belongs
to $N.$ It follows that the labels of any paths connecting $o(\Gamma)$ with $o(\Pi_i)$ belong to the
same right coset of $N$ in $F.$  Since $\Lab(p(o(\Pi_i)))\equiv v_i\in H,$ for every  $i$, the labels of all
paths going from $o(\Gamma)$ to $o(\Pi_i)$-s belong to $H$ too. Therefore $\Gamma$ is an $H$-diagram.

The words $w$ and $v$ are freely equal. So one can modify $\Gamma$ by a number of elementary transformations
and obtain a diagram $\Delta$ with boundary label $w$. Every such a transformation either identifies a boundary edge $e$ with the inverse of the subsequent edge
or inserts a pinch $ee^{-1}$ in the boundary
path. It is clear that every path from $o$ to $o(\Pi_i)$ in the modified diagrams have the same  label in $F$ as a path in
the original diagram. So $\Delta$ is an $H$-diagram too.

Suppose the diagram $\Delta$ is not reduced, i.e., it has a pair of ``mirror'' faces with the common edge $e$
as defined above. Then these two faces can be ``canceled'', which gives a diagram $\Delta'$ with the same
boundary label and with fewer faces. (One removes these two cells and the edge $e$ from $\Delta$ and then sews up the
hole using that the boundary label of the hole is trivial in $F$. More accurate details are provided in \cite{O}, Section 11.) However it is easy to see that every path $p'$
connecting $o(\Delta')$ and $o(\Pi_i)$ in $\Delta'$ has the same label in $F$ as a suitable path $p$ in $\Delta$
connecting $o(\Delta)$ with $o(\Pi_i)$ (See the end of Chapter 4 \cite{O} either, where the way to construct $p$ for a given $p'$ is described.) Hence $\Delta'$ is also an $H$-diagram. Thus in a number of steps, one
gets a reduced $H$-diagram with boundary label $w$.
\endproof

{\bf Proof of Theorem \ref{intersec}.}
{\bf Part}(1). One may assume that $g\ne 1$. Suppose that $h\in \cap_{\ell=1}^{\infty} \langle g\rangle_{\ell}^F$.
%Then $h\in \langle g\rangle_{\ell}^F\gamma_{\ell+1}(F)$ for any $\ell\ge 1$.
Note that for any $\ell\ge 1$, the subgroup $\langle g\rangle \gamma_{\ell+1} (F)$ is $\ell$-subnormal in $F$ due to the subnormal series
$$F=\langle g\rangle\gamma_1(F)\ge \langle g\rangle\gamma_2(F)\ge\dots\ge\langle g\rangle\gamma_{\ell+1}(F)$$
Therefore $\langle g\rangle_{\ell}^F \le \langle g \rangle\gamma_{\ell+1}(F)$, and so $h\in
\langle g \rangle\gamma_{\ell+1}(F)$ for every
$\ell\ge 1$.

 Hence we obtain the series of equalities $h=g^{n_{\ell}}u_{\ell}$, where $u_{\ell}\in \gamma_{\ell+1}(F)$. Thus, $g^{n_{\ell}-n_k}\in \gamma_{\ell+1}(F)$ for all $k\ge\ell\ge 1.$ But $\langle g\rangle\cap \gamma_{\ell+1}(F)=1$ for all large enough $\ell$
since $\cap_{\ell=1}^{\infty}\gamma_{\ell+1}(F)=1$ and the factor-groups
$F/\gamma_{\ell+1}(F)$ are torsion free (see \cite{MKS},
Ch. 5). Therefore there is $\ell$ such that
$n_{\ell}=n_{\ell+1}=\dots$, and so $u_{\ell}=u_{\ell+1}=\dots\in \cap_{i=\ell}^{\infty}\gamma_{i+1}(F)=1$.
Hence $h=g^{n_{\ell}}$, which proves the first part of Theorem \ref{intersec}.
\medskip

{\bf Part }(2). One may assume that the word $g=g(x_1,\dots,x_m)$ is nonempty and cyclically reduced because a conjugation can increase the
lengths of all words at most by a constant. To continue, we need

\begin{lemma}\label{l0} There is a natural number $\ell_0=\ell_0(g)$ such that the length
of arbitrary element $h\in \langle g\rangle_{\ell}^F\backslash \langle g\rangle$ is at least $|g|$ provided $\ell\ge\ell_0$.
\end{lemma}

\proof As in part (1) of the proof of Theorem \ref{intersec}, we have
   $\langle g\rangle\cap \gamma_{l+1}(F_m)=1$ if $\ell$ is large enough and $h=g^nu$
for some integer $n$ and $u\in \gamma_{\ell+1}(F).$

If $n=0$, then $h\in \gamma_{\ell+1}(F)$, and according to Fox' estimate (\cite{Fo}, Lemma 4.2) of the lengths
of nontrivial elements in the terms of the lower central series of free groups, we have
$|h|\ge (\ell+1)/2\ge |g|$ provided $\ell_0\ge 2|g|-1$.

If $n\ne 0$, we consider the epimorphism $\phi: F\to F_m$ identical on $x_1,\dots,x_m$
and trivial on other free generators of $F$. Obviously, we obtain $\phi(g)=g$,
$\phi(\langle g\rangle_{i}^F)=\langle g\rangle_{i}^{F_{m}}$ and $\phi(\gamma_i(F))=\gamma_i(F_m)$ for every $i$, and so
$\phi(\langle g\rangle)\cap \gamma_{\ell+1}(F_m)=1$. It follows that $h'=\phi(h)=\phi(g^n)\phi(u)\in \langle g\rangle_{\ell}^{F_m}\backslash \{1\}$.

Since there are finitely many elements of length $\le |g|$ in $F_m$, we can
derive from part (1) that $|h|\ge |h'|\ge |g|$ if $\ell_0$ is
sufficiently large and $h'\notin \langle g\rangle$. If $h'\in\langle g\rangle,$
then $|h'|\ge |g|$ and so $|h|\ge |h'|\ge |g|$ as well.

Thus, the desired choice of $\ell_0$ is possible, and the lemma is proved. \endproof

To complete  the proof of Theorem \ref{intersec} (2), we choose $\ell_0$
in accordance with Lemma \ref{l0}
and define $c=2^{-\ell_0}$.
Then the assertion of Theorem \ref{intersec} holds for $\ell\le \ell_0$, and proving by induction on $\ell$ one may assume that $\ell>\ell_0$.

Suppose $h$ is a reduced word from $\langle g \rangle ^F_{\ell}\backslash \langle g \rangle$. By Lemma \ref{vK}
(2), there is a reduced $H$-diagram $\Delta$ with boundary
label $h$, where $H=\langle g \rangle ^F_{\ell-1}.$

At first we assume that $\Delta$ has only one face $\Pi$, and so the
boundary path is of the form $pqp^{-1}$, where $|q|=|\partial\Pi|=|g|$
and $|p|\ge 0.$ The vertex $p_+$ can be connected with the base point
of $\Pi$ by a path $r$ of length $\le |g|/2$. Therefore the word $h$
is freely equal to $\Lab(pr)g^{\pm 1}\Lab(pr)^{-1}$. By the definition
of $H$-diagram, $\Lab(pr)\in H.$ But the word $\Lab(pr)$ does not belong to $\langle g\rangle$ since $h\notin\langle g\rangle$. Hence by the inductive hypothesis,
$|pr|\ge c 2^{l-1}$, and so $|p|\ge c 2^{l-1}-|g|/2$, whence
$|h|=2|p|+|g|\ge c2^l$, as required. Thus, one may assume further
that $\Delta$ has at least two faces.

We claim that in $\Delta$, no face $\Pi_1$ is attached to a face $\Pi_2$ along
an edge $e$. Indeed, otherwise the base points $o(\Pi_1)$ and $o(\Pi_2)$ are
connected in $\partial\Pi_1\cup\partial\Pi_2$ by a path $t$ of length $<|g|.$
Since $\Delta$ is an $H$-diagram, the choice of $\ell_0$ implies
that $\Lab(t)$ is freely equal to $1$.  In other words, the vertex $e_-$ becomes the base point of both $\Pi_1$ and $\Pi_2$ if one replace $g$ by a cyclic permutation $g'$. It follows that either the diagram $\Delta$ is not reduced or $g'$ starts with $\Lab(e)$ and ends with $\Lab(e)^{-1}$. The former case is
impossible, the latter one means that the word $g$ is not cyclically reduced,
a contradiction again.

Now we present the boundary path $p$ of $\Delta$ as $t_0r_1t_1\dots r_kt_k$,
where $r_1,\dots, r_k$ are subpaths of positive length in the boundaries
of the faces of $\Delta$, and any edge $e$ of $t_0,\dots,t_k$ is a bridge edge i.e., the
edge $e^{-1}$ also occurs in the boundary path of $\Delta$.  (Some $t_i$-s
may have zero length.) It is easy to see that $e$ and $e^{-1}$ cannot occur in the same path $t_i$. Since $\Delta$ has at least two faces and the
faces do not share boundary edges, we obtain
for any $i\le k$ that
\begin{equation}\label{tig}
|p|\ge 2|t_i|+2|g|
\end{equation}

Since every $r_i$ is a boundary arc of a face, one can find a path $\bar r_i=
r'_ir_ir''_i$ on the boundary of this face starting and ending at its base point, such that $|r'_i|,|r''_i|\le |g|/2.$ Note that $\Lab(\bar r_i)\in \langle g\rangle$. Let $\bar t_i= (r''_i)^{-1}t_i(r'_{i+1})^{-1}$, where $r''_0$ and $r'_{k+1}$ have length $0$. Then the path $\bar p=\bar t_0\bar r_1\bar t_1\dots \bar r_k\bar t_k$ has the same labels in $F$ as $p,$ and $\Lab (\bar t_i)\in H$ for every $i$ since $\Delta$ is an $H$-diagram.

If $\Lab (\bar t_i)\in \langle g \rangle $ for every $i$, then $\Lab (\bar p)$
and $\Lab(p)\equiv h$ belong to $\langle g \rangle $ too, a contradiction.
Therefore $\Lab (\bar t_i)\notin \langle g \rangle $ for some $i$, and by the
inductive hypothesis, $|\bar t_i|\ge c 2^{l-1}.$ Hence $|t_i|\ge c 2^{l-1}-
|r''_i|-|r'_i|\ge c 2^{l-1}-|g|$. Therefore by (\ref{tig}), we have $|h|=|p|\ge 2(c 2^{l-1}-|g|)+2|g|=c 2^l,$ which completes the induction. $\Box$

\medskip

\section{The growth of subnormal subgroups in $F_m$.}\label{subgrowth}

In this section, we treat
presentations with only one relation:
$G=\la X\mid x^d\ra$ , where $d\ge 1$ and
$x\in X.$ Note that if $\Delta$ is a reduced diagram over $G$, then the boundaries
of two distinct faces of $\Delta$ cannot share an edge. In other words, $\Delta$ has only
boundary edges and has no inner edges. In particular, the perimeter
$|\partial\Delta|$ is at least $fd$, where $f$ is the number of faces in $\Delta$.

Erasing all edge labels in a diagram
one obtains a {\it map} that is an (unlabeled)
finite, oriented, plane, connected, and simply-connected 2-complex with a base point and with
a base point for every face. Here we say
that two diagrams have equal {\it types} if there is an isomorphism between the
corresponding maps.

The diagrams under consideration may have non-reduced boundary paths
since we need such diagrams for inductive estimates.

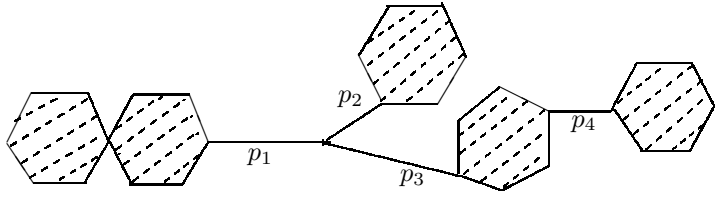
\begin{figure}[h!]
\begin{center}
%TeXCAD (http://texcad.sf.net/) Picture. File: [subsub9.pic]. Options on following lines.
%\grade{\on}
%\emlines{\off}
%\epic{\off}
%\beziermacro{\on}
%\reduce{\on}
%\snapping{\off}
%\pvinsert{% Your \input, \def, etc. here}
%\quality{8.000}
%\graddiff{0.005}
%\snapasp{1}
%\zoom{4.0000}
\unitlength 1mm % = 2.845pt
\linethickness{0.4pt}
\ifx\plotpoint\undefined\newsavebox{\plotpoint}\fi % GNUPLOT compatibility
\begin{picture}(111.25,45)(10,60)
\put(17.25,83.75){\line(1,2){3.25}}
\put(20.5,90.25){\line(1,0){7.5}}
%\emline(28,90.25)(30.75,83.75)
\multiput(28,90.25)(.03353659,-.07926829){82}{\line(0,-1){.07926829}}
%\end
%\emline(30.75,83.75)(27.5,78.25)
\multiput(30.75,83.75)(-.033505155,-.056701031){97}{\line(0,-1){.056701031}}
%\end
\put(17.25,83.5){\line(2,-3){3.5}}
\put(20.75,78.25){\line(1,0){7}}
%\emline(30.75,84)(34,90)
\multiput(30.75,84)(.033505155,.06185567){97}{\line(0,1){.06185567}}
%\end
\put(34,90){\line(1,0){7.5}}
\put(41.5,90){\line(2,-5){2.5}}
%\emline(44,83.75)(40.5,78)
\multiput(44,83.75)(-.033653846,-.055288462){104}{\line(0,-1){.055288462}}
%\end
\put(40.5,78){\line(-1,0){7}}
%\emline(30.75,84)(33.75,78.25)
\multiput(30.75,84)(.03370787,-.06460674){89}{\line(0,-1){.06460674}}
%\end
\put(59.25,83.25){\line(0,1){.75}}
\put(44,83.75){\line(1,0){16.25}}
%\emline(59.75,84)(67,88.75)
\multiput(59.75,84)(.05141844,.033687943){141}{\line(1,0){.05141844}}
%\end
\put(67,88.75){\line(1,0){7.5}}
\put(74.5,88.75){\line(3,5){3.75}}
%\emline(78.25,95)(75,102.25)
\multiput(78.25,95)(-.033505155,.074742268){97}{\line(0,1){.074742268}}
%\end
\put(67.75,101.75){\line(1,0){7.25}}
%\emline(68,101.75)(64,95.25)
\multiput(68,101.75)(-.033613445,-.054621849){119}{\line(0,-1){.054621849}}
%\end
\put(64,95.25){\line(1,-2){3}}
%\emline(59.5,83.5)(78,79)
\multiput(59.5,83.5)(.138059701,-.03358209){134}{\line(1,0){.138059701}}
%\end
%\emline(78,79)(83,77.25)
\multiput(78,79)(.09615385,-.03365385){52}{\line(1,0){.09615385}}
%\end
%\emline(83,77.25)(89.25,80.5)
\multiput(83,77.25)(.06443299,.033505155){97}{\line(1,0){.06443299}}
%\end
\put(89.25,80.5){\line(0,1){7.25}}
\put(89.25,87.75){\line(-2,1){6.5}}
%\emline(82.75,91)(77.25,86.5)
\multiput(82.75,91)(-.041044776,-.03358209){134}{\line(-1,0){.041044776}}
%\end
\put(77.25,86.5){\line(0,-1){7}}
\put(89.25,87.75){\line(1,0){8.25}}
%\emline(97.5,87.75)(101,94.25)
\multiput(97.5,87.75)(.033653846,.0625){104}{\line(0,1){.0625}}
%\end
\put(101,94.25){\line(1,0){7.25}}
%\emline(108.25,94.25)(111.25,88.5)
\multiput(108.25,94.25)(.03370787,-.06460674){89}{\line(0,-1){.06460674}}
%\end
%\emline(111.25,88.5)(108,82.5)
\multiput(111.25,88.5)(-.033505155,-.06185567){97}{\line(0,-1){.06185567}}
%\end
\put(108,82.5){\line(-1,0){6.5}}
%\emline(98,88)(101.75,82.5)
\multiput(98,88)(.033482143,-.049107143){112}{\line(0,-1){.049107143}}
%\end
%\dashline{1}(18.25,86.5)(24.75,90)
\multiput(18.18,86.43)(.0625,.0336538){13}{\line(1,0){.0625}}
\multiput(19.805,87.305)(.0625,.0336538){13}{\line(1,0){.0625}}
\multiput(21.43,88.18)(.0625,.0336538){13}{\line(1,0){.0625}}
\multiput(23.055,89.055)(.0625,.0336538){13}{\line(1,0){.0625}}
%\end
%\dashline{1}(17.75,83.25)(28,89.5)
\multiput(17.68,83.18)(.0525641,.0320513){15}{\line(1,0){.0525641}}
\multiput(19.257,84.141)(.0525641,.0320513){15}{\line(1,0){.0525641}}
\multiput(20.834,85.103)(.0525641,.0320513){15}{\line(1,0){.0525641}}
\multiput(22.41,86.064)(.0525641,.0320513){15}{\line(1,0){.0525641}}
\multiput(23.987,87.026)(.0525641,.0320513){15}{\line(1,0){.0525641}}
\multiput(25.564,87.987)(.0525641,.0320513){15}{\line(1,0){.0525641}}
\multiput(27.141,88.949)(.0525641,.0320513){15}{\line(1,0){.0525641}}
%\end
%\dashline{1}(19.25,81)(29.25,87.25)
\multiput(19.18,80.93)(.0512821,.0320513){15}{\line(1,0){.0512821}}
\multiput(20.718,81.891)(.0512821,.0320513){15}{\line(1,0){.0512821}}
\multiput(22.257,82.853)(.0512821,.0320513){15}{\line(1,0){.0512821}}
\multiput(23.795,83.814)(.0512821,.0320513){15}{\line(1,0){.0512821}}
\multiput(25.334,84.776)(.0512821,.0320513){15}{\line(1,0){.0512821}}
\multiput(26.872,85.737)(.0512821,.0320513){15}{\line(1,0){.0512821}}
\multiput(28.41,86.699)(.0512821,.0320513){15}{\line(1,0){.0512821}}
%\end
%\dashline{1}(20.25,79)(30,84.75)
\multiput(20.18,78.93)(.0541667,.0319444){15}{\line(1,0){.0541667}}
\multiput(21.805,79.888)(.0541667,.0319444){15}{\line(1,0){.0541667}}
\multiput(23.43,80.846)(.0541667,.0319444){15}{\line(1,0){.0541667}}
\multiput(25.055,81.805)(.0541667,.0319444){15}{\line(1,0){.0541667}}
\multiput(26.68,82.763)(.0541667,.0319444){15}{\line(1,0){.0541667}}
\multiput(28.305,83.721)(.0541667,.0319444){15}{\line(1,0){.0541667}}
%\end
%\dashline{1}(24.25,78.75)(29.25,81.75)
\multiput(24.18,78.68)(.0549451,.032967){13}{\line(1,0){.0549451}}
\multiput(25.608,79.537)(.0549451,.032967){13}{\line(1,0){.0549451}}
\multiput(27.037,80.394)(.0549451,.032967){13}{\line(1,0){.0549451}}
\multiput(28.465,81.251)(.0549451,.032967){13}{\line(1,0){.0549451}}
%\end
%\dashline{1}(31.75,85.5)(37.5,89.75)
\multiput(31.68,85.43)(.0449219,.0332031){16}{\line(1,0){.0449219}}
\multiput(33.117,86.492)(.0449219,.0332031){16}{\line(1,0){.0449219}}
\multiput(34.555,87.555)(.0449219,.0332031){16}{\line(1,0){.0449219}}
\multiput(35.992,88.617)(.0449219,.0332031){16}{\line(1,0){.0449219}}
%\end
%\dashline{1}(32.25,82.25)(32,82)
\multiput(32.18,82.18)(-.03125,-.03125){4}{\line(0,-1){.03125}}
%\end
%\dashline{1}(31.75,82.75)(41,89.5)
\multiput(31.68,82.68)(.0453431,.0330882){17}{\line(1,0){.0453431}}
\multiput(33.221,83.805)(.0453431,.0330882){17}{\line(1,0){.0453431}}
\multiput(34.763,84.93)(.0453431,.0330882){17}{\line(1,0){.0453431}}
\multiput(36.305,86.055)(.0453431,.0330882){17}{\line(1,0){.0453431}}
\multiput(37.846,87.18)(.0453431,.0330882){17}{\line(1,0){.0453431}}
\multiput(39.388,88.305)(.0453431,.0330882){17}{\line(1,0){.0453431}}
%\end
%\dashline{1}(40.75,89.25)(41.25,89.75)
\multiput(40.68,89.18)(.03125,.03125){8}{\line(0,1){.03125}}
%\end
%\dashline{1}(37.5,89.5)(38.25,90)
\multiput(37.43,89.43)(.046875,.03125){8}{\line(1,0){.046875}}
%\end
%\dashline{1}(32.75,80.75)(42,87.5)
\multiput(32.68,80.68)(.0453431,.0330882){17}{\line(1,0){.0453431}}
\multiput(34.221,81.805)(.0453431,.0330882){17}{\line(1,0){.0453431}}
\multiput(35.763,82.93)(.0453431,.0330882){17}{\line(1,0){.0453431}}
\multiput(37.305,84.055)(.0453431,.0330882){17}{\line(1,0){.0453431}}
\multiput(38.846,85.18)(.0453431,.0330882){17}{\line(1,0){.0453431}}
\multiput(40.388,86.305)(.0453431,.0330882){17}{\line(1,0){.0453431}}
%\end
%\dashline{1}(34,78.5)(42.75,85.25)
\multiput(33.93,78.43)(.0428922,.0330882){17}{\line(1,0){.0428922}}
\multiput(35.388,79.555)(.0428922,.0330882){17}{\line(1,0){.0428922}}
\multiput(36.846,80.68)(.0428922,.0330882){17}{\line(1,0){.0428922}}
\multiput(38.305,81.805)(.0428922,.0330882){17}{\line(1,0){.0428922}}
\multiput(39.763,82.93)(.0428922,.0330882){17}{\line(1,0){.0428922}}
\multiput(41.221,84.055)(.0428922,.0330882){17}{\line(1,0){.0428922}}
%\end
%\dashline{1}(42.5,84.75)(43,85.25)
\multiput(42.43,84.68)(.03125,.03125){8}{\line(0,1){.03125}}
%\end
%\dashline{1}(41.75,87)(42,87.5)
\multiput(41.68,86.93)(.03125,.0625){4}{\line(0,1){.0625}}
%\end
%\dashline{1}(38.25,78.5)(42.75,82)
\multiput(38.18,78.43)(.0428571,.0333333){15}{\line(1,0){.0428571}}
\multiput(39.465,79.43)(.0428571,.0333333){15}{\line(1,0){.0428571}}
\multiput(40.751,80.43)(.0428571,.0333333){15}{\line(1,0){.0428571}}
\multiput(42.037,81.43)(.0428571,.0333333){15}{\line(1,0){.0428571}}
%\end
%\dashline{1}(65.25,97)(71.5,101.25)
\multiput(65.18,96.93)(.0496032,.0337302){14}{\line(1,0){.0496032}}
\multiput(66.569,97.874)(.0496032,.0337302){14}{\line(1,0){.0496032}}
\multiput(67.957,98.819)(.0496032,.0337302){14}{\line(1,0){.0496032}}
\multiput(69.346,99.763)(.0496032,.0337302){14}{\line(1,0){.0496032}}
\multiput(70.735,100.707)(.0496032,.0337302){14}{\line(1,0){.0496032}}
%\end
%\dashline{1}(65.25,94)(75.25,101.75)
\multiput(65.18,93.93)(.0420168,.032563){17}{\line(1,0){.0420168}}
\multiput(66.608,95.037)(.0420168,.032563){17}{\line(1,0){.0420168}}
\multiput(68.037,96.144)(.0420168,.032563){17}{\line(1,0){.0420168}}
\multiput(69.465,97.251)(.0420168,.032563){17}{\line(1,0){.0420168}}
\multiput(70.894,98.358)(.0420168,.032563){17}{\line(1,0){.0420168}}
\multiput(72.323,99.465)(.0420168,.032563){17}{\line(1,0){.0420168}}
\multiput(73.751,100.573)(.0420168,.032563){17}{\line(1,0){.0420168}}
%\end
%\dashline{1}(66,91.5)(76.25,100)
\multiput(65.93,91.43)(.0401961,.0333333){17}{\line(1,0){.0401961}}
\multiput(67.296,92.563)(.0401961,.0333333){17}{\line(1,0){.0401961}}
\multiput(68.663,93.696)(.0401961,.0333333){17}{\line(1,0){.0401961}}
\multiput(70.03,94.83)(.0401961,.0333333){17}{\line(1,0){.0401961}}
\multiput(71.396,95.963)(.0401961,.0333333){17}{\line(1,0){.0401961}}
\multiput(72.763,97.096)(.0401961,.0333333){17}{\line(1,0){.0401961}}
\multiput(74.13,98.23)(.0401961,.0333333){17}{\line(1,0){.0401961}}
\multiput(75.496,99.363)(.0401961,.0333333){17}{\line(1,0){.0401961}}
%\end
%\dashline{1}(67,89.25)(77,97.75)
\multiput(66.93,89.18)(.0392157,.0333333){17}{\line(1,0){.0392157}}
\multiput(68.263,90.313)(.0392157,.0333333){17}{\line(1,0){.0392157}}
\multiput(69.596,91.446)(.0392157,.0333333){17}{\line(1,0){.0392157}}
\multiput(70.93,92.58)(.0392157,.0333333){17}{\line(1,0){.0392157}}
\multiput(72.263,93.713)(.0392157,.0333333){17}{\line(1,0){.0392157}}
\multiput(73.596,94.846)(.0392157,.0333333){17}{\line(1,0){.0392157}}
\multiput(74.93,95.98)(.0392157,.0333333){17}{\line(1,0){.0392157}}
\multiput(76.263,97.113)(.0392157,.0333333){17}{\line(1,0){.0392157}}
%\end
%\dashline{1}(70.5,89.25)(78,95.5)
\multiput(70.43,89.18)(.040107,.0334225){17}{\line(1,0){.040107}}
\multiput(71.793,90.316)(.040107,.0334225){17}{\line(1,0){.040107}}
\multiput(73.157,91.452)(.040107,.0334225){17}{\line(1,0){.040107}}
\multiput(74.521,92.589)(.040107,.0334225){17}{\line(1,0){.040107}}
\multiput(75.884,93.725)(.040107,.0334225){17}{\line(1,0){.040107}}
\multiput(77.248,94.862)(.040107,.0334225){17}{\line(1,0){.040107}}
%\end
%\dashline{1}(77.25,84.5)(84.5,90.25)
\multiput(77.18,84.43)(.0402778,.0319444){18}{\line(1,0){.0402778}}
\multiput(78.63,85.58)(.0402778,.0319444){18}{\line(1,0){.0402778}}
\multiput(80.08,86.73)(.0402778,.0319444){18}{\line(1,0){.0402778}}
\multiput(81.53,87.88)(.0402778,.0319444){18}{\line(1,0){.0402778}}
\multiput(82.98,89.03)(.0402778,.0319444){18}{\line(1,0){.0402778}}
%\end
%\dashline{1}(77.5,82.25)(86.5,89.25)
\multiput(77.43,82.18)(.0432692,.0336538){16}{\line(1,0){.0432692}}
\multiput(78.814,83.257)(.0432692,.0336538){16}{\line(1,0){.0432692}}
\multiput(80.199,84.334)(.0432692,.0336538){16}{\line(1,0){.0432692}}
\multiput(81.584,85.41)(.0432692,.0336538){16}{\line(1,0){.0432692}}
\multiput(82.968,86.487)(.0432692,.0336538){16}{\line(1,0){.0432692}}
\multiput(84.353,87.564)(.0432692,.0336538){16}{\line(1,0){.0432692}}
\multiput(85.737,88.641)(.0432692,.0336538){16}{\line(1,0){.0432692}}
%\end
%\dashline{1}(78,79.5)(88.75,88)
\multiput(77.93,79.43)(.0421569,.0333333){17}{\line(1,0){.0421569}}
\multiput(79.363,80.563)(.0421569,.0333333){17}{\line(1,0){.0421569}}
\multiput(80.796,81.696)(.0421569,.0333333){17}{\line(1,0){.0421569}}
\multiput(82.23,82.83)(.0421569,.0333333){17}{\line(1,0){.0421569}}
\multiput(83.663,83.963)(.0421569,.0333333){17}{\line(1,0){.0421569}}
\multiput(85.096,85.096)(.0421569,.0333333){17}{\line(1,0){.0421569}}
\multiput(86.53,86.23)(.0421569,.0333333){17}{\line(1,0){.0421569}}
\multiput(87.963,87.363)(.0421569,.0333333){17}{\line(1,0){.0421569}}
%\end
%\dashline{1}(80.25,78.5)(89,85.25)
\multiput(80.18,78.43)(.0428922,.0330882){17}{\line(1,0){.0428922}}
\multiput(81.638,79.555)(.0428922,.0330882){17}{\line(1,0){.0428922}}
\multiput(83.096,80.68)(.0428922,.0330882){17}{\line(1,0){.0428922}}
\multiput(84.555,81.805)(.0428922,.0330882){17}{\line(1,0){.0428922}}
\multiput(86.013,82.93)(.0428922,.0330882){17}{\line(1,0){.0428922}}
\multiput(87.471,84.055)(.0428922,.0330882){17}{\line(1,0){.0428922}}
%\end
%\dashline{1}(83.25,77.75)(88.75,82.25)
\multiput(83.18,77.68)(.0404412,.0330882){17}{\line(1,0){.0404412}}
\multiput(84.555,78.805)(.0404412,.0330882){17}{\line(1,0){.0404412}}
\multiput(85.93,79.93)(.0404412,.0330882){17}{\line(1,0){.0404412}}
\multiput(87.305,81.055)(.0404412,.0330882){17}{\line(1,0){.0404412}}
%\end
%\dashline{1}(98.5,89)(104.75,94)
\multiput(98.43,88.93)(.0416667,.0333333){15}{\line(1,0){.0416667}}
\multiput(99.68,89.93)(.0416667,.0333333){15}{\line(1,0){.0416667}}
\multiput(100.93,90.93)(.0416667,.0333333){15}{\line(1,0){.0416667}}
\multiput(102.18,91.93)(.0416667,.0333333){15}{\line(1,0){.0416667}}
\multiput(103.43,92.93)(.0416667,.0333333){15}{\line(1,0){.0416667}}
%\end
%\dashline{1}(99.25,86.75)(108.5,94)
\multiput(99.18,86.68)(.0418552,.0328054){17}{\line(1,0){.0418552}}
\multiput(100.603,87.795)(.0418552,.0328054){17}{\line(1,0){.0418552}}
\multiput(102.026,88.91)(.0418552,.0328054){17}{\line(1,0){.0418552}}
\multiput(103.449,90.026)(.0418552,.0328054){17}{\line(1,0){.0418552}}
\multiput(104.872,91.141)(.0418552,.0328054){17}{\line(1,0){.0418552}}
\multiput(106.295,92.257)(.0418552,.0328054){17}{\line(1,0){.0418552}}
\multiput(107.718,93.372)(.0418552,.0328054){17}{\line(1,0){.0418552}}
%\end
%\dashline{1}(100.5,84.75)(109.75,92)
\multiput(100.43,84.68)(.0418552,.0328054){17}{\line(1,0){.0418552}}
\multiput(101.853,85.795)(.0418552,.0328054){17}{\line(1,0){.0418552}}
\multiput(103.276,86.91)(.0418552,.0328054){17}{\line(1,0){.0418552}}
\multiput(104.699,88.026)(.0418552,.0328054){17}{\line(1,0){.0418552}}
\multiput(106.122,89.141)(.0418552,.0328054){17}{\line(1,0){.0418552}}
\multiput(107.545,90.257)(.0418552,.0328054){17}{\line(1,0){.0418552}}
\multiput(108.968,91.372)(.0418552,.0328054){17}{\line(1,0){.0418552}}
%\end
%\dashline{1}(102.25,82.75)(110.5,90)
\multiput(102.18,82.68)(.0381944,.0335648){18}{\line(1,0){.0381944}}
\multiput(103.555,83.888)(.0381944,.0335648){18}{\line(1,0){.0381944}}
\multiput(104.93,85.096)(.0381944,.0335648){18}{\line(1,0){.0381944}}
\multiput(106.305,86.305)(.0381944,.0335648){18}{\line(1,0){.0381944}}
\multiput(107.68,87.513)(.0381944,.0335648){18}{\line(1,0){.0381944}}
\multiput(109.055,88.721)(.0381944,.0335648){18}{\line(1,0){.0381944}}
%\end
%\dashline{1}(105.75,82.75)(109.75,86)
\multiput(105.68,82.68)(.0408163,.0331633){14}{\line(1,0){.0408163}}
\multiput(106.823,83.608)(.0408163,.0331633){14}{\line(1,0){.0408163}}
\multiput(107.965,84.537)(.0408163,.0331633){14}{\line(1,0){.0408163}}
\multiput(109.108,85.465)(.0408163,.0331633){14}{\line(1,0){.0408163}}
%\end
\put(49.25,81.25){$p_1$}
\put(61.25,89){$p_2$}
\put(69.5,78.5){$p_3$}
\put(92.25,85.75){$p_4$}
\end{picture}
\end{center}
\caption{Faces and bridges in a diagram over $G$}\label{bri}
\end{figure}

\begin{lemma} \label{type} Let $n\ge d\ge 12$ and $\sigma=16(\log_2 d)/d$.
Then there exists less than $2^{\sigma n}$
types of reduced diagrams $\Delta$ over $G=\la X\mid x^d\ra $ such that the boundary path $q$ of $\Delta$ is a product of at most $3$ reduced subpaths and $|q|=n$.
\end{lemma}
\proof Every diagram $\Delta$ under consideration is constructed from several
faces $\Pi_1,\Pi_2,\dots$
and bridges. By definition, a {\it bridge} is a  subpath $p$ of positive length in $q$ (more accurately, it is the pair $(p,p^{-1})$)
such that $p^{-1}$ is also a subpath of the same path $q$.
(There are $5$ faces and $4$ maximal bridges $p_1$, $p_2$, $p_3$, and $p_4$ in Fig. \ref{bri}.)

To code all possible maps up to isomorphism, we first enumerate the consecutive edges of $q$ as
$e_1,\dots,e_n.$ Then we move along $q$ and place pairs of brackets as follows. If our next edge
$e_i$ is the first edge of a maximal bridge $p$, then we place a left round bracket right before $e_i$ and a right
round bracket after the last edge of $p^{-1}.$ Similarly, if $e_i$ is the first edge of the boundary
of a face $\Pi_j$ (i.e., $e_1,\dots,e_{i-1}$ do not belong to $\partial \Pi_j$), then we place
a left square bracket before $e_i$ and the right square bracket after the last edge of $q$ belonging to $\partial \Pi_j$.

The brackets restore the map up to isomorphism. Indeed, if $q$ has a vertex of degree one different of $q_-=q_+$,
then it is the end of a bridge $p=e_i\dots e_{i+k-1}$ , and so
we should have a pair of round brackets  $\dots(e_i\dots e_{i+k-1}e_{i+k}\dots e_{i+2k-1})\dots$ and no brackets between them. Such an arrangement restores the bridge $p$.
If we remove $p$ from $\Delta$ and $(pp^{-1})$ from its boundary path, then we decrease the number of bridges
and restore the whole map by induction. If there are no such vertices of degree 1, then
since the map is simply connected and there are no inner edges in $\Delta$,
there should be a face $\Pi_j$ attached to the remaining part of the map at one vertex only. Hence the corresponding to $\Pi_j$ part of the boundary
has the form $[e_i\dots e_{i+d-1}]$ without brackets in the middle. So this pair of brackets defines
the (unlabeled) boundary of
some face.
If we remove
$\Pi_j$ from $\Delta$ and $[e_i\dots e_{i+d-1}]$ from its boundary, then we decrease the number of faces, and the whole map $\Delta$ restores by induction (but this moment, only up to the base points of the faces).

     The number of faces in $\Delta$
     is at most $n/d$. If there are no vertices of degree $1$ in $\Delta$ except for $q_-$, the number of maximal bridges is less than twice the number of the faces because $\Delta$ is a simply connected map. There can be at most one or two additional vertices of degree $1$ since $q$ is a product of at most $3$ reduced paths. But these vertices can increase the number of bridges by at most $4$. Thus we have at most  $(2n+3d)/d\le 5n/d$ pairs of round brackets placed in the sequence $1,\dots, n$.

 If the number of pairs of square brackets is $r$ then the number of the symbols after bracketing
 is $n+2r$, and so to obtain all the brackets, one declares $r$ symbols as left square brackets and $r$
 symbols as right ones. Hence the number of square bracket arrangements  is less than
 $$\frac{(n+2r)!}{(n+r)!r!}\times\frac{(n+r)!}{n!r!}=\frac{(n+2r)!}{n!(r!)^2}$$
 Similarly, if the number of pairs of round brackets is $s$, then the number of their different
 arrangements is less than $\frac{(n+2r+2s)!}{(n+2r)!(s!)^2}$. The product of these two fractions
 is $\frac{(n+2r+2s)!}{n!(r!)^2 (s!)^2},$ which is less than  $\frac{(2n)^{2r+2s}}{(r!)^2 (s!)^2},$
 since for $n\ge d\ge 12$, we have $2r+2s \le \frac{2n}{d}+\frac{2(2n+3d)}{d}\le n$. In turn,
 \begin{equation}\label{nrs}
 \frac{(2n)^{2r+2s}}{(r!)^2 (s!)^2}<\left(\frac{2ne}{r}\right)^{2r}\times\left(\frac{2ne}{s}\right)^{2s}
 \end{equation}
 since $r! >\left(\frac{r}{e}\right)^r$
 and $s! >\left(\frac{s}{e}\right)^s$ by well-known Stirling's inequality. Our next estimate of the product (\ref{nrs}) uses that the functions $\left(\frac{a}{x}\right)^x$ increase  on the intervals $(0; a/e).$ So replacing $r$ by $n/d$ and $s$ by $5n/d$ at the right-hand side of (\ref{nrs}), we obtain the upper bound $(2ed)^{2n/d}(2ed/5)^{10n/d}<90d^{12n/d}$. Therefore for fixed $r$ and $s,$ we have that the number of bracket arrangements $L(r,s)$ is less than $90d^{12n/d}$.
 Hence
 \begin{equation}\label{nd}
 \sum_{r\le n/d, s\le 5n/d} L(r,s)<90(12(n/d)^2)d^{12n/d}< d^{15n/d}=2^{15(\log_2 d) n/d}
\end{equation}
 There are at most $d$ ways to chose the base point of a face, i.e., at most $d^{n/d}=2^{nd^{-1}\log_2 d}$
 ways for all the faces. Taking (\ref{nd}) into account, we get less than $2^{16(\log_2d)n/d}$ possible types, as required.
 \endproof

Now we need a condition depending on three parameters but helpful for
the further argument since it is adjusted to $H$-diagrams and survives
after the inductive steps in the proof of Theorem \ref{base}.

Let $H$ be a subgroup of the free group $F=F(X)=F(x_1,\dots,x_m)$.
\begin{defn}
We say that the property $P(d,\mu,\rho)$
holds for $H$, where $d\ge 12$ is a positive integer and $\mu,\rho$ are positive numbers, if

(1) $H$ contains $x^d$ for some $x\in X$;

(2) for any $n\ge 1$, any reduced words $u,v\in F$ and any letter $y\in X^{\pm 1}$, there are at most $2^{\rho n}$ reduced words $w$ of length $n\ge 1$ over $X^{\pm 1}$
such that

(a) $|w|\ge \mu(|u|+|v|)$,

(b) the product $uwv$ belongs to $H$, and

(c) the first letter of $w$ is not $y$.
\end{defn}

\begin{lemma}\label{estim} If a subgroup $H\le F$ satisfies the condition $P(d,\mu,\rho),$ then
the condition $P(d,\mu',\rho')$ holds for the normal closure $N$ of $x^d$ in $H$,
where $\mu'=\mu^{1/2}$ and $$\rho'=\rho/2+\sigma(1+\mu^{-1/2})+(1+\mu^{-1/2})d^{-1}+2(\mu+\mu^{1/2})\log_2(2m-1)$$
with $\sigma$ defined in Lemma \ref{type}.
\end{lemma}
\proof We should verify the second part of $P(d,\mu',\rho')$ for $N$,
 and so the letter $y$ is now fixed. By Lemma \ref{vK} (2), for any words $u,v,$ and $w$ satisfying the  parts (2a) -- (2c) with the parameter $\mu'$, there is a reduced $H$-diagram $\Delta$ over the
presentation $G =\la X\mid x^d\ra$ whose boundary
label is $uwv$. If $|w|=n <d-|u|-|v|,$ then the perimeter of $\Delta$ is less than $d$, and so $\Delta$ has no faces at all. It follows that
$uwv=1$ in $H$, i.e., there is only one solution $w=u^{-1}v$, and we are done since $1<2^{\rho'}$.
Hence we may assume further that $n\ge d-|u|-|v|,$ and so $|\partial\Delta|\ge d$.

 The perimeter of $\Delta$ is at most $(1+(\mu')^{-1})n$ by (2a), and so the number of possible types
of such diagrams $\Delta$ is less than $2^{\sigma(1+(\mu')^{-1})n}$ by Lemma \ref{type}, because every $H$-diagram is a diagram over $G$. Below we consider an $H$-diagram of a fixed type with the reduced boundary path $q=q_1pq_2,$ where $q_1$ and
$q_2$ are labeled by $u$ and $v$, respectively, and $|p|=n.$ Our goal is to estimate the number
of possible labels for $p$ in all such $H$-diagrams of the given type.

Starting with a base point $o(\Pi)$ of a face $\Pi$ and going clockwise along $\partial\Pi$ one can read
either $x^d$ or $x^{-d}.$ Therefore
there are at most $2^{(1+(\mu')^{-1}) n/d}$ ways to label the boundaries of all the faces in $\Delta.$
It remains to estimate the number of possible ways to label the bridge edges of $\Delta$
belonging to $p$ under the assumption that the boundaries of the faces are already labeled
and their base points are fixed.

The unlabeled edges of $p$ belong to the union of maximal subpaths $p_0, p_1...,p_s$ of $q$, connecting different faces  or connecting the base point $o$ of $\Delta$ with a face. Note that $s\le 2f$, where $f$ is the number of faces in $\Delta$ since the diagram is simply-connected.   We will subsequently prescribe the labels to
the unlabel edges of $p_0,p_1,\dots$ and estimate the number of possibilities for $Lab(p_i)$ under
the condition that the paths
$p_0,...,p_{i-1}$ are completely labeled ($i\ge 0$). Some
of the edges of $p_i$ could be labeled at previous steps being edges of $q_1^{\pm 1}$  or $q_2^{\pm 1}$, or $p_j^{-1}$, where $j<i$. However since $\Delta$ is
simply-connected, we have a factorization $p_i=p'p''p'''$ (the lengths of some factors can be $0$), where every edge of $p'$ and $p'''$ has been labeled before the beginning of the $i$-th step and every edge of $p''$ is yet unlabeled at this stage. Now we consider two cases.

{\bf Case 1}: The path $p''$ is `long': $|p''|\ge \mu(|p'|+|p'''|+d)$
if $p_i$ connects two faces or $|p''|\ge \mu(|p'|+|p'''|+d/2)$ if it connects
the base point $o(\Delta)$ with a face.
If the vertex $(p_i)_-$ belongs to some face $\Pi$, then one can
find a path $z'$ of length $\le d/2$ connecting the base point $o_{\Pi}$ and
$(p_i)_-$. Similarly, the vertex $(p_i)_+$ can be connected with a base
point  $o_{\Pi'}$ by a path $z'''$ of length $\le d/2.$ Thus by the definition of
$H$-diagram, $\Lab(z'p')\Lab(p'')\Lab(p'''z''')\in H$, and the word $w_0\equiv \Lab(p'')$
satisfies the condition $ u_0w_0v_0\in H$, where $u_0$ and $v_0$ are the reduced forms of
the words $\Lab(z'p')$ and $\Lab(p'''z''')$, respectively, and so $|u_0|+|v_0|\ge |p'|+|p'''|-d$ and
$$|w_0|=|p''|\ge \mu(|p'|+|p'''|+d)\ge \mu(|u_0|+|v_0|)$$

Now, applying the condition $P(d,\mu,\rho)$ to $H$ with the triple $u_0,$ $ w_0,$ $ v_0$, one should
name a prohibited first letter $y_0$ of the word $w_0$. It is the letter $y$ distinguished above
if $p''$ is just the beginning of the path $p$. Otherwise, since the word $w$ should be reduced,
$y_0$ is the inverse letter to the last letter of the labeled (by this stage) beginning
of the path $p$.
So the number of possible
labels for $p''$ is less than $2^{\rho|p''|}$.
The same estimate works if $(p_i)_-$ or $(p_i)_+$ coincides with the base point of $\Delta$
since we have $|z'|=0$ or $|z'''|=0$ in these cases.

{\bf Case 2}: The path $p''$ is short (= not long). Then we just bound
the number of possible labels for $p''$ by the number of reduced words over $X^{\pm 1}$
with the restriction of the form $y''\ne y_0$ for the first letter $y''$ of $\Lab(p'')$.
Hence the upper bounds in the two subcases for the labels of $p''$, depending on whether $p_i$ connects
two faces or not, are
$(2m-1)^{\mu(|p'|+|p'''|+d)}$ and $(2m-1)^{\mu(|p'|+|p'''|+d/2)},$ respectively.

\medskip

The total length of all short paths $p''_i$ does not exceed $\mu(|p_0|+\dots+|p_s|+2fd)$
since $s\le 2f$. Hence it does not exceed
$2\mu|\partial\Delta|\le 2\mu(1+(\mu')^{-1})n.$ So the total number of choices of the labels we have
in all short cases does not exceed $(2m-1)^{2\mu(1+(\mu')^{-1})n}.$

The sum of length of all unlabeled $p''_i$-s over all long cases does not exceed
$|p|/2=n/2$ since if $e$ is an unlabeled edge of $p$, then $e^{-1}$ must also occur in
the path $p$ (but not in $p''_i$) and it is also unlabeled
for the $i$-th step; so this step decreases the number of unlabeled edges in $p$
by $2|p''_i|$. Therefore the total number of options over all long cases is less
than $\prod _{i=0}^s 2^{\rho |p''_i|}\le 2^{\rho n/2}.$

It remains to multiply the upper bounds obtained in the proof for the number of types, the number of
possible labels of faces, and the numbers of options
one has in long and short cases. As desired, this product is $2^{\rho'n}$, by the definition of $\rho'$.

\endproof

\begin{lemma} \label{param} Let $N=N(\ell,d)$ denote the $\ell$-subnormal closure of the power $x^d$ in a free
group $F$ of rank $m\ge 1,$ where $x$ is one of the free generators. Then for every $\varepsilon>0$ and $\mu'\in (0,\varepsilon]$, there
exists $d_0\ge 12$ such that $N$ satisfies the condition
$P(d,\mu', {2^{-\ell}}log_2(2m-1) +\varepsilon)$ provided $d\ge d_0.$
\end{lemma}

\proof To prove the lemma by induction on $\ell$, we set $N(0,d)=F.$ Then the statement of the lemma holds for
$\ell=0$ since the number of reduced words of length $n$ over $X^{\pm 1}$ with the restriction of the form $y\ne y_0$ for
the first letter is equal $(2m-1)^n = 2^{n\log_2(2m-1)}.$ Assume now that $\ell\ge 1$, the statement
holds for $\ell-1$, and we have arbitrary $\varepsilon>0$ and $\mu'\in (0,\varepsilon]$.

One can choose a number $\mu$ so that

$$0<\mu\le \min((\mu')^2,\varepsilon/2)\;\; and \;\;(\mu+\mu^{1/2})\log_2(2m-1)<\varepsilon/8.$$
Then one can choose $d_1$ so that for $\sigma=16(\log_2d_1)/d_1$, we have
$$\sigma<(1+\mu^{-1/2})^{-1}\varepsilon/8\;\; and
\;\;(1+\mu^{-1/2})d_1^{-1}<\varepsilon/8$$
By the inductive hypothesis one can choose $d_0\ge d_1$ such that the property $P(d, \mu, \rho)$ holds for the subgroup
$H=N(\ell-1,d)$ if
$$\rho = 2^{-\ell+1}log_2(2m-1) +\varepsilon\;\;  and \;\; d\ge d_0$$
Now by Lemma \ref{estim},
the property $P(d, \mu^{1/2}, \rho')$ holds for $N=N(\ell,d)$ if $d\ge d_0$ and
$$\rho'=\rho/2 + \varepsilon/8 +\varepsilon/8+\varepsilon/4 ={2^{-\ell}}\log_2(2m-1)+\varepsilon$$
Since this property allows us to increase
the second parameter, we obtain $P(d, \mu', \rho')$ for $N(\ell,d)$, as required.
\endproof

{\bf Proof of Theorem \ref{base}}. Below we keep in mind Theorem \ref{rate}.

{\bf Part} (1). It suffices to prove that there are constants $c>0$ and $\psi> 2^{-l}$ such that the relative growth function $g_N$ of $N$ in $F$ satisfies the inequality $g_N(n)> c(2m-1)^{\psi n}$ for every $n\ge 0$. For $\ell=1$ this assertion follows from Lemma \ref{gri}. Then we induct on $\ell$
and assume that $\ell\ge 2$. The subgroup $N$ is a normal subgroup of an $\ell-1$-subnormal subgroup $H$ of $F$.
It follows from the inductive hypothesis that $g_H(n)> b(2m-1)^{\theta n}$ for every $n\ge 0$ and some
constants $\theta>2^{-\ell+1}$ and $b>0.$

Let $u$ be a non trivial element of $N$ and $|u|=a$ for some $a\ge 1$. Assume that $vuv^{-1}=v'u(v')^{-1}$, where
$v,v'\in H$ and $|v|, |v'|\le n$ for some $n$. It follows that $v^{-1}v'$ belongs to the centralizer of $u,$ which is a cyclic subgroup
$\la u_0 \ra $ of $F$ (see \cite{LS}, I.2.19).  The cyclic subgroup $\la u_0 \ra $ has at most $4n+1$ elements
of length $\le 2n$. Therefore for any $v$ of length $\le n$, there are at most $4n+1$ distinct values of $v'$
with $|v'|\le n$ and  $vuv^{-1}=v'u(v')^{-1}$. It follows that the number of distinct elements of
$N$ of the form $vuv^{-1}$, where $|v|\le n$ is at least $g_H(n)/(4n+1),$ where $g_H$ is the relative growth function
of $H$ in $F$. Since $|vuv^{-1}|\le 2n+a$, we conclude that $g_N(2n+a)\ge (4n+1)^{-1} g_H(n)$ for every $n\ge 0$. In other words, $$g_N(n)\ge (2n)^{-1}g_H(\lfloor(n-a)/2\rfloor)$$ Here the right hand side is greater
than $b(2m-1)^{\theta \lfloor(n-a)/2\rfloor-log_{2m-1}(2n)}$. The exponent is greater than $\psi n$ for any $\psi \in (2^{-l},\theta/2)$ and all sufficiently large values of $n$. Hence there is a constant $c=c(\psi)>0$ such that
$g_N(n) >c(2m-1)^{\psi n}$ for every $n\ge 0$,
and the assertion (1) is proved.
\medskip

\begin{rem} The above argument implies that $\alpha_N\ge\alpha_H^{1/2}$ for any nontrivial normal
subgroup $N$ of any subgroup $H\le F_m$.
\end{rem}

{\bf Part} (2). One may assume that $m\ge 2$. Note that the restriction (2a) is empty in the condition  $P(d, \mu, \rho)$
if the words $u$ and $v$ are empty. Hence by Lemma \ref{param} (with $\mu'=\varepsilon$), for arbitrary $\ell\ge 1$ and $\varepsilon>0$, we have a nontrivial $\ell$-subnormal subgroup $N(\ell,d)$ of $F$ such that the number of reduced words of length $n$ in it with the restriction
of the form $y\ne y_0$ on the first letter does not exceed $(2m-1)^{(2^{-\ell}+\varepsilon)n}.$ Hence the number of all reduced
words of length $\le n$ in this subgroup is at most $2(n+1)(2m-1)^{(2^{-\ell}+\varepsilon)n}.$ The (upper) limit of the $n$-th
roots of these values does not exceed $(2m-1)^{(2^{-\ell}+\varepsilon)}$, and the theorem is proved.
$\Box$

\section{Transversals to subnormal subgroups in $F_m$; part (1) of Theorem \ref{cogr}.}\label{cos}

Let $H$ be a subgroup of the free group $F=F(x,y,\dots)$ and $x\in H$.
We denote by $T$ the set of all reduced words $w$ in the generators $\{x^{\pm 1}, y^{\pm 1}, \dots\}$ satisfying the following condition:
     \begin{Condition}\label{Cond1}
  If $v$ is a prefix of $w$ and $v\in H$, then $vx^{\pm 1}$ is not a prefix of $w$.
     \end{Condition}

Equivalently: If $vx^{\pm 1}$ is a prefix of $w$ then $v\not\in H$.

\bigskip

   Let $Y= T\cap H$. An easy consequence of Condition \ref{Cond1} is as follows

   \begin{Condition}\label{Cond2}
     A word from $T$ (respectively, from $Y$) does not start (resp., neither starts nor ends) with $x^{\pm 1}$.
     \end{Condition}

\begin{lemma}\label{lone} Let $F$, $H$, $T$ and $Y$ be as above. Denote by $N$ the normal closure of $x$ in $H$. Then the following are true:
\begin{enumerate}\item[\rm(a)] the set $\{wxw^{-1}|\:w\in Y\}$ is the free basis of $N$;

\item[\rm(b)] every nontrivial element $h\in N$ equals the reduced  word of the form
\begin{equation}\label{red}
u_0x^{k_1}u_1x^{k_2}\cdots x^{k_{s-1}}u_{s-1}x^{k_s}u_{s},
\end{equation}
 where $s\ge 1$,
     $u_0,\ldots,u_s\in H$, with $u_0\cdots u_s=1$ in $F$,  the words $u_1,\dots, u_{s-1}$ are all non-empty, the words $u_0,\dots,u_s$ neither start nor end with $x^{\pm 1}$, and the exponents $k_1,\dots, k_s$ are all nonzero;

    \item[\rm(c)] the reduced form of an element $h\in N\backslash\{ 1, x, x^{-1}\}$, viewed as an element of $F$, can be written as $v_1v_2v_3$,
     where $v_1$ and $v_3$ are some nonempty words representing the elements of $H$;

     \item[\rm(d)] the set $T$ is a right Schreier transversal of the subgroup $N$ in $F$; moreover, each $t\in T$ is a shortest element in the coset $Nt$.
     \end{enumerate}
\end{lemma}

     \proof (a) By definition, $N$ is generated by the conjugates $wxw^{-1}$, where $w\in H$. Let us apply induction on the  length of $w$ to prove that $N$ is generated by the elements $wxw^{-1}$ with $w\in Y$. If the reduced form of $w$ has no prefixes of the form $vx^{\pm 1}$ with $v\in H$, then $w\in Y,$ and we are done. Otherwise,  assume the reduced form of $w$ equals $vx^{\pm 1}u$, where $v\in H$. Then we have in $F$:
     $$wxw^{-1} = vx^{\pm 1}u x u^{-1}x^{\mp 1}v^{-1}= (vx^{\pm 1}v^{-1})((vu) x (vu)^{-1})(vx^{\mp 1}v^{-1}).$$
Now the reduced forms of $v$ and $vu$ are shorter than $w$. By the assumption, $v\in H$. Since  each of the words $vx^{\pm 1}u$, $v$, and $x$ represents an element of $H$, we also have $vu\in H$. Applying induction to  $vx^{\pm 1}v^{-1}$ and $(vu) x (vu)^{-1}$ proves that the set $\{wxw^{-1}|\:w\in Y\}$ generates $N$.

     Now assume that we have a nontrivial relation between the elements $wxw^{-1}$, where $w\in Y$:
     \begin{equation}\label{rel}
     w_1x^{k_1} w_1^{-1}\cdots w_sx^{k_s}w_s^{-1}=1
     \end{equation}
     in $F$, where the exponents $k_1,\dots,k_s$ are nonzero
     and $w_i\ne w_{i+1}$ ($i=1,\dots,s-1$). If for some $i$, none of the
     two factors completely cancels in $w_i^{-1}w_{i+1}$, then the reduced form of $w_i^{-1}w_{i+1}$ will not start/end with $x^{\pm 1}$ by Condition \ref{Cond2}.
      Then without loss of generality, we may assume that $w_i$ is the prefix of $w_{i+1}$. But $w_ix^{\pm 1}$ is not a prefix of $w_{i+1}$ by Condition \ref{Cond1}, and so the non-empty reduced form of $w_i^{-1} w_{i+1}$ does not start (nor end) with $x^{\pm 1}$. Hence in any case, the factors $x^{k_i}$ will not be touched by the cancelations. So the left hand side of (\ref{rel}) is a nonempty word, a contradiction.  Thus, the proof of statement (a) is complete.

     (b) This claim is a direct consequence of the proof of (a) because all the factors $x^{k_i}$ survive in (\ref{rel}), each $w_i^{-1}w_{i+1}$ ($i=1,\dots s-1$) reduces to a nonempty word $u_i$ representing an element of $H$ and having no $x^{\pm 1}$ as a prefix or a suffix, and $u_0u_1\dots u_s=1 $ in $F$ for $u_0=w_1$ and $u_s=w_s^{-1}$.

      (c) If the reduced form of $h$ equals (\ref{red}), then we set $v_1=u_0$ if this is a nonempty word, and $v_1$ is the first letter of $x^{k_1}$ otherwise. Then set $v_3=u_s$ if $u_s$ is nonempty. If empty,  then using $h\neq x^{\pm 1}$ allows us to set $v_3$ equal to the last letter of $x^{k_s}$, without overlapping with $v_1$.

(d) We first show that for any $u\in F,$ there is $t\in T$ with $|t|\le|u|$ such that $Nu=Nt$. If $u\in T$, there is nothing to prove. Otherwise  by Condition \ref{Cond1}, $u=vx^{\pm 1}v'$ where $v\in H$. In this case, $u$ is equal in the free  group $F$ to $(vx^{\pm 1}v^{-1})vv'$, where $(vx^{\pm 1}v^{-1})\in N$. It follows that $Nu = Nvv'$ with $|vv'|<|u|$. Applying induction, we find $t\in T$ with $|t|\le|vv'|<|u|$ such that $Nu=Nvv'=Nt$.

     Obviously, the set $T$ is closed under prefixes, and  it only remains to show that $t_1t_2^{-1}\notin N$
     for different $t_1, t_2\in T$, where there are no cancelations in the product $t_1t_2^{-1}$.
     But if a reduced nontrivial word $t_1t_2^{-1}$ belongs to $N$, then by (b),
     either the word $t_1$ has to start with $u_0x^{\pm 1}$
     or $t_2$ has to start with $u_s^{-1} x^{\pm 1}$. This contradiction with Condition \ref{Cond1} completes the proof.

     \endproof

     \begin{lemma}\label{maxgr} Let $H$ be a subgroup of the free group $F=F_m=F(x_1,\dots,x_m)$, $m\ge 2$, and
     let $H$ contain the generator $x=x_m$. We denote by $N$ the normal closure of $x$ in $H$.
     The number of reduced words of length $n\ge 0$ in $F,$ which
     belong to $H$ can be written as $a_n(2m-1)^n$ for some real numbers $a_n$. If the series
     $\sum_{n=0}^{\infty} a_n$ converges, then the cogrowth of $N$ with respect to the
     generators $x_1,\dots,x_m $ is maximal.

     \end{lemma}

     \proof Note that $H$ is a proper subgroup of $F$ since otherwise $a_n\ge 1$ for every $n$.
     Therefore the growth function of the compliment $F\backslash H$ is $\Theta$-equivalent to the growth function
     of $F$, and so it is greater than  $c_1(2m-1)^n$ for some $c_1>0$ and every $n\ge 1$. Note that the number of all words
     of length $\le \ell$ in $F$ is $\sum_{i=0}^{\ell} 2m(2m-1)^i\le 2(2m-1)^{\ell}$.
Hence the number of reduced words in $F\backslash H$ of length $\le n$ without prefixes of lengths $\ge k$
     from $H$ is greater than
     $$ c_1(2m-1)^n - a_k(2m-1)^k\times  2(2m-1)^{n-k}-...- a_n(2m-1)^n\times 2>c_2(2m-1)^n$$
     for some $c_2>0$ if $k$ is chosen so that $2\sum_{i=k}^{\infty} a_i<c_1$ and $n\ge k$.
     So the growth of the subset of words $S\subset F\backslash H $ without prefixes from $H$ of
     lengths $\ge k$ for some fixed $k=k(H)$, is maximal.

     %, i.e. $\ge c_3(2m-1)^n$ for some $c_3>0$
     %and .

     Let us modify the set $S$ as follows. If a word $v\in S$ has a maximal prefix $u$
     from $H$ (of length $\le k$), i.e., $v=uv'$, then we count $v'$ to the set $S'$.
     The following properties of the mapping $v\mapsto v'$ are clear:
     \begin{enumerate}
    \item[(1)]$v'$ has no non-empty
     prefixes from $H$,
     \item[(2)] $|v'|\le |v|$,
     \item[(3)] every $v'$ has at most $c_3$ preimages $v$ in $S$,
     where $c_3$ is the number of reduced words of length $\le k$.
     \end{enumerate}
 Therefore the growth of $S'$ is also maximal with a constant $c_4\ge c_2/c_3>0.$

  It follows from Lemma \ref{lone} (c) that any inclusion $ g_1g_2^{-1}\in N$,
     where $g_1$ and $g_2$ are different reduced words, implies that either $g_1$ or $g_2$ has a non-empty
     prefix belonging to $H$. Hence different elements $g_1$ and $g_2$ from $S'$ cannot belong to the
     same right coset of $N$. Thus, the cogrowth function of $N$ with respect to the generators
     $x_1,\dots,x_m$ is at least $c_4(2m-1)^n$ for every sufficiently large $n$, as required.
     \endproof

     {\bf Proof of Theorem \ref{cogr} (1).}
    Let $x=x_m \in F=F(x_1,\dots,x_m) $ and $m\ge 3.$ Since the $\ell$-subnormal closure of $x$ is contained in the $2$-subnormal closure of $x$, we may
     assume that $\ell=2$.

     Let $H$ be the normal closure of $x$ in $F$ and $N$ be the normal closure of $x$ in $H$.
     The group $F/H$ is free of rank $m-1\ge 2$,
     in particular, $F/H$ is non-amenable \cite{G}. By Grigorchuk's amenability
     criterion \cite{Gr}, the e.g. rate of $H$ in $F$ with respect to $ x_1,\dots,x_m$ is less than $2m-1-\varepsilon$ for some
     $\varepsilon>0$. In other words, the number of reduced words of $F$ belonging to
     $H$ and having length $n$ is less that $a_n(2m-1)^n,$ where $a_n=o((2m-1)^{-\varepsilon n})$.
     Since the geometric progression series $\sum_{n=0}^{\infty} (2m-1)^{-\varepsilon n}$ converges, the cogrowth of $N$ is maximal by
     Lemma \ref{maxgr}. $\Box $

\begin{rem}
By Theorem \ref{intersec},
$\lim_{\ell\to\infty} D(g,\ell)=\infty$  for every $g\in F_m$.  But if $m\ge 2$ and  $x$ is
a free generator of $F_m$, then we have a simple explicit formula. It follows from Lemma
\ref{lone} (b), that $D(x,\ell)\ge 2D(x,\ell-1)+1$, because
for
every element $h\in N\backslash\langle x\rangle$, where $N=\langle x\rangle_{\ell}^{F_m}$,
we
should have at least two $u$-factors from $H=\langle x\rangle_{\ell-1}^{F_m}$ and at
least one letter $x^{\pm 1}$ in (\ref{red}). For the opposite inequality, one can take $h=u_0xu_0^{-1}$
with $u_0\in H\backslash \langle x \rangle$.
Thus by induction, $D(x,\ell) = 2^{\ell+1}-1$.

\end{rem}

\section{Random walks associated with $F(x,y)$.}\label{rand}

     We will assume in Sections \ref{rand} - \ref{part2} that the group $F=F(x,y)$ is $2$-generated,
     $H$ is the normal closure of
      $x$ in $F$ and $N$ is the normal closure of $x$ in $H$.
      Note that a word $v=v(x,y)$ represents an element of $H$ if and only if $\sigma_y(v)=0$,
      where $\sigma_y(v)$ is the sum of the exponents at $y$ in the word $v$.

      We want to estimate the
     growth of the set $T$ introduced in Section \ref{cos} and make use of Lemma \ref{lone} (d). Fortunately,  we are able to compare $T$ with other sets whose
     growth is more adaptable to probabilistic methods.

     Let $T^+$ (respectively, $T^-$) consist of all words $w$ of $T$ starting with
     the letter $y$ (resp., with $y^{-1}$). Note that $T=T^+\sqcup T^-\sqcup\{1\}$ by Condition \ref{Cond2}, and there is an involution $T^+\leftrightarrow T^-$ given by the  rule $w\mapsto \overline{w}$, where the word $\overline{w}$ results from $w$ after the
     replacements of all the occurrences of $y$ by $y^{-1}$, and vice versa.

\begin{figure}[h!]
\begin{center}
%TeXCAD (http://texcad.sf.net/) Picture. File: [sub1.pic]. Options on following lines.
%\grade{\on}
%\emlines{\off}
%\epic{\off}
%\beziermacro{\on}
%\reduce{\on}
%\snapping{\off}
%\pvinsert{% Your \input, \def, etc. here}
%\quality{8.000}
%\graddiff{0.005}
%\snapasp{1}
%\zoom{4.0000}
\unitlength 1mm % = 2.845pt
\linethickness{0.4pt}
\ifx\plotpoint\undefined\newsavebox{\plotpoint}\fi % GNUPLOT compatibility
\begin{picture}(114,60)(0,70)
%\dashline{1}(14.75,104)(114,103.75)
\put(14.68,103.93){\line(1,0){.9925}}
\put(16.665,103.925){\line(1,0){.9925}}
\put(18.65,103.92){\line(1,0){.9925}}
\put(20.635,103.915){\line(1,0){.9925}}
\put(22.62,103.91){\line(1,0){.9925}}
\put(24.605,103.905){\line(1,0){.9925}}
\put(26.59,103.9){\line(1,0){.9925}}
\put(28.575,103.895){\line(1,0){.9925}}
\put(30.56,103.89){\line(1,0){.9925}}
\put(32.545,103.885){\line(1,0){.9925}}
\put(34.53,103.88){\line(1,0){.9925}}
\put(36.515,103.875){\line(1,0){.9925}}
\put(38.5,103.87){\line(1,0){.9925}}
\put(40.485,103.865){\line(1,0){.9925}}
\put(42.47,103.86){\line(1,0){.9925}}
\put(44.455,103.855){\line(1,0){.9925}}
\put(46.44,103.85){\line(1,0){.9925}}
\put(48.425,103.845){\line(1,0){.9925}}
\put(50.41,103.84){\line(1,0){.9925}}
\put(52.395,103.835){\line(1,0){.9925}}
\put(54.38,103.83){\line(1,0){.9925}}
\put(56.365,103.825){\line(1,0){.9925}}
\put(58.35,103.82){\line(1,0){.9925}}
\put(60.335,103.815){\line(1,0){.9925}}
\put(62.32,103.81){\line(1,0){.9925}}
\put(64.305,103.805){\line(1,0){.9925}}
\put(66.29,103.8){\line(1,0){.9925}}
\put(68.275,103.795){\line(1,0){.9925}}
\put(70.26,103.79){\line(1,0){.9925}}
\put(72.245,103.785){\line(1,0){.9925}}
\put(74.23,103.78){\line(1,0){.9925}}
\put(76.215,103.775){\line(1,0){.9925}}
\put(78.2,103.77){\line(1,0){.9925}}
\put(80.185,103.765){\line(1,0){.9925}}
\put(82.17,103.76){\line(1,0){.9925}}
\put(84.155,103.755){\line(1,0){.9925}}
\put(86.14,103.75){\line(1,0){.9925}}
\put(88.125,103.745){\line(1,0){.9925}}
\put(90.11,103.74){\line(1,0){.9925}}
\put(92.095,103.735){\line(1,0){.9925}}
\put(94.08,103.73){\line(1,0){.9925}}
\put(96.065,103.725){\line(1,0){.9925}}
\put(98.05,103.72){\line(1,0){.9925}}
\put(100.035,103.715){\line(1,0){.9925}}
\put(102.02,103.71){\line(1,0){.9925}}
\put(104.005,103.705){\line(1,0){.9925}}
\put(105.99,103.7){\line(1,0){.9925}}
\put(107.975,103.695){\line(1,0){.9925}}
\put(109.96,103.69){\line(1,0){.9925}}
\put(111.945,103.685){\line(1,0){.9925}}
%\end
\put(20,104.25){\line(0,1){11.5}}
\put(20,115.75){\line(1,0){8}}
\put(28,115.75){\line(0,1){4.25}}
\put(28,120){\line(-1,0){3.75}}
\put(24.25,120){\line(0,-1){8}}
\put(24.25,112.25){\line(1,0){12}}
\put(36,112.25){\line(0,-1){4}}
\put(36,108.25){\line(-1,0){4}}
\put(32,108.25){\line(0,1){7.75}}
\put(32.25,115.5){\line(1,0){12.25}}
\put(44.25,115.5){\line(0,1){4.25}}
\put(44.25,119.75){\line(-1,0){4}}
\put(40.25,119.75){\line(0,-1){7.75}}
\put(40.25,112.5){\line(0,1){1.5}}
\put(40.25,112){\line(1,0){8.25}}
\put(48.25,112.25){\line(0,-1){21.5}}
\put(48.25,91.5){\line(1,0){8.5}}
\put(56.25,95.5){\line(0,1){0}}
\put(56.25,95.5){\line(0,1){0}}
\put(56.75,95.5){\line(0,-1){3.75}}
\put(56.75,95.5){\line(-1,0){4}}
\put(52.75,95.5){\line(0,-1){8}}
\put(52.75,87.75){\line(1,0){12}}
\put(65,87.75){\line(0,-1){3.75}}
\put(65,84){\line(-1,0){4}}
\put(61,84){\line(0,1){14.5}}
\put(60.75,98.5){\line(1,0){.5}}
\put(61.25,98.5){\line(1,0){3.75}}
\put(65,98.75){\line(0,-1){4}}
\put(60.75,95){\line(1,0){.25}}
\put(60.75,94.75){\line(1,0){8.5}}
\put(69,94.5){\line(0,-1){3.5}}
\put(69,91){\line(1,0){3.75}}
\put(72.75,91){\line(0,1){24.5}}
\put(69,114.75){\line(1,0){3.75}}
\put(68.75,114.75){\line(0,-1){3.5}}
\put(68.75,111.25){\line(1,0){12.5}}
\put(81,114.75){\line(0,-1){3.25}}
\put(81,115){\line(1,0){4}}
\put(84.75,114.75){\line(0,-1){3.5}}
\put(84.75,111.25){\line(1,0){11.5}}
\put(96,111.25){\line(0,-1){15.75}}
\put(91.75,95.25){\line(1,0){4.25}}
\put(92.25,95.5){\line(0,-1){4.5}}
\put(88.5,91.25){\line(1,0){3.5}}
\put(89,107.25){\line(0,-1){16}}
\put(85.25,107.5){\line(1,0){3.75}}
\put(85,107.5){\line(0,-1){8.25}}
\put(84.75,99.25){\line(-1,0){3.5}}
\put(81.5,99.25){\line(0,-1){4.25}}
\put(81.5,94.25){\line(0,-1){1.25}}
\put(20,103.75){\circle{1.581}}
\put(48.25,103.5){\circle{1.5}}
\put(72.75,103.75){\circle{1.581}}
\put(85,103.75){\circle{1.803}}
\put(96,104){\circle{1.5}}
\put(89,103.75){\circle{1.581}}
\put(33.5,121.75){$w_1$}
\put(55,84.5){$w_2$}
\put(79,118){$w_3$}
\put(93.75,88.75){$w_4$}
\put(15.5,109.5){$y^{k_1}$}
\put(51.25,109.5){$y^{- \ell_1}$}
\put(42.75,96.75){$y^{- k_2}$}
\put(75.25,97){$y^{\ell_2}$}
\put(75.25,107.75){$y^{k_3}$}
\put(99,109.25){$y^{- \ell_3}$}
\put(100.5,99){$y^{- k_4}$}
\put(81.25,92){\line(-1,0){.25}}
\put(81.5,92.25){\line(0,-1){2}}
\put(19.5,107.5){\line(1,0){1.25}}
%\emline(19.75,111.5)(20.75,111.75)
\multiput(19.75,111.5)(.125,.03125){8}{\line(1,0){.125}}
%\end
\put(47.75,107.5){\line(1,0){1}}
\put(47.75,100.25){\line(1,0){1}}
\put(47.75,95.5){\line(1,0){1.25}}
\put(72.25,107){\line(1,0){1}}
\put(72.25,100.25){\line(1,0){.25}}
\put(72.25,99.75){\line(1,0){1.5}}
%\emline(72.25,95)(72.5,94.75)
\multiput(72.25,95)(.03125,-.03125){8}{\line(0,-1){.03125}}
%\end
\put(72.25,94.75){\line(1,0){1.25}}
\put(95.5,107.5){\line(1,0){1.25}}
\put(95.25,99.75){\line(1,0){1.5}}
\put(46.25,37){\line(-1,0){.5}}
\put(88.75,99.25){\line(1,0){.25}}
\put(88.25,99.25){\line(1,0){1.25}}
\put(88.5,95.5){\line(1,0){1.25}}
\end{picture}
\end{center}
\caption{The word $w$ in Equation (\ref{fact})}\label{fi}
\end{figure}
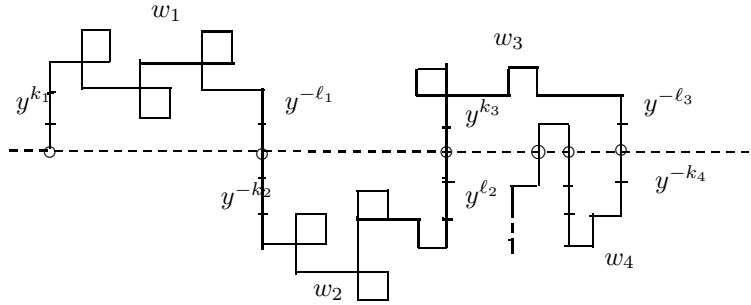

     Suppose $w\in T^+$ and $v$ is the first non-empty prefix of $w$ such that $\sigma_y(v)=0$ (if any exists).
     Then it must end with $y^{-1}$, and by Condition \ref{Cond1}, the next letter of $w$ after $v$ (if it exists)
     must be equal to $y^{-1}$  since $\sigma_y(v)=0$ and so $v\in H$.  If $w\equiv vv'$,
     then $v'$ starts with $y^{-1}$, and $v'$ has no prefixes $v''x^{\pm 1}$ with $\sigma_y(v'')=0$
     because otherwise $vv''$ would be a prefix of $w$ with
     $$\sigma_y(vv'')=\sigma_y(v)+\sigma_y(v'')=0+0=0$$
     contrary Condition \ref{Cond1}. Hence $v'\in T^-$, and by induction we obtain the irreducible factorization
     \begin{equation}\label{fact}
     w\equiv (y^{k_1}w_1y^{-\ell_1})(y^{-k_2}w_2y^{\ell_2})(y^{k_3}w_3y^{-\ell_3})\dots (y^{\pm k_s}w_sy^{\mp \ell_s})(y^{\mp k_{s+1}}w_{s+1}),
     \end{equation}
     where $s\ge 0$, all the exponents $k_i,\ell_i$ are strictly positive, all the words $w_1,\dots,w_{s}$ start and  end with the letters $x^{\pm 1}$, and $w_{s+1}$, if nonempty, starts with $x^{\pm 1}$. We also have $\sigma_y(W_i)=0$ for  every word $W_i$ enclosed in the $i$-th pair of the parentheses of (\ref{fact}) if $i\le s$,
     and $\sigma_y(v)>0$ ($\sigma_y(v)<0$) if $v$ is a non-empty proper prefix of $W_i$, $1\le i\le s+1$,
      and $i$ is odd (respectively, even). See fig. \ref{fi}, where a word $w\in T^+$ is pictured
     as a path of length $|w|$ in the lattice $ {\mathbb Z}^2$: its edges going to the right, to the left,
     up, and down correspond to the letters $x$, $x^{-1}$, $y$, and $y^{-1}$ of the word $w$.

      There are $4\cdot 3^{n-1}$ reduced words of length $n\ge 1$ in the alphabet ${\cal A}=\{x^{\pm 1},y^{\pm 1}\}$.
     Therefore we get any particular reduced word $a_1\dots a_n$ of length $n$ with probability $(4\cdot 3^{n-1})^{-1}$
     if we choose any letter $a_1\in\cal A$ with probability $1/4$ as the first letter, and for $i\ge 2$, the letter
     $a_i$ is chosen from ${\cal A}\backslash\{a_{i-1}^{-1}\}$ with equal probabilities $1/3$. This probabilistic model
     can be regarded as the following 2-dimensional random walk on the integer lattice ${\mathbb Z}^2$ with the standard basis $\{e_1, e_2\}$.

     The correlated random variables ${\bf X}_1, {\bf X}_2,\dots$ take values in $\{\pm e_1, \pm e_2\}$ and their distribution
     is given by:
     $${\bf P (X}_1=e_1)={\bf P (X}_1=-e_1)={\bf P (X}_1=e_2)={\bf P (X}_1=-e_2)=1/4,$$
     and for $i\ge 2$ and any $e\in \{\pm e_1, \pm e_2\}$, by conditional probabilities
     $${\bf P(X}_i=e\mid {\bf X}_{i-1}=-e)=0, \;\;{\bf P(X}_i=e\mid {\bf X}_{i-1}\ne -e)=1/3$$

     Let $\mathbf{S}_0=0$ and for $n\ge 1$, $\mathbf{S}_n=\mathbf{S}_{n-1}+{\bf X}_n$. Then $\{\mathbf{S}_n\}_{n\ge 0}$ is a
     correlated random walk on ${\mathbb Z}^2$. First we want to associate with this random walk,
     an uncorrelated (that is Markov) 1-dimensional random walk $\{{\bf s}_m\}_{m\ge 1}$.
     (A replacement of $\{{\bf S}_n\}_{n\ge 0}$ by $\{{\bf s}_m\}_{m\ge 1},$ but without lowering
     of the dimension, was used by Gillian Iossif in \cite{I}.)
     The process $\{{\bf s}_m\}_{m\ge 1}$ is obtained by observing the random walk $\{\mathbf{S}_n\}_{n\ge 0}$ only at
     the times when ${\bf X}_n=e_2$ and taking the projection $(\mathbf{S}_n)_y$ of $\mathbf{S}_n$ on the $y$-axis.

     \begin{lemma} \label{21} The process $\{{\bf s}_m\}_{m\ge 2}$ is a random walk
     on $\mathbb Z$ with the following equal distributions of independent integer-valued variables ${\bf Y}_m = {\bf s}_m-{\bf s}_{m-1}$:
 \begin{eqnarray*}
     &&{\bf P (Y}_m = k)=0,\mbox{ \emph{for} }k\ge 2,\\
     &&{\bf P (Y}_m =1)=2/3,\\
     &&{\bf P (Y}_m =k)=2^{-k}\cdot 3^{k-2},\mbox{ \emph{for} }k\le 0.
     \end{eqnarray*}
     The mean value ${\bf E}({\bf Y}_m)$ is $0$ and the variance $\sigma^2$ of ${\bf Y}_m$ does exist.
     \end{lemma}

     \proof  Let us consider an auxiliary discrete Markov chain $\{{\bf\xi}_i\}_{i\ge 1}$ with two states $e_2$ and $-e_2$. It is obtained
     by observing the process ${\bf X}_n$ when ${\bf X}_n=\pm e_2.$  Denote by $(p_{\alpha,\beta})$ the $2\times 2$-matrix of transition probabilities, where $p_{11}$ is the probability of transition from $e_2$ to $e_2$, and so on.
If $X_n=e_2$ is observed at some moment $n$, then
      \begin{eqnarray*}
     &&p_{11}=\sum_{j=1}^{\infty} {\bf P(X}_{n+1}=\dots=  {\bf X}_{n+j}=e_1,
      {\bf X}_{n+j+1}=e_2\mid{\bf X}_n=e_2)\\
    && +\sum_{j=1}^{\infty} {\bf P(X}_{n+1}=\dots=  {\bf X}_{n+j}=-e_1,
      {\bf X}_{n+j+1}=e_2\mid {\bf X}_n=e_2)\\
     &&+ {\bf P(X}_{n+1}=e_2 \mid {\bf X}_n=e_2)
      = 2\sum_{j=1}^{\infty}\left(\frac{1}{3}\right)^j\frac{1}{3} +\frac{1}{3} =2 \cdot\frac{1}{6}+\frac{1}{3}=\frac{2}{3}
     \end{eqnarray*}
Similarly we have $p_{22}=2/3$ and $p_{12}=p_{21}=1/3$.

Returning to the statement of the lemma we assume that ${\bf Y}_{m-1}$ is  observed at some moment $r$, i.e., $ {\bf s}_{m-1}=\mathbf({S}_r)_y$. Then we obviously have ${\bf P (Y}_m = k)=0$ if $k\ge 2$ since the walker cannot miss the values $\mathbf{S}_t$ with ${\bf X}_t=e_2$.
%  Next, ${\bf P (Y}_m =1)=p_{11}$

If ${\bf Y}_m=1$, then  in the corresponding Markov chain $\{{\bf \xi}_i\}$, we see exactly one transition, namely from $e_2$ to $e_2$,
and so ${\bf P( Y}_m=1)=p_{11}=2/3.$

 If $k\le 0$, the number of times such that the variable  ${\bf X}_t$, with $t>n$, has to take the value $-e_2$
      before the next value $e_2$ is taken and the jump ${\bf Y}_m$ happens, equals  $-k+1$.
In other words, after the value ${\bf \xi}_i=e_2$ is taken, we observe $-k+1$ values $-e_2$ and then again
$e_2$ in the Markov chain $\{{\bf \xi}_i\}$. Therefore
$${\bf P( Y}_m=k)=p_{12}p_{22}^{-k}p_{21}=\frac{1}{3}\left(\frac{2}{3}\right)^{-k}\frac{1}{3}=2^{-k}\cdot 3^{k-2}$$
if $k\le 0$. The first claim of the lemma is proved since the obtained distribution of ${\bf Y_m}$ does not depend on
     the value taken by ${\bf Y}_{m-1}$.

Finally,
\begin{equation}\label{mean}
    {\bf E(Y}_m)=\sum_{k=-\infty}^{\infty}k{\bf P(Y}_m=k)=2/3+\sum_{k=1}^{\infty}{-k}\cdot 2^k\cdot 3^{-k-2}=0
    \end{equation}
    and there exists the variance
    \begin{equation}\label{var}
    \sigma^2=Var({\bf Y}_m)= \sum_{k=-\infty}^{\infty}k^2{\bf P(Y}_m=k)<\infty
    \end{equation}
\endproof

     The next lemma gives the estimate of the probability $$p_m={\bf P(s}_1>0,{\bf s}_2> 0,\dots, {\bf s}_m> 0)$$
      that the entire walk of length $m$ belongs to the positive semiaxis.

    \begin{lemma}\label{pm} We have $p_m=\Theta(\frac{1}{\sqrt{m}})$.
    \end{lemma}

\proof
    It follows from conditions (\ref{mean}) and (\ref{var}) that the series
    \begin{equation}\label{Fe}
    \sum_{j=1}^{\infty}\frac1j (({\bf P(s}_j>0)-\frac12 )
    \end{equation}
    converges at least conditionally (see \cite{F}, Theorem 1 in
    Subsection XVIII.5). In turn, the convergence of (\ref{Fe}) implies that $p_m=\Theta(\frac{1}{\sqrt{m}})$
    (\cite{F}, Theorem 1a in Subsection XII.7).
    \endproof

    Now we will apply Lemma \ref{21} and obtain an asymptotic property of the original correlated walk $\{{\bf S}\}_{n\ge 0}$ or, equivalently, of a random reduced word $w$ of length $n$ in the alphabet $\{x^{\pm 1}, y^{\pm 1}\}$.
    Below we denote by $W_n$ the set of all reduced words of length $n$ in this alphabet, and the subset $U_n$ consists of all the words $u$ such that $\sigma_y(v)\ge 0$ for every  prefix $v$ of $u$.

    \begin{lemma}\label{UU} Let $P_n={\bf P}(u\in U_n\mid u\in W_n),$
i.e., \\$P_n={\bf P((S}_1)_y \ge 0, ({\bf S}_2)_y \ge 0,\dots, ({\bf S}_n)_y \ge 0)$. Then $P_n =\Theta(\frac{1}{\sqrt{n}})$.
    \end{lemma}
    \proof The subset $U_n$ is partitioned in two subsets, namely,
    $U'_n$ contains all the words $u$ such that the number $m$ of the occurrences of $y$ in $u$ belongs to the
    segment $I_n=[\frac{n}4- \frac{n}8, \frac{n}4+ \frac{n}8]$, and $U''_n$ is the compliment of $U'_n$ in $U_n$.

    Recall that the number of all reduced words in the alphabet $\{x^{\pm 1}, y^{\pm 1}\}$ is $4\cdot 3^{n-1}$.
    The mean number of the occurrences of $y$ in arbitrary reduced word $w$ of length $n$ is $n/4$. It follows from the Large Deviation Theory that
    $$\lim_{n\to\infty}\frac{\# U''_n}{4\cdot 3^{n-1}}=0,$$ and the rate of convergence is exponential. In fact this is true if one replace $1/8$ by any $\varepsilon>0$ in the definition of $I_n$ (for instance, see \cite{KSS}). Therefore
    \begin{equation}\label{exp}
   {\bf P}(u\in U''_n\mid u\in W_n)=O(\lambda^n)
    \end{equation}
    for some positive $\lambda<1$.

    Now for a random walk ${\bf S}_n={\bf X}_1+\dots+{\bf X}_n$, we denote by $m({\bf S}_n)$ the number of steps $e_2$.
    Then
    $${\bf P}(u\in U'_n\mid u\in W_n)
     = {\bf P((S}_1)_y\ge 0,\dots, {\bf (S}_n))_y\ge 0, m({\bf S}_n)\in I_n)$$
    $$={\bf P(s}_1>0,\dots, {\bf s}_{m({\bf S}_n)}>0, m({\bf S}_n)\in I_n)$$
     $$=\sum_{j\in I_n} {\bf P(s}_1>0,\dots, {\bf s}_{m({\bf S}_n)}>0, m({\bf S}_n)=j)$$
    $$=\sum_{j\in I_n} {\bf P(s}_1>0,\dots, {\bf s}_{m({\bf S}_n)}>0\mid m({\bf S}_n)=j)\cdot{\bf P}(m({\bf S}_n)=j)$$
    \begin{equation}\label{last}
    =\sum_{j\in I_n} {\bf P(s}_1>0,\dots, {\bf s}_j >0)\cdot{\bf P}(m({\bf S}_n)=j)
    \end{equation}

    By Lemma \ref{pm}, $${\bf P(s}_1>0,\dots {\bf s}_j >0)=\Theta\left(\frac{1}{\sqrt{j}}\right)=\Theta\left(\frac{1}{\sqrt{n}}\right)$$
    because $j\in I_n=[\frac{n}4- \frac{n}8, \frac{n}4+ \frac{n}8]$. Therefore by (\ref{exp}) and by the definition of $U''_n$, the sum (\ref{last}) can be rewritten as
    $$\Theta\left(\frac{1}{\sqrt{n}}\right)\sum_{j\in I_n}{\bf P}(m({\bf S}_n)=j)=
    \Theta\left(\frac{1}{\sqrt{n}}\right){\bf P}(m({\bf S}_n)\in I_n)$$ $$=\Theta\left(\frac{1}{\sqrt{n}}\right)(1-O(\lambda^n))=\Theta\left(\frac{1}{\sqrt{n}}\right)$$
    Thus
    \begin{equation}\label{U'}
    {\bf P}(u\in U'_n\mid u\in W_n)=\Theta\left(\frac{1}{\sqrt{n}}\right)
    \end{equation}

    The comparison of (\ref{exp}) and (\ref{U'}) shows that
    ${\bf P}( u\in U''_n\mid u\in W_n)=o({\bf P} (u\in U'_n\mid u\in W_n)$ and so $P_n=\Theta(\frac{1}{\sqrt{n}})$, as required.
    \endproof

\begin{lemma}\label{ptwo} The cogrowth function $f_N(n)$ with respect to the free generators $x$ and $y$ is at least  $\Theta(\frac{3^n}{\sqrt{n}})$.
\end{lemma}
\proof The function $f_N(n)$ is equal to the growth function $g_T(n)$ of the shortest Schreier transversal
of $N$ in $F$ (see Lemma \ref{lone} (d)). Note that $yU\subset T^+$ by the definitions of $H$, $T^+$, and $U$. So we have $g^s_{T^+}(n)\ge g^s_U(n-1)$ ($n\ge 2$) for the spherical growth functions of the sets $T^+$ and $U$,
which count only the words of length $n$ in these subsets.
By the definitions of the (correlated) random walk $\{{\bf S}_n\}$ and the probabilities $P_n,$ we
have $g^s_U(n-1)=P_{n-1}\cdot4\cdot 3^{n-2}$ since $4\cdot3^{n-2}$ is the number of all reduced words of length $n-1$,
 whence by Lemma \ref{UU}, $$g_T(n)\ge g_{T^+}(n)=2+\sum_{i=2}^{n} g^s_{T^+}(i)\ge 2+\Theta\left(\sum_{i=1}^{n-1} g^s_U(i)\right)=$$ $$2+\Theta\left(\sum_{i=1}^{n-1} P_i\cdot 4\cdot 3^{i-1}\right)=\Theta\left(\sum_{i=1}^{n-1} \frac{3^i}{\sqrt{i}}\right)=\Theta\left(\frac{3^{n-1}}{\sqrt{n-1}}\right)=\Theta\left(\frac{3^{n}}{\sqrt{n}}\right),$$ %which completes the proof since $\#(T^+\cap W_n)=\#(T^-\cap W_n).$
 \endproof

\section{Upper bound in Theorem \ref{cogr} (3).}\label{upp}

To replace ``at least'' by ``at most'' in the formulation of Lemma \ref{ptwo}, we need two different
combinatorial estimates for the number of words $w\in T_n=T\cap W_n$. We denote by $cr(w)$ the number of times the
corresponding to $w$ walk $p=p(w)$  crosses the $x$-axis. ($cr(w)=5$ for the word $w$ pictured at Fig.\ref{fi}.)
We will not count the vertex $p_-$ to $cr(w)$ but will count $p_+$ if it belongs to the $x$-axis and $|p|>0$. Let us start with the easier estimate.

\begin{lemma} \label{over} Let $\overline T_n$ be the subset of $T_n=W_n\cap T$ consisting of the words $w$ with $cr(w)\ge \frac12\log_3 n.$ Then $\#\{\overline T_n\}=O\left(\frac{3^n}{\sqrt n}\right)$.
\end{lemma}

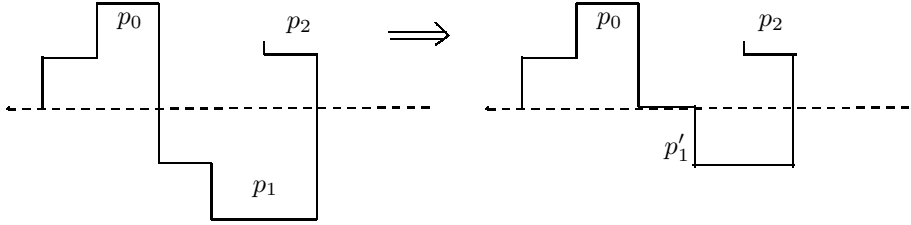
\begin{figure}[h!]
\begin{center}
%TeXCAD (http://texcad.sf.net/) Picture. File: [subsub5.pic]. Options on following lines.
%\grade{\on}
%\emlines{\off}
%\epic{\off}
%\beziermacro{\on}
%\reduce{\on}
%\snapping{\off}
%\pvinsert{% Your \input, \def, etc. here}
%\quality{8.000}
%\graddiff{0.005}
%\snapasp{1}
%\zoom{4.0000}
\unitlength 1mm % = 2.845pt
\linethickness{0.4pt}
\ifx\plotpoint\undefined\newsavebox{\plotpoint}\fi % GNUPLOT compatibility
\begin{picture}(132.75,52.5)(12,110.00)
%\dottedline(12.75,138.5)(13,138.75)
\multiput(12.68,138.43)(.125,.125){3}{{\rule{.4pt}{.4pt}}}
%\end
%\dottedline(76.5,138.5)(76.75,138.75)
\multiput(76.43,138.43)(.125,.125){3}{{\rule{.4pt}{.4pt}}}
%\end
%\dashline{1}(13,138.5)(69,138.75)
\put(12.93,138.43){\line(1,0){.9825}}
\put(14.895,138.438){\line(1,0){.9825}}
\put(16.86,138.447){\line(1,0){.9825}}
\put(18.824,138.456){\line(1,0){.9825}}
\put(20.789,138.465){\line(1,0){.9825}}
\put(22.754,138.474){\line(1,0){.9825}}
\put(24.719,138.482){\line(1,0){.9825}}
\put(26.684,138.491){\line(1,0){.9825}}
\put(28.649,138.5){\line(1,0){.9825}}
\put(30.614,138.509){\line(1,0){.9825}}
\put(32.579,138.517){\line(1,0){.9825}}
\put(34.544,138.526){\line(1,0){.9825}}
\put(36.509,138.535){\line(1,0){.9825}}
\put(38.474,138.544){\line(1,0){.9825}}
\put(40.438,138.553){\line(1,0){.9825}}
\put(42.403,138.561){\line(1,0){.9825}}
\put(44.368,138.57){\line(1,0){.9825}}
\put(46.333,138.579){\line(1,0){.9825}}
\put(48.298,138.588){\line(1,0){.9825}}
\put(50.263,138.596){\line(1,0){.9825}}
\put(52.228,138.605){\line(1,0){.9825}}
\put(54.193,138.614){\line(1,0){.9825}}
\put(56.158,138.623){\line(1,0){.9825}}
\put(58.123,138.631){\line(1,0){.9825}}
\put(60.088,138.64){\line(1,0){.9825}}
\put(62.053,138.649){\line(1,0){.9825}}
\put(64.017,138.658){\line(1,0){.9825}}
\put(65.982,138.667){\line(1,0){.9825}}
\put(67.947,138.675){\line(1,0){.9825}}
%\end
%\dashline{1}(76.75,138.5)(132.75,138.75)
\put(76.68,138.43){\line(1,0){.9825}}
\put(78.645,138.438){\line(1,0){.9825}}
\put(80.61,138.447){\line(1,0){.9825}}
\put(82.574,138.456){\line(1,0){.9825}}
\put(84.539,138.465){\line(1,0){.9825}}
\put(86.504,138.474){\line(1,0){.9825}}
\put(88.469,138.482){\line(1,0){.9825}}
\put(90.434,138.491){\line(1,0){.9825}}
\put(92.399,138.5){\line(1,0){.9825}}
\put(94.364,138.509){\line(1,0){.9825}}
\put(96.329,138.517){\line(1,0){.9825}}
\put(98.294,138.526){\line(1,0){.9825}}
\put(100.259,138.535){\line(1,0){.9825}}
\put(102.224,138.544){\line(1,0){.9825}}
\put(104.188,138.553){\line(1,0){.9825}}
\put(106.153,138.561){\line(1,0){.9825}}
\put(108.118,138.57){\line(1,0){.9825}}
\put(110.083,138.579){\line(1,0){.9825}}
\put(112.048,138.588){\line(1,0){.9825}}
\put(114.013,138.596){\line(1,0){.9825}}
\put(115.978,138.605){\line(1,0){.9825}}
\put(117.943,138.614){\line(1,0){.9825}}
\put(119.908,138.623){\line(1,0){.9825}}
\put(121.873,138.631){\line(1,0){.9825}}
\put(123.838,138.64){\line(1,0){.9825}}
\put(125.803,138.649){\line(1,0){.9825}}
\put(127.767,138.658){\line(1,0){.9825}}
\put(129.732,138.667){\line(1,0){.9825}}
\put(131.697,138.675){\line(1,0){.9825}}
%\end
\put(17.5,138.5){\line(0,1){7}}
\put(81.25,138.5){\line(0,1){7}}
\put(17.75,145.25){\line(1,0){7}}
\put(81.5,145.25){\line(1,0){7}}
\put(24.75,145.25){\line(0,1){7.25}}
\put(88.5,145.25){\line(0,1){7.25}}
\put(24.75,152.5){\line(1,0){8.25}}
\put(88.5,152.5){\line(1,0){8.25}}
\put(33,152.5){\line(0,-1){21.25}}
\put(33,131.25){\line(1,0){7}}
\put(40,131.25){\line(0,-1){7.5}}
\put(40,123.75){\line(1,0){14}}
\put(54,123.75){\line(0,1){22}}
\put(54,145.75){\line(-1,0){7}}
\put(117.75,145.75){\line(-1,0){7}}
\put(47,145.75){\line(0,1){1.75}}
\put(110.75,145.75){\line(0,1){1.75}}
\put(96.75,152.5){\line(0,-1){13.75}}
\put(96.75,138.75){\line(1,0){7.25}}
\put(104.25,139){\line(0,-1){8.25}}
\put(104,131){\line(1,0){13.5}}
\put(117.25,145.75){\line(0,-1){15}}
\put(27.5,149.75){$p_0$}
\put(45.5,127.5){$p_1$}
\put(50,149.25){$p_2$}
\put(91.25,149.75){$p_0$}
\put(100.25,132.5){$p'_1$}
\put(112.75,149.5){$p_2$}
\put(63.5,148.75){\line(1,0){7.75}}
\put(63.75,147.75){\line(1,0){7}}
%\emline(69.75,150)(71.5,148.25)
\multiput(69.75,150)(.03365385,-.03365385){52}{\line(0,-1){.03365385}}
%\end
%\emline(69.75,147)(71.5,148.25)
\multiput(69.75,147)(.04605263,.03289474){38}{\line(1,0){.04605263}}
%\end
\end{picture}
\caption{Transformation $p\mapsto p'$ in Lemma \ref{over}}\label{trans}
\end{center}
\end{figure}

\proof It suffices to obtain the similar estimate for $\#\{\overline T^+_n\}=\#\{\overline T_n\cap T^+\}$.
The path $p=p(w)$ factorizes as
$p=p_0p_1\dots p_t,$ where all the edges of $p_0, p_2,\dots$ (of $p_1, p_3,\dots$) belong to
the upper (respectively, to the lower) half-plane, and $t\ge \frac12\log_3 n.$ Recall that
the paths $p_1, p_3,\dots$ start and terminate (except for $p_t$) with the edges parallel to $\pm e_2,$ and the ordinates of the vertices in all {\it other} their edges are at most $-1$. So if we remove the original and the terminal
edges from each of $p_1, p_3,\dots$ and lift the remaining subpaths of these paths by one unit,
we obtain a new path $p'=p_0p'_1p_2p'_3\dots$ of length  $n-t$, where $p'_1, p'_3,\dots$ are still in the lower
half-plane, and $p'=p(w')$ corresponds to some reduced word $w'$ as well
(may be $w'\notin \overline T_{n-t}$, see fig. \ref{trans}).

Note that given $p'$, one can restore the path $p$ since the only edges of $p'$ starting or ending at
the $x$-axis and belonging to the upper half-plane are the original and terminal edges of $p_0,p_2,\dots$.
Therefore the mapping $p\mapsto p'$ is injective, and so $\#\{\overline T^+_n\}$ cannot exceed the number of all reduced
words of length $\le n-\frac12\log_3 n,$ which is $\Theta(3^{n-\frac12\log_3 n})=\Theta\left(\frac{3^n}{\sqrt n}\right)$.
\endproof

\bigskip

Now we want to obtain an upper bound for the growth of the compliment $\overline{\overline {T_n}}=T_n\backslash \overline  T_n$. As a preliminary, we prove a statement on large deviations. We do this combinatorially since unlike the classical formulations, the walk ${\bf S}_n$ is correlated. We need a stronger estimate for the``tail'', which is aggravated  by the demand that the deviation interval should be $o(n)$. (One can say that the deviations are medium rather than large in the next lemma.)

For a random word $w\in W_n$, we denote by $n_x$ (by $n_{x^{-1}}, n_{y},$ and $ n_{y^{-1}}$ ) the number of letters
    $x$ (respectively, $x^{-1}$, $y$, and $y^{-1}$) in  $w$. Let $n_{ab}$  be the number of the subwords $ab$ in the word $w$, where $a,b\in\{x^{\pm 1},y^{\pm 1}\}$ and $a\ne b^{-1}$. We also denote by $s=s(w)$ the number
    of the $y$-syllabi in the factorization
    \begin{equation}\label{s}
    w\equiv (x^{k_0}) y^{\ell_1}\dots x^{k_{s-1}}y^{\ell_s}(x^{k_s})
    \end{equation}
    where only $k_0$ and $k_s$ can be zero. The number of positive (negative) exponents $k_i$ is denoted by $s_+$
    (respectively, by $s_-$).

\begin{lemma} \label{devi} For every sufficiently large $n$, for  each $a\in\{x^{\pm 1},y^{\pm 1}\}$ and each reduced $2$-letter word $ab$, we have

\begin{equation}\label{de1}
 {\bf P}\left(s\in \left(\frac{n}{3}- n^{3/4};\; \frac{n}{3}+n^{3/4}\right)\right)> 1-\exp({-n^{1/2}})
 \end{equation}
\begin{equation}\label{de4}
 {\bf P}\left(s_+\in \left(\frac{s}{2}- n^{3/4};\; \frac{s}{2}+n^{3/4}\right)\right)> 1-\exp({-n^{1/2}})
 \end{equation}
\begin{equation}\label{de2}
{\bf P}\left(n_a\in \left(\frac{n}{4}- 2n^{3/4};\; \frac{n}{4}+2n^{3/4}\right)\right)> 1-\exp({-n^{1/2}})
 \end{equation}
\begin{equation}\label{de3}
{\bf P}\left(n_{ab}\in \left(\frac{n}{12}-3n^{3/4};\; \frac{n}{12}+ 3n^{3/4}\right)\right)> 1-\exp({-n^{1/2}})
 \end{equation}
\end{lemma}

\proof (1) It suffices to proof the analogs of the formulae  (\ref{de1}) - (\ref{de3}) under each of the following
four conditions: (a) $k_0=0$, $k_s\ne 0$; (b) $k_0\ne 0$, $k_s= 0$; (c) $k_0=0$, $k_s=0$; (d) $k_0\ne 0$, $k_s\ne 0$.
Respectively, the set $W_n$ is the disjoint union of four subsets $W_n(a),\dots,W_n(d).$
We will assume that the condition (a) holds since the other cases are similar.

Let $n_1=n_x+n_{x^{-1}}$ and $n_2=n_y+n_{y^{-1}}$. Thus $n_1+n_2=n$. To obtain an arbitrary word
$w$ with fixed $n_1$, $n_2$ and $s$ we first assign the signs to each of the exponents $\ell_1,\dots, k_s$
(their modules are not yet assigned). There are $2 ^{2s}$ ways to do this. So we should put one $x^{\pm 1}$
(one $y^{\pm 1}$) to determine the signs for all even (all odd) syllabi of $w$ in the factorization (\ref{s}).
To continue the creation of the product (\ref{s}),  we should distribute the remaining  $n_1-s$ letters $x$ ($n_2-s$ letters $y$)
among the $s$  syllabi (with already determined exponents $\pm 1$ for every syllabi). We have ${n_1-1\choose s-1}$ (respectively, ${n_2-1\choose s-1}$) possible distributions.
Hence the number of the words in $W_n(a)$ with fixed $n_1, n_2$ and $s$ is equal
to
\begin{equation} \label{12s}
N(n_1,n_2,s)= 2^{2s}{n_1-1\choose s-1} {n_2-1\choose s-1}
\end{equation}

If $n_1>n_2+1$, then $$N(n_1,n_2,s) N(n_1-1,n_2+1,s)^{-1}= \frac{(n_1-1)(n_2+1-s)}{(n_1-s)n_2}\le 1$$
for every $s\ge 1$. Hence for a fixed $s,$ the value $N(n_1,n_2,s)$ is maximal when $|n_1-n_2|\le 1$.
Therefore $$\sum_{n_1+n_2=n}N(n_1,n_2,s)\le (n+1)N(m,n-m,s),$$ where $m=\lfloor n/2\rfloor.$ Thus, we obtain
the inequality (\ref{de1}) of the lemma if we prove for large enough $n$ the inequality
\begin{equation}\label{des}
{\bf P}\left(\left|\frac{n}{3} -s\right|\ge n^{3/4}\mid w\in W_n(a), n_1=\lfloor n/2 \rfloor\right)< \exp(-n^{1/2})/(n+1)
\end{equation}

Let $s_1=\lfloor n/3\rfloor$ and $s_2\ge s_1 +n^{3/4}$. Then for $m=\lfloor n/2 \rfloor$, we have
$$N(m,n-m,s_2) N(m,n-m,s_1)^{-1}$$ $$=2^{2(s_2-s_1)} \prod_{i=0}^{s_2-s_1-1} \frac{m-s_1-i}{s_2-1-i}
\prod_{i=0}^{s_2-s_1-1} \frac{n-m-s_1-i}{s_2-1-i}$$
\begin{equation}\label{s1s2}
= \prod_{i=0}^{s_2-s_1-1} \frac{2(m-s_1-i)}{s_2-1-i}
\prod_{i=0}^{s_2-s_1-1} \frac{2(n-m-s_1-i)}{s_2-1-i}
\end{equation}

Since $m=\lfloor n/2\rfloor$,  $s_1=\lfloor n/3\rfloor$, and $s_2-s_1\ge n^{3/4},$ we have $$\frac{2(m-s_1-i)}{s_2-1-i}\le \frac{\frac{n}{3}+1-2i}{\frac{n}{3}+n^{3/4}-i}<1-n^{-1/4}$$
for every factor of the first product in (\ref{s1s2}) and big enough $n$. Since $s_2-s_1\ge n^{3/4}$, this product does not exceed
$$(1-n^{-1/4})^{n^{3/4}}= ((1-n^{-1/4})^{n^{1/4}})^{n^{1/2}}\le \exp(-n^{1/2})$$
since $(1-x^{-1})^x <e^{-1}$ for $x>1.$
The same estimate works for the factors of the second product in (\ref{s1s2}). Hence
\begin{equation}\label{s21}
N(m,n-m,s_2) N(m,n-m,s_1)^{-1}\le \exp(-2n^{1/2})<\exp(-n^{1/2})/(n+1)
\end{equation}
for every sufficiently large $n.$

If $s_1=\lceil n/3\rceil$ and $s_2\le s_1-n^{3/4}$, then
$$N(m,n-m,s_2) N(m,n-m,s_1)^{-1}$$ $$= \prod_{i=0}^{s_1-s_2-1} \frac{s_1-1-i}{2(m-s_2+1-i)}
\prod_{i=0}^{s_1-s_2-1} \frac{s_1-1-i}{2(n-m-s_2-i)},$$
and again we have the same upper estimate $1-n^{-1/4}$ for the factors of both products,
which again produces the estimate (\ref{s21})
for large $n$. So the inequality (\ref{s21}) holds for any $s\notin \left(\frac{n}{3}- n^{3/4};\; \frac{n}{3}+n^{3/4}\right)$, which implies (\ref{des}), as desired.
Thus the inequality (\ref{de1}) of the lemma is proved for large enough $n$.

\medskip

(2) Note that one  may continue proving under the additional condition that $s\in (\frac{n}{3}-n^{3/4}, \frac{n}{3}+n^{3/4})$ (The probability
that $s$ does not belong to this interval is even $100$ times less then given by (\ref{de1}) since
$\exp(-2n^{1/2})<\exp(-n^{1/2})/(100n+100)$ if $n$ is large enough.) Under this condition, one  may switch the parameters $n_1$ and $s$ in the argument of part (1) and
obtain inequality
$${\bf P}\left(n_1\notin \left(\frac{n}{2}- n^{3/4};\; \frac{n}{2}+n^{3/4}\right)\right)< \exp(-n^{1/2})$$
(and even with additional factor $1/100$ in the right-hand side).

So it suffices to obtain the inequality  (\ref{de4})-(\ref{de3}) (or even $100$ times stronger) under the condition that
$n_1\in \left(\frac{n}{2}- n^{3/4};\; \frac{n}{2}+n^{3/4}\right)$.

The number of possible distributions of $s_+$ signs '$+$', denoted as  $M(s_+)=M(n,n_1,s,s_+)$, between the $x$-syllabi
is ${s \choose s_+}$, and if $s_+- s_0\ge n^{3/4}$ for $s_0=\lfloor s/2\rfloor,$ then

  \begin{equation}\label{splus}
  M(s_+)M(s_0)^{-1} =\prod_{i=1}^{s_+-s_0} \frac{s-s_0-i+1}{s_+-i+1}
  \end{equation}
  Under our assumption on $s$, each factor of the right-hand side of (\ref{splus}) is less than $1-5n^{-1/4}$ for large enough $n$. Since the number of factors is at least $n^{3/4}$, the argument of part (1) shows that  $M(s_+)M(s_0)^{-1}< e^{-5n^{1/2}}$. The case $\lceil s/2\rceil -s_+\ge n^{3/4}$ is similar, and so
 the inequality (\ref{de4}) is proved, and we may further assume that
$|s_+-s/2|\le n^{3/4}$.

(3) To prove (\ref{de2}), we distribute $n_x-s_+$ letters $x$ among $s_+$ syllabi (recall that $s_+$ letters $x$ are already
assigned to each of these $s_+$ syllabi), and $n_{x^{-1}}-s_-=n_1-n_x-s+s_+$ letters $x^{-1}$ to be distributed among $s_-=s-s_+$ syllabi. The number of such distributions is equal to
$L(n_x)={n_x-1 \choose s_+-1}{n_1-n_x-1 \choose s-s_+-1}$. Again we obtain the inequality of the form
 $L(n_x)L(\lfloor n_1/2\rfloor)^{-1}< \exp(-n^{1/2})/100$ (under the above assumption on $n_1, s$ and $s_+$) for large values of $n$ if $|n_x-\lfloor n_1\rfloor /2| \ge n^{3/4}$.
 (Now we leave the verification to the reader.) This proves the inequality (\ref{de2}) for $n_x$ and, by the symmetry,
 for $n_{x^{-1}}, n_y$, and $n_{y^{-1}}$.

 \medskip

 (4) Note that for a fixed $s$ and any distribution
 of the signs of  $k_i$-s and $\ell_i$-s in (\ref{s}), we have the same number of words
 with  prescribed $2s$ signs $(\pm,\pm,\dots,\pm)$ to the exponents. The number $r$ of the subwords
 $yx$ in $w$ depends only of the sign distribution. Therefore we now forget on the modules of the exponents
 and look at the number of distributions of signs $K(r)$ having exactly $r$ subwords $++$, where the
 second $+$ occupies an even position (i.e. corresponds to some $k_i$). Our goal is to uniformly
 estimate the numbers $K(r)K(\lfloor n/12\rfloor)^{-1}$ if $|r-(n/12)|\ge n^{3/4}.$

 In the above notation, we have $s_+$ pluses at the even positions, i.e., there are
 ${s\choose s_+}$ choices of signs at the even positions for fixed $s$ and $s_+$ . Then we should select $r$ of them (${s_+\choose r}$ possibilities)  and
 put $r$ pluses before each of them, put $s_+-r$ minuses before other $s_+-r$ positions and put
 any signs into the remaining $s_-=s-s_+$ odd positions ($2^{s_-}$ variants).  Hence
 $$K(r)=K(r,s,s_+)= \frac{2^{s_-} s!}{(s-s_+)!r!(s_+-r)!}$$  Repeating the approach we used
 a few times earlier, we get
 $$K(r)K(\lfloor s_+/2\rfloor )^{-1}=\frac{\lfloor s_+/2\rfloor !(s_+-\lfloor s_+/2\rfloor )!}{r!(s_+-r)!}< \exp(-n^{1/2})$$
 for big enough $n$
 if $|r-s_+/2|\ge n^{3/4},$ since again we may assume that $s$ and $s_+$ belong to the
 intervals $(\frac{n}{3}-n^{3/4}, \frac{n}{3}+n^{3/4})$ and $(\frac{s}{2}-n^{3/4}, \frac{s}{2}+n^{3/4})$,
 respectively. Thus the inequality (\ref{de3}) is proved for $n_{yx}$, and similarly, for all $n_{(x^{\pm 1}y^{\pm 1})^{\pm 1}}.$

 The inequality (\ref{de3}) for $n_{xx}$ (and, similarly, for $n_{x^{-1}x^{-1}},$ $n_{yy},$ and $n_{y^{-1}y^{-1}}$),
 follows now from the equality $n_{xx}=n_x-s$ and the inequalities for $n_x$ and $s$ obtained above. \endproof

 \vskip50ex

 Let us denote by $\hat W_n$ the set of all words $w$ from $W_n$ with all the parameters
 $s=s(w),$ $s_+=s_+(w)$, $n_a=n_a(w)$, and $n_{ab}=n_{ab}(w)$ belonging to the intervals introduced in the
 formulation of Lemma \ref{devi}. It follows that for large enough $n$, we have
 \begin{equation} \label{hat}
 {\bf P}(w\notin \hat W_n \mid w\in W_n) < 20 \exp(-n^{1/2})
\end{equation}

Let by definition, $\hat T_n=\overline{\overline T}_n\cap \hat W_n.$ We want to compare this set with the
set $S_n$ of all words $w\in W_n$ satisfying the following condition.
If $w\equiv v_1v_2$, where $\sigma_y(v_1)=0$ and the prefix $v_1$ ends with a letter $y^{\eta}$ ($\eta=\pm 1$),
then
\begin{itemize}
\item the word $v_2$, if non-empty, must start with $x^{\pm 1}$
\item the first occurrence of $y^{\pm 1}$ in $v_2$ (if any exists) must be $y^{\eta}$
\end{itemize}

\begin{lemma} \label{hatT} We have $\#\hat T_n=O(\# S_n)$.
\end{lemma}

\proof Let us say that two words $w$ and $w'$ from $\hat T_n$ have the same scheme $\Sigma$ if

(1) $s(w)=s(w')=s$ in their factorizations (\ref{s}),

 (2) these two words have the same $s$-tuples $\ell_1,\dots,\ell_s$,

 (3) every exponent $k_i$ is non-zero iff $k'_i$ is non-zero,

(4) we have $s_+(w)=s_+(w')=s_+$ for the numbers of positive $x$-syllabi in $w$ and in $w'$, and

(5) $n_x(w)=n_x(w')=n_x$  and $n_{x^{-1}}(w)=n_{x^{-1}}(w')=n_{x^{-1}}.$

Let us consider a scheme $\Sigma$ with  $k_0=0$ and $k_s\ne 0.$
To construct any word from $\Sigma,$ we start with a pattern word $y^{\ell_1}x^{\pm 1}\dots y^{\ell_s}x^{\pm 1}.$
At first we should choose  $s_+$ $x$-syllabi from $s$ ones and call them positive. We shall have one letter $x$
in each of the chosen syllabi and one $x^{-1}$ in each of the remaining $s_-=s-s_+$ $x$-syllabi. Then the remaining
$n_x-s_+$ letters $x$ and $n_{x^{-1}}-(s-s_+)$ letters $x^{-1}$ to be distributed among $s_+$  positive syllabi
and among $s-s_+$ negative ones, respectively. Therefore

\begin{equation} \label{Si} \#\Sigma= {s\choose s_+}{n_x-1 \choose s_+-1} {n_{x^{-1}}-1 \choose s-s_+-1}
\end{equation}

Note that all the words of $\Sigma$ have the same crossing number $cr(w)=t=t(\Sigma)$
since $cr(w)$ does not depend on the positions of letters $x^{\pm 1}$. Moreover,
the vertices, where the paths  $p(w)$ cross the $x$-axis, splits the corresponding syllabi $y^{\ell_i}$
identically for all $w\in \Sigma$, i.e. $y^{\ell_i}\equiv (y^{\ell'_i})(y^{\ell''_i})$ with nonzero
exponents ${\ell'_i}$ and ${\ell''_i}$ of the same sign by the definition of the set $T$ (see the factorization (\ref{fact}) and fig. \ref{fi}).
Therefore we can construct a new scheme $\Psi$
if we replace these syllabi $y^{l_i}$ by the products $y^{\ell'_i}x^{k'_i}y^{\ell''_i}$ with some nonzero $k'_i$.
Hence we have $\Psi\subset S_n$,
where, by definition of $\Psi$,  the parameters $n$, $n_x$, and $n_{x^{-1}}$ are preserved and the parameters $s$
and $s_+$ of $\Sigma$ are replaced by $s+t$ and $s_++t$, respectively. This substitution
turns the equality (\ref{Si}) in

\begin{equation} \label{Ps} \#\Psi= {s+t\choose s_++t}{n_x-1 \choose s_++t-1} {n_{x^{-1}}-1 \choose s-s_+-1}
\end{equation}

The restrictions imposed on $n_x, n_{x^{-1}}, s$ and $s_+$ imply that the binomial coefficients in (\ref{Ps}) are well-defined for large enough $n$ since $t$ is bounded by logarithm according to the definition of $\overline {\overline T}_n$.
Observe also that different schemes $\Sigma$ and $\Sigma'$ produce different $\Psi$ and $\Psi'$; indeed, given
$w\in \Psi$, one restores the exponents $\ell_1,\dots \ell_s$ corresponding to $\Sigma$ by erasing every occurrence
of $x^{\pm 1}$ in $w\equiv v_1x^{\pm 1}v_2$ if $\sigma_y(v_1)=0$. Therefore to prove the lemma, it suffices
to assume that $n$ is large enough and to obtain an inequality
\begin{equation}\label{C}
\#\Sigma\le C (\#\Psi)
\end{equation}
for every (non-empty) scheme $\Sigma$ and the corresponding $\Psi,$ where the positive constant $C$ must be independent of $n$ and $\Sigma$.
After cancelations in factorials, we obtain from (\ref{Si}) and (\ref{Ps}):
$$\frac{\#\Sigma}{\#\Psi}=\prod_{i=1}^{t}\frac{s_++i}{s+i}\prod_{i=1}^{t}\frac{s_+ +i-1}{n_x-s_+-i+1}$$
\begin{equation}\label{SP}
=\prod_{i=1}^{t}\frac{2(s_++i)}{s+i}\prod_{i=1}^{t}\frac{s_++i-1}{2(n_x-s_+-i+1)}
\end{equation}
Since $s_+<\frac{n}{6}+\frac32 n^{3/4}$, $s> \frac{n}{3}-n^{3/4}$, and $i\le t\le \frac12\log_3n< n^{3/4}$,
 we see that $\frac{2(s_++i)}{s+i}<1+20n^{-1/4}$
for every $i$ and large enough $n$. Similarly, every fraction $\frac{s_++i-1}{2(n_x-s_+-i+1)}$ is less than $1+40n^{-1/4}$ since $n_x>\frac{n}{4}-n^{3/4}$. Thus the inequality (\ref{SP}) gives us: $\frac{\#\Sigma}{\#\Psi}\le (1+40n^{-1/4})^{2t}$.
Since $(1+z)^{1/z}<e$ for positive $z$, it follows that $\frac{\#\Sigma}{\#\Psi}< \exp({40n^{-1/4}\log_3n})$. Here
the right-hand side is a bounded function of $n$, and so the desired estimate (\ref{C}) is obtained if
$k_0=0$ and $k_s\ne 0$. The other three cases are similar,
and the lemma is proved.
\endproof

\begin{lemma} \label{oover}The following is true: $\#\overline{\overline T}_n=O(3^n/\sqrt{n}).$
\end{lemma}

\proof By inequality (\ref{hat}), $\#(W_n\backslash \hat W_n)<20\exp(-n^{1/2})(4\times 3^{n-1})$,
and since $\exp(-n^{1/2})=o(n^{-1/2})$, to prove the lemma we may first replace $\overline{\overline T}_n$ by $\hat T_n$ in its formulation, and then Lemma \ref{hatT} allows us to replace $\hat T_n$ by $S_n$ in the formulation.
Moreover, we may replace $S_n$ by its "half" $S'_n$ consisting of the words, where the first $y$-syllabus
has positive exponent.

For the next replacement, we want to show that $\#S'_n\le \# U_n$, where the set $U_n$ is defined in Section \ref{rand}. To see this we exploit a 2-dimensional modification  of the known Reflection Principle.

Suppose $w\in S'_n$ and $p(w)=p_0q_1p_1q_2\dots q_tp_t$ is the factorization of the corresponding 2-dimensional walk, where the factors $p_0, p_2,\dots $, except for their initial and terminal vertices,  belong to the open upper half-plane, $p_1,p_3,\dots $ are, similarly,  in the lower half-plane, and $q_1,q_2,\dots$ lie
on the $x$-axis.

\begin{figure}[h!]
\begin{center}
%TeXCAD (http://texcad.sf.net/) Picture. File: [subsub6.pic]. Options on following lines.
%\grade{\on}
%\emlines{\off}
%\epic{\off}
%\beziermacro{\on}
%\reduce{\on}
%\snapping{\off}
%\pvinsert{% Your \input, \def, etc. here}
%\quality{8.000}
%\graddiff{0.005}
%\snapasp{1}
%\zoom{4.0000}
\unitlength 1mm % = 2.845pt
\linethickness{0.4pt}
\ifx\plotpoint\undefined\newsavebox{\plotpoint}\fi % GNUPLOT compatibility
\begin{picture}(125.5,42.5)(10,25)
\put(17.5,46.75){\line(0,1){12.75}}
\put(81.25,46.75){\line(0,1){12.75}}
\put(17.5,59.5){\line(1,0){13.5}}
\put(81.25,59.5){\line(1,0){13.5}}
\put(31,59.5){\line(0,-1){6}}
\put(94.75,59.5){\line(0,-1){6}}
\put(25,53.5){\line(1,0){6}}
\put(88.75,53.5){\line(1,0){6}}
\put(24.5,53.75){\line(0,-1){6.75}}
\put(88.25,53.75){\line(0,-1){6.75}}
\put(24.5,47){\line(1,0){12.5}}
\put(88.25,47){\line(1,0){12.5}}
\put(37,47){\line(0,-1){6.25}}
\put(37,40.75){\line(1,0){19}}
\put(56.25,46.75){\line(0,-1){6.25}}
\put(43.5,46.75){\line(1,0){12.5}}
\put(107.25,46.75){\line(1,0){12.5}}
\put(43.5,46.75){\line(0,1){12.75}}
\put(107.25,46.75){\line(0,1){12.75}}
\put(43.5,59.25){\line(1,0){6.25}}
\put(107.25,59.25){\line(1,0){6.25}}
\put(50.5,59.25){\line(1,0){1.5}}
\put(114.25,59.25){\line(1,0){1.5}}
\put(52.75,59.25){\line(1,0){1.5}}
\put(116.5,59.25){\line(1,0){1.5}}
\put(28.5,44.75){$q_1$}
\put(92.25,44.75){$q_1$}
\put(43.5,38){$p_1$}
\put(50.5,49.5){$q_2$}
\put(45.25,62){$p_2$}
\put(109,62){$p_2$}
\put(21.75,62.5){$p_0$}
\put(85.5,62.5){$p_0$}
%\dashline{1}(14.25,47)(24.25,47)
\put(14.18,46.93){\line(1,0){.9091}}
\put(15.998,46.93){\line(1,0){.9091}}
\put(17.816,46.93){\line(1,0){.9091}}
\put(19.634,46.93){\line(1,0){.9091}}
\put(21.452,46.93){\line(1,0){.9091}}
\put(23.271,46.93){\line(1,0){.9091}}
%\end
%\dashline{1}(78,47)(88,47)
\put(77.93,46.93){\line(1,0){.9091}}
\put(79.748,46.93){\line(1,0){.9091}}
\put(81.566,46.93){\line(1,0){.9091}}
\put(83.384,46.93){\line(1,0){.9091}}
\put(85.202,46.93){\line(1,0){.9091}}
\put(87.021,46.93){\line(1,0){.9091}}
%\end
%\dashline{1}(37,46.75)(43.5,47)
\put(36.93,46.68){\line(1,0){.9286}}
\put(38.787,46.751){\line(1,0){.9286}}
\put(40.644,46.823){\line(1,0){.9286}}
\put(42.501,46.894){\line(1,0){.9286}}
%\end
%\dashline{1}(100.75,46.75)(107.25,47)
\put(100.68,46.68){\line(1,0){.9286}}
\put(102.537,46.751){\line(1,0){.9286}}
\put(104.394,46.823){\line(1,0){.9286}}
\put(106.251,46.894){\line(1,0){.9286}}
%\end
%\dashline{1}(56,46.75)(61.75,46.75)
\put(55.93,46.68){\line(1,0){.9583}}
\put(57.846,46.68){\line(1,0){.9583}}
\put(59.763,46.68){\line(1,0){.9583}}
%\end
%\dashline{1}(119.75,46.75)(125.5,46.75)
\put(119.68,46.68){\line(1,0){.9583}}
\put(121.596,46.68){\line(1,0){.9583}}
\put(123.513,46.68){\line(1,0){.9583}}
%\end
\put(101.25,53.25){\line(0,-1){.25}}
\put(101,53.25){\line(0,-1){6.25}}
\put(101,53){\line(1,0){19}}
\put(120,53){\line(0,-1){6}}
\put(122.5,52.75){$\bar p_1$}
\put(113.25,43.75){$q_2$}
\put(64,51.5){\line(1,0){5.5}}
\put(64,50.5){\line(1,0){5}}
%\emline(68.5,52.75)(70,51.25)
\multiput(68.5,52.75)(.03333333,-.03333333){45}{\line(0,-1){.03333333}}
%\end
%\emline(70,51.25)(69.75,51.5)
\multiput(70,51.25)(-.03125,.03125){8}{\line(0,1){.03125}}
%\end
%\emline(70.25,50.75)(68.25,49.5)
\multiput(70.25,50.75)(-.05263158,-.03289474){38}{\line(-1,0){.05263158}}
%\end
\end{picture}
\caption{Transformation $p(w)\mapsto p(\bar w)$ in Lemma \ref{oover}}\label{refl}
\end{center}
\end{figure}
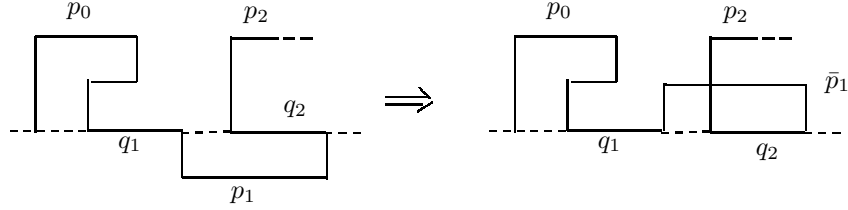

Denote by $p(\bar w)$ the path obtaining by mirror reflection of the
factors $p_1,p_3,\dots$ with respect to the $x$-axis (see fig. \ref{refl}).
The path $p(\bar w)$ is reduced because all the paths $q_1,\dots, q_t$ have positive lengths by the definition of the
set $S_n$. Hence  $\bar w\in U_n$. It remains to observe that the mapping $w\mapsto\bar w$ is injective because
the series of inverse reflections restoring $p(w)$ from $p(\bar w)$ is obvious: One should leave fixed the `even' pieces $p_0,p_2,...$ and reflect the remaining ones.  This proves the inequality $\#S'_n\le \# U_n.$

Thus, we may replace $S'_n$ by $U_n$, and the statement of the lemma follows now from Lemma \ref{UU}.
\endproof

{\bf Proof of Theorem \ref{cogr} (3).} By lemma \ref{lone} (d), the function $f_N(n)$ is the growth function
of the transversal $T$ of the subgroup $N$. Therefore by Lemmas \ref{over} and \ref{oover},
$f_N(n)= O(\sum_{i=0}^n \frac{3^i}{\sqrt{i}}) = O(\frac{3^n}{\sqrt{n}})$. This estimate together
with the opposite estimate from Lemma \ref{ptwo} prove the statement (3) of Theorem \ref{cogr}. $\Box$

\section{Part (2) of Theorem \ref{cogr}.}\label{part2}

Here we use the symbols  $F, H, N$, and $\{{\bf S}_n\}$ in the sense of Sections \ref{cos} and \ref{rand}.
    To apply Lemma \ref{maxgr} to the pair $(N, \langle x\rangle^N)$, we need two different (rough) estimates of the numbers of reduced words $w$ of length $n$ in $N,$ according to whether the number of factors $s$ in the factorizations (\ref{red}) of $w$ is small or not too small in comparison with $n$.

    We start with a simple $1$-dimensional random walk, where the variables ${\bf X}_i$ take two values $1$ and $-1$ with equal probabilities. Let ${\bf z}_0=0$ and ${\bf z}_i={\bf z}_{i-1}+{\bf X}_i$ for $i\ge 1$.

    \begin{lemma} \label{Q} Let $q_m$ denote the probability that the series ${\bf z}_1,\dots,{\bf z}_m$ has at least
    $m^{2/3}$ zeros. Then we have  $q_m =o \left(\exp\left(-\nu m^{1/6}\right)\right)$ for some positive $\nu$.
    \end{lemma}
    \proof
    Let $f_j(0\mid k)$ denote the probability that the random walker reaches $0$ for the $k$'th time at step $j$
    and let $V_j(k)$ be the expected number of points visited exactly $k$ times by the $j$-step walk.
    We clearly have $f_j(0\mid k)\le V_j(k),$ and the formula (4.62) of \cite{W} gives us
    \begin{equation}\label{sqrt}
    \sum_{j=1}^m f_j(0\mid k)\le \sum_{j=1}^m V_j(k)=O\left(m\left(1-\sigma\sqrt\frac2m\right)^k\right)
    \end{equation}
    If $k\ge m^{2/3}$, then $k\cdot \sigma\sqrt\frac2m\ge m^{1/6}\cdot\sigma\sqrt{2}$, and so the right-hand side of (\ref{sqrt}) is $o(\exp(-c m^{1/6}))$ for any positive $c<\sigma\sqrt{2}$.
    Therefore
    $$q_m\le \sum_{k\ge m^{2/3}}^{m} \sum_{j=1}^m f_j(0\mid k)= m\cdot o(\exp(-c m^{1/6}))=o(\exp(-\nu m^{1/6}))$$
    for every positive $\nu< c.$ \endproof

    We denote by $M_n$ the set of all words $w$ (reduced or not) of length $n$ and by $L_n$ the subset
    containing the words with $s\ge n^{2/3}$ prefixes $v$ ending with $y^{\pm 1}$ and satisfying $\sigma_y(v)=0$ .

    \begin{lemma}\label{Qn} We have ${\bf P}(w\in L_n\mid w\in M_n)
    = o (\exp(-\nu n^{1/9}))$ for some positive $\nu$.
    \end{lemma}

     \proof Consider the simple $2$-dimensional random walk associated with the sets $M_n$, where every variable ${\bf Y}_i$
     takes each of the vectors $\pm e_1$, $\pm e_2$ with probability $1/4$, ${\bf Z}_0=0 $ and ${\bf Z}_i = {\bf Z}_{i-1}+{\bf Y}_i$ for $i\ge 1$. Also consider the series ${\bf z}_1,\dots, {\bf z}_m$ (where $m=m({\bf Z}_n))$
     obtained by observing the random walk
$\{{\bf Z}_n\}$ only at the times when ${\bf Y}_n = \pm e_2$ and taking the projection
$({\bf Z}_n)_y$ of ${\bf Z}_n$ on the $y$-axis. Clearly, this  is a simple $1$-dimensional random walk.
     As in Lemma \ref{UU}, we now obtain:

    $${\bf P}(w\in L_n\mid w\in M_n)
    %$$={\bf P}(there \;are\;\ge n^{2/3}\; zeros\; in \; the\; sequence \;({\bf Z}_0)_y,\dots, ({\bf Z}_n)_y)$$
     ={\bf P}(there \;are\; \ge n^{2/3}\; zeros \; among \;{\bf z}_1,\dots, {\bf z}_m)$$
    $$ \le \sum_{j \ge n^{2/3}} {\bf P}(there \;are\; \ge j^{2/3} \; zeros \; among\; {\bf z}_1,\dots, {\bf z}_m, m({\bf Z}_n)=j)=$$
    $$=\sum_{j\ge n^{2/3}} {\bf P}(there \;are\; \ge j^{2/3} \; zeros \; among \; {\bf z}_1,\dots, {\bf z}_m \mid m({\bf Z}_n)=j)\times$$ $$\times {\bf P}(m({\bf Z}_n)=j)=$$
    \begin{equation}\label{notlast}
    =\sum_{j\ge n^{2/3}} {\bf P}(there \;are\; \ge j^{2/3} \; zeros \; among \; {\bf z}_1,\dots, {\bf z}_j)\cdot{\bf P}(m({\bf Z}_n)=j)
    \end{equation}
    By Lemma \ref{Q},
$${\bf P}(there \;are\; \ge j^{2/3} \; zeros \; among \; {\bf z}_1,\dots, {\bf z}_j)=$$ $$ o (\exp(-\nu j^{1/6}))= o (\exp(-\nu n^{1/9}))$$
     if $j\ge n^{2/3}$. The statement of the lemma lemma follows from this inequality and from (\ref{notlast})
    because $\sum_{j\ge n^{2/3}} {\bf P}(m({\bf Z}_n)=j)\le 1$.
    \endproof

    To return to reduced words, we need a simple lemma:

    \begin{lemma} \label{uni} Let $0\le k\le n$ and $w$ be a reduced word of length $k$ in an alphabet $\{x_1^{\pm 1},\dots,x_{\ell}^{\pm 1}\}.$ Then the number of words of length $n$  with reduced form equal to $w$
    does not depend on $w$ (but depends on $k, n,$ and $\ell$).
    \end{lemma}

    \proof Let $\Gamma$ be the Cayley graph of $F_{\ell}$. Then there is a natural bijection between
    the words (or reduced words) and the paths (respectively, reduced paths) starting at $1$.
    This bijection preserves the length. A path $p$ is the reduced form of a path $q$ if and only if
    the label of $p$ is the reduced form of the label of $q$. Observe that any two reduced paths of the
    same length starting at $1$ belong to the same orbit under $Aut(\Gamma)$, where $\Gamma$ is regarded
    as an {\it unlabeled} rooted tree. This implies the assertion of the lemma.
    \endproof

    For a constant $b>0$, the subset $R_n^b$ consists, by definition, of the reduced words $w$ of length $n$ such that
    $w$ has $s\ge bn^{2/3}$ prefixes $v$, where $v$ ends with $y^{\pm 1}$ and $\sigma_y(v)=0.$

    \begin{lemma} \label{bn} There exist  $b,\nu>0$ such that ${\bf P}(w\in R_n^b\mid w\in W_n)
    = o(\exp(-\nu n^{1/9}))$.
    \end{lemma}
    \proof We denote by $M_{n,k}$ the subset of $M_n$, where every word has the reduced form of length $k$,
    and we set $L_{n,k}=M_{n,k}\cap L_n$.  Also we will use $R_{k,n}$ for the subset of $W_k$ whose
    words  have $s\ge n^{2/3}$ prefixes $v$ which (1) end with $y^{\pm 1}$ and (2) satisfy the equation $\sigma_y(v)=0$.

    Note that if a reduced form $w'$ of a word $w\in M_n$ belongs to $R_{k,n}$, then $w\in L_{n,k}$
    because every prefix of $w'$ with the properties (1) and (2) is a reduced form of a prefix of $w$ with the same
    properties.
    Since by Lemma \ref{uni}, every word $w'$ from $R_{k,n}$ is a reduced form of $a(n,k)$ words from
    $M_{n,k}$, where the coefficient $a(n,k)$ does not depend on $w'$, it follows that
    \begin{equation}\label{RL}
    \frac{\# R_{k,n}}{\# W_k} = \frac{a(n,k)(\# R_{k,n})}{a(n,k)(\# W_k)}\le
    \frac{\#L_{n,k}}{a(n,k)(\# W_k)}=\frac{\#L_{n,k}}{\#M_{n,k}}
    \end{equation}
    Also we observe that
    \begin{equation}\label{RR}
    \frac{\# R_{k,n}}{\# W_k} \le \frac{\# R_{k+1,n}}{\# W_{k+1}}
    \end{equation}
    if $1\le k<n$. Indeed, $\# W_{k+1}=3(\# W_k)$ while $\# R_{k+1,n}\ge 3(\# R_{k,n})$ since
    one has $3$ possibilities to add a letter to a word from $R_{k,n}$ and obtain a word from $R_{k+1,n}$.

    There is a constant $\beta\in (0,1)$ such that for every $n\ge 0$, at least a half of the words from $M_n$
    have reduced form of length $>\beta n$ (see \cite{Sa}, Theorem 1), that is, $\sum_{k=m}^{n}\#M_{n,k}\ge \frac12 \sum_{k=0}^n\#M_{n,k}=
    M_n/2$, where $m=\lfloor\beta n\rfloor$. This inequality together with (\ref{RR}) and (\ref{RL}) provide us with

     $$\frac{\#R_{m,n}}{\#W_m}=\min_{m\le k\le n}
     \frac{\#R_{k,n}}{\#W_k }\le \min_{m\le k\le n}
     \frac{\#L_{n,k}}{\#M_{n,k}}\le $$
     \begin{equation}\label{skill}
     \le \frac{\sum_{m\le k\le n}
     \#L_{n,k}}{\sum_{m\le k\le n}\#M_{n,k}}\le 2\cdot\frac{\sum_{k=m}^n
     \#L_{n,k}}{\#M_n}\le 2\cdot\frac{\#L_n}{\#M_n},
     \end{equation}
    whence $\frac{\#R_{m,n}}{\#W_m}=o (\exp(-\nu n^{1/9}))$
     by Lemma \ref{Qn}. But it follows from the definitions that $R_m^b\subset R_{m,n}$
     if $b=2\beta^{-1}$ and $n\ge \beta^{-1}$. Therefore the inequalities (\ref{skill}) prove the lemma (up to the substitution of $m$ for $n$)
     because  $m\le n$ and so $ o (\exp(-\nu n^{1/9}))$ is also $o (\exp(-\nu m^{1/9})).$
    \endproof

\medskip

    From now we fix  $b$ and $\nu$ given in Lemma \ref{bn}, and so let $R_n=R_n^b$. It remains to estimate the size of  $V_n=(W_n\cap N)\backslash R_n$.

    \begin{lemma}\label{number} For arbitrary $\varepsilon>0$, the number of words in $V_n$ is $O(2^{(1+\varepsilon)n})$.
    \end{lemma}
    \proof If $v\in V_n\subset N$, then by Lemma \ref{lone} (b), $v$  has the reduced
    factorization (\ref{red}). Since the subwords $u_1,\dots,u_{s-1}$ start and end with $y^{\pm 1}$ and belong to $H$
    (i.e. $\sigma_y(u_i)=0$ for $i=1,\dots, s-1$), we have at least $s-1$ prefixes $v_i$ of $v$ such that $v_i$ ends with $y^{\pm 1}$ and $\sigma_y(v_i)=0$. It follows that $s<bn^{2/3}+1$ because $v\notin R_n$.

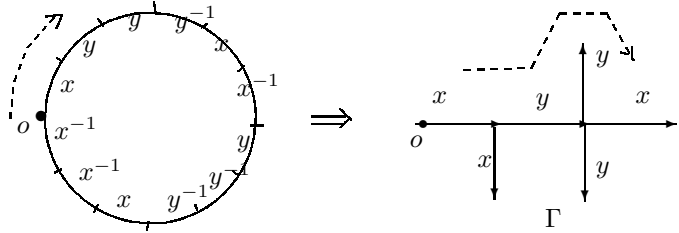
\begin{figure}[h!]
\begin{center}
    %TeXCAD (http://texcad.sf.net/) Picture. File: [subsub7.pic]. Options on following lines.
%\grade{\on}
%\emlines{\off}
%\epic{\off}
%\beziermacro{\on}
%\reduce{\on}
%\snapping{\off}
%\pvinsert{% Your \input, \def, etc. here}
%\quality{8.000}
%\graddiff{0.005}
%\snapasp{1}
%\zoom{4.0000}
\unitlength 1mm % = 2.845pt
\linethickness{0.4pt}
\ifx\plotpoint\undefined\newsavebox{\plotpoint}\fi % GNUPLOT compatibility
\begin{picture}(98.5,35)(5,65)
%\circle(28.5,83.5){28.004}
\put(42.502,83.5){\line(0,1){.7092}}
\put(42.484,84.209){\line(0,1){.7074}}
\put(42.43,84.917){\line(0,1){.7037}}
\multiput(42.341,85.62)(-.03129,.17457){4}{\line(0,1){.17457}}
\multiput(42.216,86.319)(-.032071,.138209){5}{\line(0,1){.138209}}
\multiput(42.055,87.01)(-.032525,.113672){6}{\line(0,1){.113672}}
\multiput(41.86,87.692)(-.032777,.095896){7}{\line(0,1){.095896}}
\multiput(41.631,88.363)(-.032893,.082349){8}{\line(0,1){.082349}}
\multiput(41.368,89.022)(-.032908,.071624){9}{\line(0,1){.071624}}
\multiput(41.071,89.666)(-.032844,.062879){10}{\line(0,1){.062879}}
\multiput(40.743,90.295)(-.032715,.055577){11}{\line(0,1){.055577}}
\multiput(40.383,90.906)(-.032531,.049362){12}{\line(0,1){.049362}}
\multiput(39.993,91.499)(-.0322974,.0439853){13}{\line(0,1){.0439853}}
\multiput(39.573,92.071)(-.0320206,.0392721){14}{\line(0,1){.0392721}}
\multiput(39.125,92.62)(-.031704,.0350933){15}{\line(0,1){.0350933}}
\multiput(38.649,93.147)(-.0334407,.0334425){15}{\line(0,1){.0334425}}
\multiput(38.147,93.648)(-.0350915,.0317059){15}{\line(-1,0){.0350915}}
\multiput(37.621,94.124)(-.0392703,.0320228){14}{\line(-1,0){.0392703}}
\multiput(37.071,94.572)(-.0439835,.0322999){13}{\line(-1,0){.0439835}}
\multiput(36.499,94.992)(-.04936,.032533){12}{\line(-1,0){.04936}}
\multiput(35.907,95.383)(-.055576,.032718){11}{\line(-1,0){.055576}}
\multiput(35.296,95.743)(-.062877,.032848){10}{\line(-1,0){.062877}}
\multiput(34.667,96.071)(-.071623,.032912){9}{\line(-1,0){.071623}}
\multiput(34.022,96.367)(-.082347,.032898){8}{\line(-1,0){.082347}}
\multiput(33.364,96.63)(-.095895,.032783){7}{\line(-1,0){.095895}}
\multiput(32.692,96.86)(-.113671,.032531){6}{\line(-1,0){.113671}}
\multiput(32.01,97.055)(-.138207,.032079){5}{\line(-1,0){.138207}}
\multiput(31.319,97.215)(-.17457,.0313){4}{\line(-1,0){.17457}}
\put(30.621,97.341){\line(-1,0){.7037}}
\put(29.917,97.43){\line(-1,0){.7074}}
\put(29.21,97.484){\line(-1,0){1.4184}}
\put(27.792,97.484){\line(-1,0){.7074}}
\put(27.084,97.43){\line(-1,0){.7037}}
\multiput(26.381,97.341)(-.17457,-.03128){4}{\line(-1,0){.17457}}
\multiput(25.682,97.216)(-.13821,-.032063){5}{\line(-1,0){.13821}}
\multiput(24.991,97.055)(-.113674,-.032518){6}{\line(-1,0){.113674}}
\multiput(24.309,96.86)(-.095898,-.032772){7}{\line(-1,0){.095898}}
\multiput(23.638,96.631)(-.082351,-.032888){8}{\line(-1,0){.082351}}
\multiput(22.979,96.368)(-.071626,-.032904){9}{\line(-1,0){.071626}}
\multiput(22.334,96.072)(-.062881,-.032841){10}{\line(-1,0){.062881}}
\multiput(21.706,95.743)(-.055579,-.032712){11}{\line(-1,0){.055579}}
\multiput(21.094,95.383)(-.049364,-.032528){12}{\line(-1,0){.049364}}
\multiput(20.502,94.993)(-.0439871,-.0322949){13}{\line(-1,0){.0439871}}
\multiput(19.93,94.573)(-.0392739,-.0320184){14}{\line(-1,0){.0392739}}
\multiput(19.38,94.125)(-.0350951,-.031702){15}{\line(-1,0){.0350951}}
\multiput(18.854,93.65)(-.0334444,-.0334388){15}{\line(-1,0){.0334444}}
\multiput(18.352,93.148)(-.0317079,-.0350897){15}{\line(0,-1){.0350897}}
\multiput(17.876,92.622)(-.032025,-.0392685){14}{\line(0,-1){.0392685}}
\multiput(17.428,92.072)(-.0323024,-.0439816){13}{\line(0,-1){.0439816}}
\multiput(17.008,91.5)(-.032536,-.049358){12}{\line(0,-1){.049358}}
\multiput(16.618,90.908)(-.032721,-.055574){11}{\line(0,-1){.055574}}
\multiput(16.258,90.296)(-.032851,-.062876){10}{\line(0,-1){.062876}}
\multiput(15.929,89.668)(-.032916,-.071621){9}{\line(0,-1){.071621}}
\multiput(15.633,89.023)(-.032902,-.082345){8}{\line(0,-1){.082345}}
\multiput(15.37,88.364)(-.032788,-.095893){7}{\line(0,-1){.095893}}
\multiput(15.14,87.693)(-.032537,-.113669){6}{\line(0,-1){.113669}}
\multiput(14.945,87.011)(-.032086,-.138205){5}{\line(0,-1){.138205}}
\multiput(14.785,86.32)(-.03131,-.17457){4}{\line(0,-1){.17457}}
\put(14.659,85.622){\line(0,-1){.7037}}
\put(14.57,84.918){\line(0,-1){.7073}}
\put(14.516,84.211){\line(0,-1){2.1257}}
\put(14.569,82.085){\line(0,-1){.7037}}
\multiput(14.659,81.381)(.03127,-.17457){4}{\line(0,-1){.17457}}
\multiput(14.784,80.683)(.032055,-.138212){5}{\line(0,-1){.138212}}
\multiput(14.944,79.992)(.032512,-.113676){6}{\line(0,-1){.113676}}
\multiput(15.139,79.31)(.032766,-.0959){7}{\line(0,-1){.0959}}
\multiput(15.369,78.639)(.032884,-.082353){8}{\line(0,-1){.082353}}
\multiput(15.632,77.98)(.0329,-.071628){9}{\line(0,-1){.071628}}
\multiput(15.928,77.335)(.032837,-.062883){10}{\line(0,-1){.062883}}
\multiput(16.256,76.706)(.032709,-.055581){11}{\line(0,-1){.055581}}
\multiput(16.616,76.095)(.032525,-.049365){12}{\line(0,-1){.049365}}
\multiput(17.006,75.503)(.0322925,-.043989){13}{\line(0,-1){.043989}}
\multiput(17.426,74.931)(.0320162,-.0392758){14}{\line(0,-1){.0392758}}
\multiput(17.874,74.381)(.0317,-.0350969){15}{\line(0,-1){.0350969}}
\multiput(18.35,73.854)(.0334369,-.0334463){15}{\line(0,-1){.0334463}}
\multiput(18.851,73.353)(.0350879,-.0317099){15}{\line(1,0){.0350879}}
\multiput(19.378,72.877)(.0392667,-.0320273){14}{\line(1,0){.0392667}}
\multiput(19.928,72.429)(.0439798,-.0323049){13}{\line(1,0){.0439798}}
\multiput(20.499,72.009)(.049356,-.032539){12}{\line(1,0){.049356}}
\multiput(21.092,71.618)(.055572,-.032724){11}{\line(1,0){.055572}}
\multiput(21.703,71.258)(.062874,-.032855){10}{\line(1,0){.062874}}
\multiput(22.332,70.93)(.071619,-.03292){9}{\line(1,0){.071619}}
\multiput(22.976,70.633)(.082344,-.032907){8}{\line(1,0){.082344}}
\multiput(23.635,70.37)(.095891,-.032793){7}{\line(1,0){.095891}}
\multiput(24.306,70.141)(.113667,-.032544){6}{\line(1,0){.113667}}
\multiput(24.988,69.945)(.138203,-.032094){5}{\line(1,0){.138203}}
\multiput(25.679,69.785)(.17456,-.03132){4}{\line(1,0){.17456}}
\put(26.377,69.66){\line(1,0){.7037}}
\put(27.081,69.57){\line(1,0){.7073}}
\put(27.788,69.516){\line(1,0){1.4184}}
\put(29.207,69.516){\line(1,0){.7074}}
\put(29.914,69.569){\line(1,0){.7037}}
\multiput(30.618,69.659)(.17457,.03126){4}{\line(1,0){.17457}}
\multiput(31.316,69.784)(.138214,.032047){5}{\line(1,0){.138214}}
\multiput(32.007,69.944)(.113678,.032505){6}{\line(1,0){.113678}}
\multiput(32.689,70.139)(.095902,.032761){7}{\line(1,0){.095902}}
\multiput(33.361,70.368)(.082355,.032879){8}{\line(1,0){.082355}}
\multiput(34.019,70.632)(.07163,.032896){9}{\line(1,0){.07163}}
\multiput(34.664,70.928)(.062885,.032833){10}{\line(1,0){.062885}}
\multiput(35.293,71.256)(.055583,.032706){11}{\line(1,0){.055583}}
\multiput(35.904,71.616)(.049367,.032522){12}{\line(1,0){.049367}}
\multiput(36.497,72.006)(.0439908,.03229){13}{\line(1,0){.0439908}}
\multiput(37.069,72.426)(.0392776,.0320139){14}{\line(1,0){.0392776}}
\multiput(37.619,72.874)(.0350987,.031698){15}{\line(1,0){.0350987}}
\multiput(38.145,73.349)(.0334482,.033435){15}{\line(1,0){.0334482}}
\multiput(38.647,73.851)(.0317119,.0350861){15}{\line(0,1){.0350861}}
\multiput(39.122,74.377)(.0320295,.0392649){14}{\line(0,1){.0392649}}
\multiput(39.571,74.927)(.0323074,.043978){13}{\line(0,1){.043978}}
\multiput(39.991,75.499)(.032542,.049354){12}{\line(0,1){.049354}}
\multiput(40.381,76.091)(.032728,.05557){11}{\line(0,1){.05557}}
\multiput(40.741,76.702)(.032858,.062872){10}{\line(0,1){.062872}}
\multiput(41.07,77.331)(.032924,.071617){9}{\line(0,1){.071617}}
\multiput(41.366,77.975)(.032912,.082342){8}{\line(0,1){.082342}}
\multiput(41.63,78.634)(.032799,.095889){7}{\line(0,1){.095889}}
\multiput(41.859,79.305)(.03255,.113665){6}{\line(0,1){.113665}}
\multiput(42.054,79.987)(.032102,.138201){5}{\line(0,1){.138201}}
\multiput(42.215,80.678)(.03133,.17456){4}{\line(0,1){.17456}}
\put(42.34,81.377){\line(0,1){.7037}}
\put(42.43,82.08){\line(0,1){.7073}}
\put(42.484,82.788){\line(0,1){.7123}}
%\end
\put(63.75,82.75){\vector(1,0){11.5}}
\put(75.25,82.75){\vector(1,0){11.75}}
\put(87,82.75){\vector(1,0){11.5}}
\put(74.25,82.25){\vector(0,-1){9.5}}
\put(86.25,82.75){\vector(0,-1){10.25}}
\put(86,82.5){\vector(0,1){11}}
\put(14,83.75){\circle*{1.414}}
\put(41.75,82.5){\line(1,0){1.75}}
%\emline(29,99)(29.25,97.25)
\multiput(29,99)(.03125,-.21875){8}{\line(0,-1){.21875}}
%\end
%\emline(28.25,69.75)(28,68.5)
\multiput(28.25,69.75)(-.03125,-.15625){8}{\line(0,-1){.15625}}
%\end
%\emline(16,91.5)(16.5,91.25)
\multiput(16,91.5)(.0625,-.03125){8}{\line(1,0){.0625}}
%\end
%\emline(21.5,96.75)(22.25,95.75)
\multiput(21.5,96.75)(.0326087,-.0434783){23}{\line(0,-1){.0434783}}
%\end
%\emline(36,96)(35.25,95.5)
\multiput(36,96)(-.05,-.0333333){15}{\line(-1,0){.05}}
%\end
%\emline(40.25,89.75)(41,90.5)
\multiput(40.25,89.75)(.0326087,.0326087){23}{\line(0,1){.0326087}}
%\end
%\emline(39.5,76.5)(40,76)
\multiput(39.5,76.5)(.0333333,-.0333333){15}{\line(0,-1){.0333333}}
%\end
%\emline(34.5,72)(34.75,71.25)
\multiput(34.5,72)(.03125,-.09375){8}{\line(0,-1){.09375}}
%\end
%\emline(21.5,72)(21,71)
\multiput(21.5,72)(-.0333333,-.0666667){15}{\line(0,-1){.0666667}}
%\end
%\emline(15.75,76.25)(16.75,76.75)
\multiput(15.75,76.25)(.0666667,.0333333){15}{\line(1,0){.0666667}}
%\end
\put(66,85.5){$x$}
\put(79.75,85.25){$y$}
\put(87.75,91){$y$}
\put(93,85.5){$x$}
\put(87.75,76.25){$y$}
\put(72,77){$x$}
\put(16.5,87){$x$}
\put(19.5,92.25){$y$}
\put(25.5,95.75){$y$}
\put(31.75,95.5){$y^{-1}$}
\put(37,92.5){$x$}
\put(40,86.5){$x^{-1}$}
\put(40,79.75){$y$}
\put(36.25,74.75){$y^{-1}$}
\put(30.75,71.75){$y^{-1}$}
\put(24,71.75){$x$}
\put(19,75.25){$x^{-1}$}
\put(15.75,80.75){$x^{-1}$}
\put(50,83.75){\line(1,0){4.5}}
\put(50,82.75){\line(1,0){4.25}}
%\emline(53.75,84.75)(55.25,83.25)
\multiput(53.75,84.75)(.03333333,-.03333333){45}{\line(0,-1){.03333333}}
%\end
%\emline(53.75,82.25)(55,83.5)
\multiput(53.75,82.25)(.03289474,.03289474){38}{\line(0,1){.03289474}}
%\end
%\emline(16.75,97.25)(14.5,97)
\multiput(16.75,97.25)(-.28125,-.03125){8}{\line(-1,0){.28125}}
%\end
%\emline(16.75,97.25)(16,95.5)
\multiput(16.75,97.25)(-.0326087,-.076087){23}{\line(0,-1){.076087}}
%\end
\put(64.75,82.75){\circle*{1.118}}
\put(92.5,91.5){\line(0,1){0}}
%\emline(92.5,91.5)(92.75,91.25)
\multiput(92.5,91.5)(.03125,-.03125){8}{\line(0,-1){.03125}}
%\end
\put(92.5,93){\line(0,-1){1.5}}
%\emline(91,92.25)(92,92)
\multiput(91,92.25)(.125,-.03125){8}{\line(1,0){.125}}
%\end
%\dashline{1}(10,83.5)(10.5,88.5)
\put(9.93,83.43){\line(0,1){.8333}}
\put(10.096,85.096){\line(0,1){.8333}}
\put(10.263,86.763){\line(0,1){.8333}}
%\end
%\dashline{1}(10.5,88.5)(12.25,93)
\multiput(10.43,88.43)(.032407,.083333){9}{\line(0,1){.083333}}
\multiput(11.013,89.93)(.032407,.083333){9}{\line(0,1){.083333}}
\multiput(11.596,91.43)(.032407,.083333){9}{\line(0,1){.083333}}
%\end
%\dashline{1}(12.25,93)(16.25,97.5)
\multiput(12.18,92.93)(.0333333,.0375){15}{\line(0,1){.0375}}
\multiput(13.18,94.055)(.0333333,.0375){15}{\line(0,1){.0375}}
\multiput(14.18,95.18)(.0333333,.0375){15}{\line(0,1){.0375}}
\multiput(15.18,96.305)(.0333333,.0375){15}{\line(0,1){.0375}}
%\end
%\dashline{1}(70.25,90)(79.25,90.25)
\put(70.18,89.93){\line(1,0){.9}}
\put(71.98,89.98){\line(1,0){.9}}
\put(73.78,90.03){\line(1,0){.9}}
\put(75.58,90.08){\line(1,0){.9}}
\put(77.38,90.13){\line(1,0){.9}}
%\end
%\dashline{1}(79.25,90.25)(83,97.5)
\multiput(79.18,90.18)(.0320513,.0619658){13}{\line(0,1){.0619658}}
\multiput(80.013,91.791)(.0320513,.0619658){13}{\line(0,1){.0619658}}
\multiput(80.846,93.402)(.0320513,.0619658){13}{\line(0,1){.0619658}}
\multiput(81.68,95.013)(.0320513,.0619658){13}{\line(0,1){.0619658}}
\multiput(82.513,96.624)(.0320513,.0619658){13}{\line(0,1){.0619658}}
%\end
%\dashline{1}(83,97.5)(89.25,97.5)
\put(82.93,97.43){\line(1,0){.8929}}
\put(84.715,97.43){\line(1,0){.8929}}
\put(86.501,97.43){\line(1,0){.8929}}
\put(88.287,97.43){\line(1,0){.8929}}
%\end
%\dashline{1}(89.25,97.5)(92.25,92.25)
\multiput(89.18,97.43)(.032967,-.0576923){13}{\line(0,-1){.0576923}}
\multiput(90.037,95.93)(.032967,-.0576923){13}{\line(0,-1){.0576923}}
\multiput(90.894,94.43)(.032967,-.0576923){13}{\line(0,-1){.0576923}}
\multiput(91.751,92.93)(.032967,-.0576923){13}{\line(0,-1){.0576923}}
%\end
%\dashline{1}(92.25,92.25)(92.25,92.25)
\put(92.18,92.18){\line(0,1){0}}
%\end
\put(81,69){$\Gamma$}
\put(63,79.75){$o$}
\put(10.75,81.5){$o$}
\end{picture}
\caption{The cancellation tree $\Gamma$ for the word $xyyy^{-1}xx^{-1}yy^{-1}y^{-1}xx^{-1}x^{-1}$ }\label{canc}
\end{center}
\end{figure}

    At first we investigate the word $u\equiv u_0\dots u_s$ . It is
    equal to $1$ in $F$ by Lemma \ref{lone} (b), and therefore there is an oriented finite, planar tree $\Gamma$
    such that (1) every edge $e$ of $\Gamma$ is labeled by a letter from $\{x^{\pm 1}, y^{\pm 1}\}$,
    (2) $\Lab(e^{-1})=\Lab(e)^{-1}$, and (3) starting with a distinguished vertex $o,$ and
    going around  $\Gamma$ we read the word $u.$  (See fig. \ref{canc}: The word $u$
     is written on a circle with $12$ edges; then after $6$ outer pinches one obtains the tree
     $\Gamma$ with $6$ edges.)  The number of (pairs of mutual inverse) edges of
    $\Gamma$ is $|u|/2\le n/2.$

    Note that every word  $u_i$ is reduced, and so $u$ has at most $s$ pairs of {\it neighbor}
    mutually inverse letters. It follows that the graph $\Gamma$ has at most $s+1$
    vertices of valency $1$ (taking $o$ into account too), i.e., at most $s+1$ leaves. Every
    unlabeled tree $\Gamma$ with $t+1\le s+1$ leaves
    can be constructed as follows. We take a tree $\Gamma'$ with
    $t$ leaves, chose an  edge $e$ of $\Gamma',$ and attach a new leave $p$
    to $e_+$ (i.e., we insert $pp^{-1}$ in the boundary path
     of the tree after $e$). We therefore have
    the following estimate on the numbers $N_{m,t+1}$ of unlabeled trees with a distinguished vertex
    (up to isomorphism), with  $\le m\le n/2$ positive edges, and with $t+1$  leaves:
    $N_{m,t+1}\le (2m)m N_{m,t}$, whence by induction on $t$ we obtain
    $N_{m,t+1}\le (2m^2)^{t+1}.$ Summing over $m\le n/2$ and taking into account that
    $t\le s< bn^{2/3}+1$ we obtain the upper bound of the form
    \begin{equation}\label{O}
    O(n^{2bn^{2/3}})=O(2^{2bn^{2/3}\log_2n})=O(2^{\varepsilon n/2})
    \end{equation}
    for the number of the unlabeled trees under consideration, where $\varepsilon$ is arbitrary fixed positive number.
    Then for each of the unlabeled edges $e$ of $\Gamma$,
    we choose the direction and one of two possible labels: $x$ or $y$. This gives at most
    $4^{n/2}$ labelings of the tree. Taking into account the estimate (\ref{O}), we obtain
    $O(2^{(1+\varepsilon/2)n})$ possible words $u\equiv u_1\dots u_s$ if $s< bn^{2/3}+1.$

    Now to obtain the word $v\in V_n$, we consequently insert the powers of $x$ in $u$. We start with the choice
    of the the first cut to insert $x^{k_1}$ with $|k_1|\le n$. Clearly, we have at most $n(2n)$ options.
    Then we similarly insert $x^{k_2},...,  x^{k_s}$. so the total number of options $\le (2n^2)^s$
    which again is $O(2^{\varepsilon n/2})$ since $s< bn^{2/3}+1.$  Together with the estimate of the previous paragraph, this gives  $O(2^{(1+\varepsilon)n})$ as the upper bound for the cardinality of $V_n$.
    \endproof

    {\bf Proof of Theorem \ref{cogr} (2).}
Since the $\ell$-subnormal closure of $x$ is contained in the $3$-subnormal closure of $x$, we may
     assume that $\ell=3$.
     To apply Lemma \ref{maxgr}, we rename: $H_1=N$, and $N_1$ to be the normal closure of $x$ in $H_1$,
     i.e., $N_1$ is the $3$-subnormal closure of $x$ in $F$. It follows from lemmas \ref{bn} and \ref{number}
      that the number of elements of length $n$ in $H_1$ with respect to the generators of $F$
      can be expressed as $a_n 3^n$ with $a_n\le C(\exp(-\nu n^{1/9}))$ for some $C, \nu>0.$
      Since
      $$\sum_{n=1}^{\infty}\exp(-\nu n^{1/9})<\int_0^{\infty}\exp(-\nu x^{1/9})dx= 9!\nu^{-9},$$
      %Therefore the series
      %$\sum_{n=1}^{\infty}a_n$ converges, and
      the cogrowth of $N_1$ in $F$ is maximal by Lemma \ref{maxgr}
      applied to the pair $(H_1, N_1)$.
      $\Box$

      \bigskip

{\bf Acknowledgements.} Theorem \ref{cogr} answers a question asked by Yuri Bahturin.
  I am thankful to him for the lasting discussions. I also thank Mark Ellingham and
  Rostislav Grigorchuk for the information relevant to Lemma \ref{graph} and Theorem
  \ref{rate} and thank Victor Guba and Vladimir Piterbarg for conversations.

\bigskip

{\bf Alexander A. Olshanskii:} Department of Mathematics, Vanderbilt University, Nashville
37240, USA, and Moscow State University, Moscow 119991, Russia.

{\it E-mail}: alexander.olshanskiy@vanderbilt.edu

\end{document}